\input amstex
\documentstyle{amsppt}
\magnification1200
\tolerance=500
\overfullrule=0pt
\def\n#1{\Bbb #1}
\def\p{\Bbb C_{\infty}}

\def\im{\hbox{im }}

\def\Gal{\hbox{Gal }}
\def\Aut{\hbox{Aut }}

\def\Hom{\hbox{Hom}}

\def\Ker{\hbox{Ker }}

\def\Aut{\hbox{Aut}}

\def\det{\hbox{det }}

\def\Id{\hbox{Id}}

\def\Proj{\hbox{Proj }}
\def\Sing{\hbox{Sing}}

\def\rank{\hbox{ rank }}

\def\s{\sigma}
\def\si{\sigma}
\def\z{\zeta}

\def\e11{E_{11}}

\def\g{\goth }
\def\ga{\goth A}

\def\ve{\varepsilon}
\def\vf{\varphi}
\def\bh{\beth}

\def\de{\delta}

\def\ga{\gamma}
\def\la{\lambda}
\def\La{\Lambda}
\def\Ga{\Gamma}
\def\be{\beta}
\def\z{\zeta}
\def\al{\alpha}
\def\om{\omega}
\def\th{\theta}
\def\vk{\varkappa}
\def\vt{\vartheta}

\topmatter
\title
L-functions of Carlitz modules, resultantal varieties and rooted binary trees, II.
\endtitle
\author
A. Grishkov, D. Logachev, A. Zobnin\footnotemark \footnotetext{E-mails: shuragri{\@}gmail.com; logachev94{\@}gmail.com (corresponding author); al.zobnin{\@}yandex.ru\phantom{*******************}}
\endauthor
\NoRunningHeads
\address
First author: Departamento de Matem\'atica e estatistica
Universidade de S\~ao Paulo. Rua de Mat\~ao 1010, CEP 05508-090, S\~ao Paulo, Brasil, and Omsk State University n.a. F.M.Dostoevskii. Pr. Mira 55-A, Omsk 644077, Russia.
\medskip
Second author: Departamento de Matem\'atica, Universidade Federal do Amazonas, Manaus, Brasil
\medskip
Third author: Faculty of Computer Science, National Research University Higher School of Economics, Moscow, Russia
\endaddress
\thanks Thanks: The first and second authors are grateful to FAPESP, S\~ao Paulo, Brazil for a financial support (process No. 2017/19777-6). The first author is grateful to SNPq, Brazil, to RFBR, Russia, grant 16-01-00577a (Secs. 1-4), and to Russian Science Foundation, project 16-11-10002 (Secs. 5-8) for a financial support. The authors are grateful to Aleksandr Esterov for some important remarks, and to (an) anonymous reviewer(s) who indicated a simple proof of Theorem 8.1, refinements of some statements of Section 0.1.3 concerning the $L$-functions of $t$-motives, as well as of some other remarks.
\endthanks
\abstract We continue study of some algebraic varieties (called resultantal varieties) started in a paper of A. Grishkov, D. Logachev
"Resultantal varieties related to zeroes of L-functions of Carlitz modules". These varieties are related with the Sylvester matrix for the resultant of two polynomials, from one side, and with the L-functions of twisted Carlitz modules, from another side. Surprisingly, these varieties are described in terms of finite weighted rooted binary trees. We give a (conjecturally) complete description of them, we find parametrizations of their irreducible components and their invariants: degrees, multiplicities, Jordan forms, Galois actions. Proof of the fact that this description is really complete is a subject of future research. Maybe a generalization of these results will give us a solution of the problem of boundedness of the analytic rank of twists of Carlitz modules.
\endabstract
\endtopmatter
\document
For a list of notations see p. 17. For a short version of the present paper see [GL21]. 
\medskip
{\bf 0. Introduction.} There is a famous
\medskip
{\bf Open problem 0.0.1.} (a) Are ranks of elliptic curves over $\n Q$ bounded, or not? If they are bounded, what is the upper bound of ranks; what is the maximal number $r$ such that there exists infinitely many elliptic curves over $\n Q$ of rank $r$?
\medskip
(b) The same questions for the ranks of all $\bar \n Q/\n Q$-twists of a fixed elliptic curve $E$ over $\n Q$.
\medskip
Earlier it was thought that the ranks are unbounded. There are examples of curves of rank $\ge 28$ and infinite series of curves of rank $\ge 18$. Nevertheless, a recent paper [PPVW] gives evidence that the ranks of all elliptic curves over $\n Q$ are bounded, moreover, it gives evidence that there exists only finitely many elliptic curves over $\n Q$ of rank $\ge 22$.
\medskip
The purpose of the present paper is to continue the study (started in [GL16]) of an analog of this problem for the function field case. A function field case analog of an elliptic curve is a Drinfeld module of rank 2. They are particular cases of Anderson t-motives\footnotemark \footnotetext{The terminology "t-motive" and "module" reflects two (equivalent) definitions of these objects. Traditionally, the Drinfeld and Carlitz objects are called "modules" while the Anderson objects are called "t-motives": Drinfeld module is an Anderson t-motive, up to some equivalence.}. Really, we shall consider a simpler example of Anderson t-motives --- the Carlitz module $\g C$ ( = the Drinfeld module of rank 1) over a finite field $\n F_q$ (where $q$ is a power of a prime $p$), its tensor powers $\g C^n$ and their twists. The twists of $\g C^n$ are parameterized by polynomials $P\in \n F_q[\th]$ (here $\th$ is an independent variable). We denote the degree of $P$ by $m$,
$$P=\sum_{\iota=0}^m a_\iota\theta^\iota \hbox{ where } a_\iota\in \n F_q\eqno{(0.0.1c)}$$
and we denote the corresponding twist of $\g C^n$ by $\g C^n_P$ (see 0.1.2 for details).
\medskip
We consider the analytic rank of $\g C^n_P$. There are various versions of $L$-functions of an Anderson t-motive $M$, see 0.1.3B for a discussion, and respectively various versions of the notions of the analytic rank and of the statements of the Birch --- Swinnerton-Dyer conjecture. We shall consider in the present paper the $L$-function $L(M,T)$ which is defined, for example, in [L], upper half of page 2603 (warning: $\tau$ of the present paper is $u$ of [L]), or in [Ga], Chapter 4, III, Definition 4.10, or in [G92], 3.2.15, or in [TW], formula (2.1); this $L$-function is called a na\"{\i}ve $L$-function in [B12]. Also, a very simple and explicit formula for $L$-function is given in [GL16], I.1 (caution: this formula of [GL16] is erroneous, therefore it should be treated only as a first approximation to the correct formula. Clearly these errors do not have any influence to the results of [GL16]). In 0.1.3A of the present paper we give a definition of a slightly weaker object, namely a $L$-function up to finitely many Euler factors.
\medskip
For $M=\g C^n_P$ we have: $L(\g C^n_P,T)$ belongs to $(\n F_q[t])[T]$, i.e. it should be considered as a polynomial in $T$ with coefficients in $\n F_q[t]$ where $\n F_q[t]$ is a subring of the Anderson ring $\n F_q(\th)\{t,\tau\}$ (see (0.1.0) below, case $\Phi=\n F_q(\th)$).
\medskip
The analytic rank of $\g C^n_P$ is, by definition, the order of zero of $L(\g C^n_P,T)$ at $T=1$ (recall that for other definitions of $L$-functions we shall get other definitions of the analytic rank, see 0.1.3B). Its study was started in [GL16]. It is denoted by $r_1(n,P)$ or simply by $r_1(P)$ if $n$ is fixed (subscript 1 because of the order of zero at $T=1$).
\medskip
Hence, this analytic rank of t-motives is an analog of the analytic rank of elliptic curves. For example, [L], Proposition 2.1, p. 2604 can be considered as an analog of the strong form of the Birch and Swinnerton-Dyer conjecture for $L(M,T)$ at $T=1$ (here $M$ can be even a generalization of an Anderson t-motive like in [L]): "l'ordre d'annulation de $L(X, (\Cal E,u), T)$ en $T=1$ est \'egal \`a (dimension of a cohomology group) et la valeur sp\'eciale $L^*(X, (\Cal E,u),1)$ est \'egale \`a (an expression) dans $F$".
\medskip
{\bf 0.0.1d.} The results of [GL16] show that the behavior of $r_1(P)$, while $P$ varies, resembles the behavior of ranks of twists of a fixed elliptic curve.
\medskip
First, let us recall this behavior. Let $E$ be an elliptic curve over $\n Q$ defined by the equation $y^2=x^3+g_2x+g_3$, where $g_2, \ g_3\in \n Q$ (to simplify notations, we do not consider its minimal model). For a squarefree $d\in \n Z$ let $E_d$ be its $d$-twist defined by the equation $dy^2=x^3+g_2x+g_3$. The parity of $E_d$ ( = the sign of its functional equation) depends on a residue of $d$ by some module (depending on the conductor of $E$ and on a chosen model), i.e. its behavior is regular. Moreover, the half of residues are even and the half of residues are odd. The parity of the rank of $E_d$ is equal to the the parity of $E_d$. For almost all $d$ the rank of $E_d$ takes the minimal possible value, i.e. it is 0 for even curves and 1 for the odd ones. Rare jumps of the rank of $E_d$ occur, i.e. occasionally the rank of the even $E_d$ can be 2 (rarely), 4 (much more rarely) etc, and respectively the rank of the odd $E_d$ can be 3 (rarely), 5 (much more rarely) etc. There jumps of rank seem to be irregular, i.e. we do not know a formula for $d$ such that the rank of $E_d$ has a jump. 
\medskip
The behavior of $r_1(P)$ is the following. For most $P$ we have $r_1(P)=0$. There exists a coset of a subgroup of the index $q^2$ in the group of twists (see Lemma 0.1.2.3) such that for $P$ belonging to this coset we have $r_1(P)\ge1$ ([GL16], (4.4), (4.5)). Further, rare jumps of $r_1(P)$ occur.
\medskip
An analog of the parity in characteristic $p$ could be the residue modulo $p$. Unfortunately, for the functional field case there is no information either on the parity of $r_1(P)$ or its residue modulo $p$, because for Anderson t-motives there is no functional equation for their $L$-series.
\medskip
Therefore, we see that the behavior of $r_1(P)$, while $P$ varies, is as similar as possible to the behavior of ranks of twists of a fixed elliptic curve, taking into consideration the absence of a functional equation.
\medskip
The function field analog of Open problem 0.0.1 (b) is
\medskip
{\bf Open problem 0.0.2.} Problem of boundedness of ranks of twists of Carlitz modules: Are the analytic ranks of $\g C^n_P$ bounded? ($q$, $n$ are fixed, $m$, $P$ vary). If it is bounded, what is the maximal value of rank; what is the maximal value of rank that occurs for infinitely many $P$?
\medskip
Contents of the present paper should be considered as an approach (second step) to a solution of 0.0.2 (the first step was [GL16]).
\medskip
Let us describe what was made in [GL16]. Theorem 3.3 of [GL16] gives an explicit formula for $L(\g C^n_P,T)$: it is the characteristic polynomial of a matrix $\Cal M(P,n,\bar \kappa)_{fc}$ (roughly speaking; see 0.1.4 for $\bar \kappa$, Theorem 0.1.6 for the exact statement. The subscript $fc$ always means "finite characteristic") whose entries are polynomials in $a_0,\dots, a_m$ from (0.0.1c), and in $t$, with coefficients in $\n F_p$. The matrix $\Cal M(P,n,\bar \kappa)_{fc}$ is a generalization of the Sylvester matrix of the resultant of two polynomials.
\medskip
For the case $q=3$, $n=1$ (this is the simplest non-trivial case, because there are no non-trivial twists for $q=2$, see 0.1.2.7) elementary considerations of the type "dimension of a generic variety is equal to the quantity of variables minus the quantity of equations" show that the maximal value of $r_1(P)$ should be 3. In reality, there are examples of $P$ such that $r_1(P)=6$ ([GL16], Table 6.11). Hence, we do not know even a conjectural answer to 0.0.2.
\medskip
{\bf 0.0.2a.} The coefficients $a_0, \dots, a_m$ of $P$ of (0.0.1c) form a point of $\n A^{m+1}(\n F_q)$ (the affine space of dimension $m+1$). Let us fix $i$ and consider all $P$ of degree $m$ such that $r_1(P)\ge i$. The coefficients of these $P$ form a subset of $\n A^{m+1}(\n F_q)$ denoted by $X_1(q,n,m,i)_{set; \ fc}$ (subscript $set$ indicates that we consider it as a set, without any additional structure).
\medskip
Formula for $L(\g C^n_P,T)$ and the definition of the analytic rank show that there exists an affine algebraic variety denoted by $X_1(q,n,m,i)_{fc}$ which is the set of $\bar \n F_q$-zeroes of some (non-homogeneous) polynomials $D_{**fc}\in\n F_p[a_0,\dots, a_m]$. This variety has the property:

$$X_1(q,n,m,i)_{set; \ fc}=X_1(q,n,m,i)_{fc}(\n F_q).$$

See [GL16], between (6.2) and (6.3) for the definition of $D_{**fc}$ (in [GL16] polynomials $D_{**fc}$ are denoted by $D_{**}$); two subscripts ** mean that $D_{**fc}$ depend on 2 integer parameters. See (0.1.13) for details.
\medskip
Hence, we can state
\medskip
{\bf Open problem 0.0.3.} What are the dimensions of $X_1(q,n,m,i)_{fc}$? Let $q, \ n$ be fixed. Are there numbers $i$ such that for all $m$ the varieties $X_1(q,n,m,i)_{fc}$ are empty sets? (If yes then the answer to 0.0.2 is: "bounded", but not the converse: a non-empty variety can have no $\n F_q$-points).
\medskip
{\bf Remark.} To state a question on dimensions, we have to consider varieties $X_1(q,n,m,i)_{fc}$, because there is no notion of dimension of $X_1(q,n,m,i)_{set; \ fc}$ which is only a subset of $\n A^{m+1}(\n F_q)$.
\medskip
{\bf 0.0.4.} If $D_{**fc}$ were independent then the dimensions of $X_1(q,n,m,i)_{fc}$ would be so small that they would be empty for all sufficiently large $i$. For example, for $q=3$, $n=1$ in this case we would have $X_1(q,n,m,i)_{fc}=\emptyset$ $\forall \ i\ge4$, $\forall \ m$, and the maximal value of $r_1(P)$ would be 3, as it was mentioned above. Nevertheless, calculations of [GL16], Section 6 (especially 6.7 -- 6.10) give evidence that for $q=3$, $n=1$, $i\le3$ the codimensions of $X_1(3,1,m,i)_{fc}$ are less than the quantity of $D_{**fc}$, i.e. apparently $D_{**fc}$ are dependent.
\medskip
It turns out that it is much easier to consider the behavior (as $P$ varies) of $L(\g C^n_P, T)$ not at $T=1$ but at $T=\infty$. The deficiency of the order of pole of $L(\g C^n_P, T)$ at $T=\infty$ (see (0.1.11) for the exact definition) is called the analytic rank of $P$ at infinity, it is denoted by $r_\infty(P)$. Analogously to the case of $r_1(P)$, we denote the set of $(a_0,\dots,a_m)\in \n A^{m+1}(\n F_q)$ such that $r_\infty(P)\ge i$ by $X_\infty(q,n,m,i)_{set; \ fc}$. Again like the case of $r_1(P)$, there exists an algebraic variety $X_\infty(q,n,m,i)_{fc}$ which is the set of $\bar \n F_q$-zeroes of some homogeneous polynomials $H_{**fc}\in\n F_p[a_0,\dots, a_m]$ such that

$$X_\infty(q,n,m,i)_{set; \ fc}=X_\infty(q,n,m,i)_{fc}(\n F_q).$$

It is easy to show that $D_{**fc}$ are linear combinations of $H_{**fc}$ with integer coefficients, so we can expect that study of $X_\infty(q,n,m,i)_{fc}$ can shed light to study of $X_1(q,n,m,i)_{fc}$ and hence to a solution of the Problem 0.0.3. Numerical data of [GL16], 6.16 and 6.11 show that there exists a correlation between high values of $r_1(P)$ and $r_\infty(P)$ (for example, polynomials $P_2, \ P_3, \ P_4$ from [GL16], Table 6.11 having high $r_1(P)$, also have high $r_\infty(P)$ ).
\medskip
Since $H_{**fc}$ are homogeneous, we can consider a projective variety of their zeroes; it is denoted by $X_\infty(q,n,m,i)_{fc}$ as well.
\medskip
Polynomials $H_{**fc}$ are highly dependent ([GL16], Theorem 8.6), i.e. $X_\infty(q,n,m,i)_{fc}$ are of low codimension.
\medskip
{\bf 0.0.6.} Recall that the entries of $\Cal M(P,n,\bar \kappa)_{fc}$ are polynomials in $a_0,\dots, a_m, t$. The coefficients of these polynomials are images of some binomial coefficients (integer numbers) under the residue map $\n Z \to \n F_p$, see (0.1.4a). This means that there are matrices $\Cal M(P,n,\bar \kappa)\in \n Z[a_0,\dots, a_m, t]$ whose entries are given by the same formulas (0.1.4a) but in $\n Z$ (see (0.2.0)), i.e. $\Cal M(P,n,\bar \kappa)$ is a distinguished $\n F_p\to \n Z$-lift of $\Cal M(P,n,\bar \kappa)_{fc}$.
\medskip
Polynomials $H_{**fc}$ are defined using the entries of $\Cal M(P,n,\bar \kappa)_{fc}$. The same formulas applied to the entries of $\Cal M(P,n,\bar \kappa)$ give us polynomials $H_{**}$. We see that $H_{**}$ can be considered as canonical $\n F_p\to \n Z$-lifts of $H_{**fc}$. Hence we can consider the varieties of zeroes of $H_{**}$ in $\n P^m(\n C)$. They  are denoted by $X_\infty(q,n,m,i)$. These varieties are defined over $\n Q$, because coefficients of $H_{**}$ are in $\n Z$. It turns out that for this case of characteristic 0 the simplest possible value of $q$ is 2: although for $q=2$ there are no twists of $\g C^n$ (see (0.1.2.7)), the theory of $X_\infty(2,n,m,i)$ --- unlike the theory of $X_\infty(2,n,m,i)_{set; \ fc}$ over $\n F_2$ --- is non-trivial. Moreover, it turns out that for $q=2$ the varieties $X_\infty(2,n,m,i)$ (conjecturally) do not depend on $n$, see Conjecture 0.2.4, Remark 0.2.5 (1).
\medskip
{\bf 0.0.7.} Further, we can consider formally polynomials $H_{**}$ for $n=0$: although $\g C^0$ does not exist, the polynomials $H_{**}$, and hence varieties $X_\infty(2,0,m,i)$ have meaning for $n=0$, and (again conjecturally) $X_\infty(2,0,m,i)$ coincides with $X_\infty(2,n,m,i)$ for all $n>0$. Finally, the case $n=0$ is especially simple, because for this case the set of $H_{**}$ depends on 1 parameter and not on 2 parameters that simplifies their study. See Section 0.2 for a more detailed discussion.
\medskip
Now we can formulate the subject of the present paper: study of varieties $X_\infty(2,0,m,i)$ which are denoted simply by $X(m,i)$. The above arguments show their importance. We prove some theorems and state some conjectures on their irreducible components, their rational parameterizations, degrees, Jordan forms, Galois actions and multiplicities.
\medskip
It is necessary to emphasize that formally for reading of the present paper it is not necessary to know what is a Carlitz module, $L$-functions etc. --- definitions of $X(m,i)$ are elementary, and methods are purely combinatorial. The reader can start reading from Section 1; the notion of the Carlitz modules is used only in the Introduction, in order to show the importance of the subject.
\medskip
{\bf 0.0.8.} It is interesting that even in this interpretation varieties $X(m,i)$ "remember" their origin. For example, one of the important constructions (Proposition 7.1) depends on Conjecture 0.2.4 which indicates that $X_\infty(2,0,m,i)=X_\infty(2,n,m,i)$ for all $n>0$. Recall that the varieties $X_\infty(2,n,m,i)$ come from the $L$-functions of $\g C^n_P$.
\medskip
The subsequent sections of the Introduction contain a more detailed exposition of the theory of $L$-functions of twisted Carlitz modules (Section 0.1) and more detailed definitions of the above objects (Section 0.2). Further, we give the contents of the present paper (Section 0.3) and the possibilities of further research (Section 0.4).
\medskip
\newpage
{\bf 0.1. Exposition of the theory of $L$-functions of twisted $n$-th tensor powers of Carlitz modules.}
\medskip
Carlitz modules, their tensor powers and twists are particular cases of Anderson t-motives, so we recall their definition. Standard reference for t-motives is [A86], we use its notations. See also [GL20].

Let $q$ be a power of $p$ and $\n F_q$ the finite field of order $q$. The field $\n F_q(\theta)$ is its field of rational functions, it is the function field analog of $\n Q$. The field of the Laurent series $\n F_q((1/\theta))$ is the function field analog of $\n R$. By definition, $\p$ is the completion of the algebraic closure of $\n F_q((1/\theta))$, it is the function field analog of $\n C$.
\medskip
{\bf 0.1.0.} Let $\Phi$ be a subfield of $\p$ (we shall consider only the cases $\Phi=\n F_q(\th)$ and $\Phi=\p$). The Anderson ring $\Phi\{t, \tau\}$ is the ring of
non-commutative polynomials in two variables $t$, $\tau$ over $\Phi$ satisfying the following relations:
$$t\tau = \tau t; \ \forall a\in \Phi \hbox{ we have }\ ta=a t; \ \tau a = a^q \tau$$
We need the following (less general than in [A86]) version of the definition of Anderson t-motives over $ \Phi$ (throughout the whole paper $M$ will mean an Anderson t-motive):
\medskip
{\bf Definition 0.1.1.} An Anderson t-motive $M$ is a $\Phi\{t, \tau\}$-module such that
\medskip
1. $M$ considered as a $\Phi[ t]$-module is free of finite rank $r$;
\medskip
2. $M$ considered as a $\Phi\{\tau\}$-module is free of finite rank $n$;
\medskip
3. $\exists \varkappa>0$ such that $(t-\theta)^\varkappa M/\tau M=0$.
\medskip
An isomorphism of t-motives is an isomorphism of modules.
\medskip
Numbers $r$, resp. $n$ are called the (ordinary) rank of $M$, resp. the dimension of $M$. The ordinary rank of $M$ should not be confused with $r_1$ --- the analytic rank of $M$ considered below.
\medskip
The Carlitz module $\g C$ is an Anderson t-motive over $\p$ having $r=n=1$. Let $\{e\}=\{e_1\}$ be the only element of a basis of $M$ over $\p\{\tau\}$. $\g C$ is given by the relation $te=\theta e +\tau e$. We have: $e$ is also the only element of a basis of $\g C$ over $\p[t]$, and the multiplication by $\tau$ is given by the formula $$\tau e=(t-\th)e\eqno{(0.1.1a)}$$
\medskip
The $n$-th tensor power of $\g C$, denoted by $\g C^n$, is the $n$-th tensor power of $\g C$ over the ring $\p[t]$. Hence, it has the ordinary rank $r=1$. We denote the only element $e\otimes e \otimes e \otimes ... \otimes e$ of a basis of $\g C^n$ over $\p[t]$ by $e_n$. The action of $\tau$ on $e_n$ is given by the formula
$$\tau e_n=(t-\th)^ne_n\eqno{(0.1.1b)}$$
It is easy to check that the $\Phi\{t, \tau\}$-module defined by this formula satisfies 0.1.1.2-3, i.e. it is an Anderson t-motive of the ordinary rank $r=1$ and dimension $n$.
\medskip
{\bf 0.1.2.} The Carlitz module $\g C_\Phi$ over any $\Phi$ (not necessarily $\Phi=\p$) is an Anderson t-motive over $\Phi$ defined by the same formula (0.1.1a), its $n$-th tensor power $\g C^n_\Phi$ is defined by (0.1.1b). A twist of $\g C^n_\Phi$ is an Anderson t-motive $M$ over $ \Phi$ which is isomorphic to $\g C^n$ over $\p$ (twists of Carlitz modules are considered also in [GL16] (any $n$); [Ge] ($n=1$)). Let us recall their description. There is a general theory. Let $X$ be an object defined over a field $K$, and $L$ an extension of $K$. The set of $L/K$-twists of $X$ is\footnotemark \footnotetext{This is a classical fact. Idea of the proof: Let $X'$ be a $L/K$-twist of $X$, i.e. there exists an isomorphism $\vf: X\otimes_KL\to X'\otimes_KL$. Let $g\in \Aut(L/K)$. We have $\vf^{-1}\circ g(\vf)\in \Aut_{L}(X\otimes_KL)$. The map $g\mapsto \vf^{-1}\circ g(\vf)$ is a cocycle defining an element of $H^1(\Aut(L/K), \Aut_{L}(X\otimes_KL))$. We denote it by $c_\vf$. If instead of $\vf$ we choose another isomorphism $\psi: X\otimes_KL\to X'\otimes_KL$ then $\vf^{-1}\circ \psi\in \Aut_{L}(X\otimes_KL)$. We denote $\vf^{-1}\circ \psi$ by $\al$.  The cocycle $c_\psi:= \{g\mapsto \psi^{-1}\circ g(\psi)\}$ satisfies $c_\psi(g)=\al^{-1}\circ c_\vf\circ g(\al)$ hence the image of $c_\vf$ in $H^1(\Aut(L/K), \Aut_{L}(X\otimes_KL))$ does not depend on $\vf$. Details can be found in any textbook of algebraic geometry.} $$H^1(\Aut(L/K), \Aut_{L}(X\otimes_KL))\eqno{(0.1.2.1)}$$

We apply this formula to $X=\g C^n_\Phi$, $L=\p$, $K=\Phi$. We have $\g C^n_\Phi\otimes_{\Phi}\p=\g C^n$.
\medskip
{\bf Lemma 0.1.2.2.} $\Aut(\g C^n)=\n Z/(q-1)$. The action of $\Aut(\p/\Phi)$ on $\Aut(\g C^n)$ is trivial.
\medskip
{\bf Proof.} Let $g\in\Aut(\g C^n)$. Since $\g C^n$ is a free $\p[t]$-module of rank 1 with a $\tau$-action and $e_n$ is its basis, $g$ is defined by its matrix in $GL_1(\p[t])$, i.e. by an element $\om\in\p[t]^*=\p^*$ (the set of invertible elements in $\p[t]$). Further, $g$ preserves the multiplication by $\tau$ ($= \tau$-action): $$\tau\cdot g(e_n)=g(\tau e_n)$$ We have:
$$\tau\cdot g(e_n)=\tau \om e_n=\om^q\tau e_n=\om^q(t-\th)^n e_n; \ \ \ g(\tau e_n)= \om  (t-\th)^n e_n$$
hence $\om^{q-1}=1$, i.e. $\om\in \n F_q^*$ and $\Aut(\g C^n)=\n Z/(q-1)$. Since $\n F_q^*\subset \Phi$ the action of $\Aut(\p/\Phi)$ on $\Aut(\g C^n)$ is trivial. $\square$
\medskip
{\bf Lemma 0.1.2.3.} The set of twists $H^1(\Aut(\p/\Phi), \Aut(\g C^n))$ is $\Phi^*/\Phi^{*(q-1)}$.
\medskip
{\bf Proof.} The action of $\Aut(\p/\Phi)$ on $\Aut(\g C^n)$ is trivial, hence $H^1$ is Hom. The lemma follows from the Kummer theory. $\square$
\medskip
{\bf 0.1.2.4.} Let $P\in\Phi^*$ be a representative of an element of $\Phi^*/\Phi^{*(q-1)}$. We denote the corresponding twist by $\g C^n_P$. It is given by the formula ($e_n$ is as earlier)

$$\tau e_n=P(t-\th)^ne_n\eqno{(0.1.2.5)}$$

This follows immediately from the proof of (0.1.2.1), see the footnote, and the explicit formula for the Kummer isomorphism $$\Hom (\Aut(\p/\Phi), \n Z/(q-1)) \to \Phi^*/\Phi^{*(q-1)}$$ Two twists $\g C^n_{P}$, $\g C^n_{P_1}$ are isomorphic over $\Phi$ iff $P_{1}/P\in \Phi^{*(q-1)}$. Explicitly, it is defined by the following formula.
\medskip
{\bf 0.1.2.5a.} Let $P_1=PQ^{q-1}$ for $Q\in \Phi$. We denote the element $e_n$ from (0.1.2.5) for $P_1$ by $e_{n,1}$. The isomorphism between $\g C^n_{P}$ and $\g C^n_{P_1}$ is given by the formula
$$e_{n,1}=Qe_n\eqno{(0.1.2.6)}$$
\medskip
From now we shall consider the case $\Phi=\n F_q(\theta)$. In this case we can always choose $P$ a polynomial (i.e. $P$ from (0.0.1c)) $(q-1)$-th power free.
\medskip
{\bf (0.1.2.7)} In particular, there are no non-trivial twists for $q=2$.
\medskip
It is possible to define a t-motive $M$ not over a field $\Phi$, but also over a ring. For example, if the above $P,\ P_1\in \n F_q[\th]$ then $\g C^n_{P}$, $\g C^n_{P_1}$ are over $\n F_q[\th]$. Moreover, (0.1.2.6) shows that they are not isomorphic over $\n F_q[\th]$ (except $Q\in \n F_q$, i.e. $Q$ is a constant), i.e. they are different $\n F_q[\th]$-forms of the same t-module over $\n F_q(\th)$.
\medskip
{\bf 0.1.3A. $L$-functions.} As it was mentioned above, we consider in the present paper the $L$-function $L(M,T)$, and we give here (for the reader's convenience) a definition of a weaker object, namely $L_{fmef}(M,T)$ --- a $L$-function up to finitely many Euler factors.
\medskip
For a definition of local $L$-factors in the bad primes (i.e. for a definition of $L(M,T)$ itself) we refer the reader to [Ga], Chapter 4, III, Definition 4.10. We do not need in the present paper the definition of $L(M,T)$, because we shall use exclusively the formula for $L(\g C^n_P,T)$ given in Theorem 0.1.6 ( = [GL16], Theorem 3.3) which, in its turn, comes from the Lefschetz trace formula (see, for example, [L], p.2603, formula (1) for an explicit form).
\medskip
{\bf Remark 0.1.3A.0.} It is necessary to emphasize that $L(M,T)$ depends not on $M$ as a t-motive over $\n F_q(\th)$ (or an extension of $\n F_q(\th)$) but on its model over $\n F_q[\th]$. Namely, let $M_1$, $M_2$ be two t-motives over $\n F_q[\th]$ which are isomorphic as t-motives over $\n F_q(\th)$. In this case we have $L_{fmef}(M_1,T)=L_{fmef}(M_2,T)$ but not necessarily $L(M_1,T)=L(M_2,T)$.
\medskip
As an illustration, we can consider $\g C^n_{P}$, $\g C^n_{P_1}$ for $P$, $P_1$ from 0.1.2.5a. The Theorem 0.1.6 ( = [GL16], Theorem 3.3) gives us both $L(\g C^n_{P},T)$, $L(\g C^n_{P_1},T)$, and the formula [GL16], (5.6.1) gives us a relation between them, i.e. it indicates explicitly which Euler factors enter in their ratio\footnotemark \footnotetext{$P, \ P_1, \ Q, \ \g C^n_{*}$ of [GL16], Section 5.6 are the same as in the present paper. The first line of [GL16], Section 5.6 contains a typographic error: $\g C^n_{P}$, $\g C^n_{P_1}$ are different $\n F_q[\th]$-models.}.
\medskip
The formula [GL16], (5.6.1) is given without proof. It can be proved either using [Ga], Chapter 4, III, Definition 4.10, or comparing the values of determinants of $\g M(P,n,t)$, $\g M(P_1,n,t)$ (notations of [GL16]). If deg $Q=1$ this comparison is immediate; for deg $Q>1$ this can be a difficult problem.
\medskip
The definition of $L_{fmef}(M,T)$ is the following. Let $M$ be an Anderson t-motive over $ \n F_q(\theta)$. Let $f_*=(f_1,\dots,f_r)^{tr}$ be a basis of $M$ over $\n F_q(\theta)[ t]$ (here $tr$ means transposition: we consider elements of a basis as a column vector). Let $Q\in M_{r\times r}(\n F_q(\theta)[ t])$ be the matrix of multiplication by $\tau$ in this basis. Let $\goth P$ be an irreducible polynomial in $ \n F_q[\theta]$. $\g P$ is called good for ($M$, $f_*$) if all coefficients of all entries of $Q$ are integer at $\goth P$. For fixed ($M$, $f_*$) almost all $\g P$ are good for them.
\medskip
We need the following notation. For $a\in (\n F_q[\theta]/\goth P)[ t]$, $a=\sum c_i t^i$ where $c_i\in \n F_q[\theta]/\goth P$, we denote $a^{(k)}:=\sum c_i^{q^k} t^i$, for a matrix $A=(a_{ij})\in M_{r\times r}((\n F_q[\theta]/\goth P)[ t])$ \ $A^{(k)}:=(a_{ij}^{(k)})$ and $A^{[k]}:=A^{(k-1)}\cdot\dots\cdot A^{(1)}\cdot A$.
\medskip
Let us fix $f_*$ (hence $Q$ is also fixed). For a good $\g P$ the local $\goth P$-factor $L_\goth P(M,T)$ is defined as follows. Let $d$ be the degree of $\goth P$ and $\tilde Q\in  M_{r\times r}((\n F_q[\theta]/\goth P)[ t])$ the reduction of $Q$ at $\goth P$. We define:
$$L_\goth P(M,T):= \det(I_r-\tilde Q^{[d]}T^d)^{-1}\in \n F_q[t][[T^d]]\eqno{(0.1.3A.1)}$$ (because obviously det $(I_r-\tilde Q^{[d]}T)\in \n F_q[t,T]$ ) and the global $L$-function is, as usual, the product of local factors:
$$L_{fmef}(M,T):=\prod_{\goth P} L_\goth P(M,T)\in \n F_q[t][[T]]\eqno{(0.1.3A.2)}$$
This function is defined up to finitely many Euler factors (subscript $fmef$). It is easy to check that it does not depend on $f_*$. Really, let $g_*=(g_1,\dots,g_r)^{tr}$ be another basis of $M$ over $\n F_q(\theta)[t]$. Let $C\in M_{r\times r}(\n F_q(\theta)[t])$ be the transition matrix from $f_*$ to $g_*$ (we use its definition as $g_*=Cf_*$). We consider $\g P$ such that all coefficients of all entries of $C$ are integer at $\goth P$ and $\tilde C\in GL_r((\n F_q[\th]/\g P)[t])$ (almost all $\g P$ satisfy this condition). Let $Q_f$, $Q_g$ be the matrices of multiplication by $\tau$ in $f_*$, $g_*$ respectively. We have $Q_g=C^{(1)}Q_fC^{-1}$ and hence $Q_g^{[d]}=C^{(d)}Q_f^{[d]}C^{-1}$, $\tilde Q_g^{[d]}=\tilde C^{(d)}\tilde Q_f^{[d]}\tilde C^{-1}$.
\medskip
Since $\tilde C \in M_{r\times r}(\n F_{q^d}[t])$ we have $\tilde C^{(d)}=\tilde C$ and hence $$\det(I_r-\tilde Q_g^{[d]}T^d)=\det(I_r-\tilde Q_f^{[d]}T^d),$$ i.e. $L_\goth P(M,T)$ does not depend on a basis. This means that $L_{fmef}(M,T)$ for $f_*$ and $g_*$ coincide for almost all Euler factors.
\medskip
{\bf Example.} Let us find $L_{fmef}(\g C^n,T)$. Let $\goth P=\sum _{\vk=0}^d \al_\vk\th^\vk$ above be monic (i.e. its mayor coefficient $\al_d=1$). We have $Q=(t-\th)^n$ (see 0.1.1b), hence all $\g P$ are good for $\g C^n$ and its basis $e$. We have $\tilde Q=(t-\tilde \th)^n$ where $\tilde \th$ is the image of $\th\in \n F_q[\th]$ in $\n F_q[\th]/\g P$, and $\tilde Q^{(k)}=(t-\tilde \th^{q^k})^n$. Hence, $\tilde Q^{[d]}=\g P_{(t)}^n$ where $\g P_{(t)}:=\sum _{\vk=0}^d \al_\vk t^\vk$ (because $\g P$ is monic). The local Euler factor is $(1-\g P_{(t)}^nT^d)^{-1}$ and
$$L_{fmef}(\g C^n,T)=\sum_{\g P\in \n F_q[t] \hbox{ monic }}\g P^nT^{\deg \g P}\eqno{(0.1.3A.3)}$$

{\bf Remark 0.1.3A.4.} The above calculations show that the expression (0.1.3A.3) gives the value for $L_{fmef}(\g C^n,T)$. It is possible to prove a more strong theorem: the expression (0.1.3A.3) gives the value for $L(\g C^n,T)$. Since we do not give a formula for $L(\g C^n,T)$ (as it was written above, we do not need it), we do not give a proof of this theorem as well.
\medskip
We see that a  priory $L(\g C^n,T)\in \n F_q[t][[T]]$. Really, $L(\g C^n,T)\in \n F_q[t][T]$, see 0.1.7.1 below. A simple (without using the Lefschetz trace formula) proof of this fact is given in [G79], Section 4, Theorem 4.2. In notation of [G79] we have
$$\zeta(0, A, (s_0, -n)) = L(\g C^n, s_0^{-1}) \hbox{ (here } n > 0)$$
with $s_0^{-1}$ treated as a formal variable, see the formula at the bottom of [G79, page 110].
The vanishing of power series coefficients is established by [G79, Proposition 4.1].\footnotemark \footnotetext{The authors are grateful to an anonymous reviewer who indicated them this fact; we quote his text.} Also, this proof is given in [Th], proof of Theorem 2, and in [GL16], 7.1 and 7.2.
\medskip
{\bf Remark 0.1.3A.5.} In reality, not only for $M=\g C^n$ we have $L(M,T)\in \n F_q[t][T]$, but also for all Anderson t-motives $M$ (or, more generally, for effective t-motives in the meaning of [HJ]). This follows from [L], Section 2, formula (1).
\medskip
{\bf 0.1.3B.} The above $L$-function $L(M,T)$ considered in the present paper is a specialization (or a version) of the Goss $L$-function $L_G(M,s)$. More exactly, $L(M,T)$ contains strictly more information than the value of $L_G(M,0)$, i.e. it refines the value of $L_G(M,0)$.\footnotemark \footnotetext{The authors are grateful to an anonymous reviewer who indicated them this statement, as well as some other statements of Sections 0.1.3A, B.} The function $L_G(M,s)$ or its versions for $\tau$-sheaves and crystals are defined for example in [B05], Definition 15, [B02], Definition 2.19, the original definition is due to [G92], 3.4.2a.
\medskip
{\bf Remark 0.1.3B.0.} In most earlier cited papers $L_G(M,s)$ is considered as a function in variable $s$, case $s\le0$ (for $s\in \n Z$, $s>0$ we have a different theory). In terms of $\goth C^n$, this is the same as to consider $n$ as a variable: dilatation of $s$ to an integer $n$ corresponds to the tensor multiplication of $M$ by $\goth C^n$: see for example [B05], Proposition 9, or [GL16], (7.1) for $M=\g C^n$. Unlike this approach, we consider $n$ fixed, and we consider $T$ as a variable.
\medskip
Let us recall the definition of $L_G(M,s)$ given in [G96]. Following [G96], Definition 8.1.2.1, we let $S_\infty:=\p^*\times \n Z_p$. For $\al\in \n F_q((1/\th))^*$, $s\in S_\infty$ there exists $\al^s\in \p^*$ ([G96], Definition 8.1.2.2).
\medskip
Let $M$ be an Anderson t-motive over $\n F_q[\th]$ (we simplify the situation as possible), and $E=E(M)$ the corresponding $t$-module ([G96], below Remark 5.4.6). Let $v$ be a non-infinity valuation of $\n F_q(\th)$ ( = a prime ideal of $\n F_q[\th]$), and $T_v=T_v(E)$  the $v$-adic Tate module of $M$. Recall that it is the inverse limit by $\vk$ of $v^\vk$-torsion points of $E$ (here by $v$ we mean a generator of a prime ideal $v\subset \n F_q[\th]$).
\medskip
 We have: coordinates of torsion points of $E$ are algebraic over $\n F_q(\th)$, hence the Galois group Gal $(\overline{\n F_q(\th)}/\n F_q(\th)$ acts on $T_v(E)$.
\medskip
Now, let $w$ be another non-infinity valuation of $\n F_q(\th)$, i.e. $w\ne v$. For almost all $w$ its inertia subgroup $I_w$ acts trivially on $T_v(E)$. We consider for these $w$ the element $\s_w$ --- the arithmetic Frobenius element of $w$. It acts on $T_v(E)$. The inverse characteristic polynomial $P_w(M,U):=\det(1-\s_w|T_v(E)\cdot U)$ (here $U$ is an independent variable) belongs to $\n F_q[\th][U]$, it does not depend on a choice of $\s_w$ (it is defined up to $I_w$) and on a choice of $v$
\medskip
{\bf Remark.} For $w$ such that $I_w$ acts non-trivially on $T_v(E)$, this formula is modified as follows: $P_w(M,U):=\det(1-\s_w|(T_v(E))^{I_w}\cdot U)$.
\medskip
Finally, for $s\in S_\infty$ we have $$L_G(M,s)=\prod_{w}P_w(M,w^{-s})^{-1}$$ (the product is over the monic irreducible $w\in \n F_q[\th]$).
\medskip
For the Carlitz module $\g C$ the $L$-function $L_G(\g C, s)$ coincides with the zeta function of $\n F_q[\th]$ defined by the formula ([G96], end of (8.1)) $$\z_{\n F_q[\th]}(s):=\sum_{a\in \n F_q[\th] \hbox{ monic }}a^{-s}\eqno{(0.1.3B.1)}$$ Finally, for $n\in \n Z$ we define $s_n\in S_\infty$ (integer element in $S_\infty$) as in [G96], after 8.1.3.
\medskip
The relation between $L(M,T)$ and $L_G(M,s)$ for $M=\g C^n$ is the following.
Let $\g u: (\n F_q[t])[T]\to \p$ be the ring homomorphism defined by the formulas $\g u(t)=\th, \ \g u(T)=1$. Comparing (0.1.3A.3) and (0.1.3B.1) we get $$\g u(L(\g C^n,T))=\z_{\n F_q[\th]}(s_{-n})\eqno{(0.1.3B.2)}$$ (it holds on the level of summands of terms of both $L(\g C^n,T), \ \ \z_{\n F_q[\th]}(s_{-n})$ ). This formula also is given in [GL16], (7.1) and [Th], p. 233, middle of the page; $T$ of the present paper is $X$ of [Th].
\medskip
Further, there exists a Taelman $L$-function (or, better, a special value, see below) $L_{TA}(M/R)\in 1+T^{-1}\n F_q[[T^{-1}]]$ where $M$ is a Drinfeld module, $R$ is a finite extension of $\n F_q[\th]$, see [Ta], p. 371, below Remark 2. It is related with $L_G(M,s)$ by the following formula:
$L_{TA}(M/\n F_q[\th])$ is the value at $s=0$ (this explains the above terminology "a special value") of a version of $L_G(M,s)$ obtained by consideration of the dual of the Tate module of $M$ ([Ta], Remark 5).
\medskip
Let us compare possibilities to apply various definitions of $L$-functions to the statements of Birch --- Swinnerton-Dyer conjecture. The dual Goss $L$-function has a disadvantage that the analytic rank of all $M$ is 0 ((because $L_{TA}(M/R)$ --- which is its value at $s=0$ --- is non-zero). From another side, this value contains an algebraic and transcendental part, as it should be by the analogy with the number field case.
\medskip
The analytic rank $r_1$ defined as the order of zero at $T=1$ of the $L$-function $L(M,T)$ behaves like the analytic rank of elliptic curves (except the parity property), see 0.0.1d. From another side, $L(M,T)$ are polynomials in $T$, i.e. they are "too simple objects". Particularly, the coefficient $L^{(r_1)}(M,1)$ (entering to the Birch --- Swinnerton-Dyer conjecture; $L^{(r_1)}$ means the $r_1$-th derivative) belongs to $\n F_q[t]$, i.e. it has no transcendental part.
\medskip
Finally, we indicate that
\medskip
(a) Some results on vanishing of $L(\goth C^n,T)$ at $T=1$ were obtained in [Th], [L].
\medskip
(b) The zeta function of [AT90] is a specialization of the above $\z_{\n F_q[\th]}$ (i.e. its domain is a subset of $S_\infty$).

\medskip

{\bf 0.1.4.} An explicit formula for $L(\g C^n_P,T)$ (based on the general formula) is given in [GL16], Theorem 3.3. Let $P$ be from (0.0.1c). We denote $\bar \kappa=[\frac{m+n}{q-1}]$ (here $[x]$ means the integer part of $x$, i.e. the maximal integer which is $\le x$), and let $\Cal M(P,n,\bar \kappa)_{fc}=\Cal M(a_*,n,\bar \kappa)_{fc}$ be the matrix in $M_{\bar \kappa\times \bar \kappa}(\n F_q[t])$ whose $(\g i,\g j)$-th entry is defined by the formula

$$\Cal M(P,n,\bar \kappa)_{fc,\g i,\g j}=\sum_{l=0}^n (-\bar 1)^l\overline{\binom{n}{l}}a_{\g jq-\g i-l}\ t^{n-l}\eqno{(0.1.4a)}$$ (here $a_*=0$ if $*\not\in [0,\dots, m]$, and $\bar 1$, $\overline{\binom{n}{l}}$ are images of integer numbers $1$, $\binom{n}{l}$ in $\n F_p\subset \n F_q$).
\medskip
{\bf Remark.} Later (for example, in the proofs of Theorem 5.6 and Lemma 5.14.1) we shall consider linear transformations defined by matrices $\Cal M(P,n,\bar \kappa)_{fc}$ (and its characteristic 0 versions). Nevertheless, we do not know any natural interpretation of these linear transformations.
\medskip

In particular, for $n=1$ we have $\Cal M(P,1,\bar \kappa)_{fc,\g i,\g j}=a_{\g jq-\g i}t-a_{\g jq-\g i-1}$ and
\medskip
\noindent
$\Cal M(P,1,\bar \kappa)_{fc}=\left(\matrix a_{q-1}t-a_{q-2} &  a_{2q-1}t-a_{2q-2}   &   \dots & a_{\bar \kappa q-1}t-a_{\bar \kappa q-2}  \\ a_{q-2}t-a_{q-3}
 &  a_{2q-2}t-a_{2q-3} &   \dots & a_{\bar \kappa q-2}t-a_{\bar \kappa q-3} \\
a_{q-3}t-a_{q-4} &  a_{2q-3}t-a_{2q-4}  &   \dots & a_{\bar \kappa q-3}t-a_{\bar \kappa q-4} \\
\dots & \dots & \dots  & \dots \\ a_{q-\bar \kappa}t-a_{q-\bar \kappa-1}
 &  a_{2q-\bar \kappa}t-a_{2q-\bar \kappa-1}   &  \dots & a_{\bar \kappa q-\bar \kappa}t-a_{\bar \kappa q-\bar \kappa-1} \endmatrix \right) \ \ \ (0.1.5)$
\medskip
{\bf Theorem 0.1.6.} ([GL16], Theorem 3.3). $L(\g C^n_P,T)=\det(I_{\bar \kappa} - \Cal M(P,n,\bar \kappa)_{fc}T)$ (if $\bar \kappa=0$ then $L(\g C^n_P,T)=1$).
\medskip
{\bf 0.1.7. Remarks.} 1. Theorem 0.1.6 implies that $L(\g C^n_P,T)\in (\n F_q[t])[T]$ is of degree $\le \bar \kappa$ in $T$.
\medskip
2. $\g M(P,n,k)$ of [GL16], (3.1), (3.2) is $\Cal M(P,n,\bar \kappa)^{tr}_{fc}$. Authors apologize for this non-concordance of notations. Really, the transposition of these matrices is not important, because we consider their determinants.
\medskip
3. $\Cal M(P,n,\bar \kappa)^{tr}_{fc}$ is (up to a non-essential change of indices) a particular case of the matrix from [FP], (1.5).
\medskip
4. Formula for $L(\g C^n_P,T)$ is concordant with the natural inclusion of the set of polynomials of degree $\le m$ to the set of polynomials of degree $\le m'$, where $m'>m$.
\medskip
{\bf 0.1.8. Non-trivial part.} If $\frac{m+n}{q-1}$ is integer then the last column of $\Cal M(P,n,\bar \kappa)_{fc}$ has only one non-zero element, namely its lower element is equal to $(-1)^na_m$. Hence, for this case we denote $\kappa:=\bar \kappa-1=\frac{m+n}{q-1}-1$, we consider a $\kappa\times \kappa$-submatrix of $\Cal M(P,n,\bar \kappa)_{fc}$ formed by elimination of its last row and last column. We denote this submatrix by $\Cal M_{nt}(P,n,\kappa)_{fc}=\Cal M_{nt}(a_*,n,\kappa)_{fc}$ and $L_{nt}(\g C^n_P,T):=\det(I_{\kappa} - \Cal M_{nt}(P,n,\kappa)_{fc}T)$ (the subscript $nt$ means the non-trivial part). We have

$$L(\g C^n_P,T)=L_{nt}(\g C^n_P,T)\cdot (1-(-1)^na_mT)$$

{\bf 0.1.9. Analytic rank.} Let $(a_0,\dots, a_m)\in \n F_q^{m+1}$ be fixed, and $P$ from (0.0.1c). Let $n\ge1$ and $\g c\in \n F_q[t]$ be fixed.
\medskip
{\bf Definition 0.1.10.} The analytic rank of $\g C^n_P$ at $\g c$ is the order of 0 of $L(\g C^n_P,T)$ (considered as a function in $T$) at $T=\g c$. It is denoted by $r_\g c=r_\g c(n,P)=r_\g c(n,a_*)$.
\medskip
Since for all $P$ we have $L(\g C^n_P,T)=1$ at $T=0$, the simplest values of $\g c$ to study the analytic rank are $\g c$ = const. According to [GL16], Corollary 2.4.4, the case of any constant $\g c$ can be easily reduced to the case $\g c=1$, so $r_1$ is the most natural object to study.

Also, we can consider the order of pole of $L(\g C^n_P,T)$ at $T=\infty$ (i.e. the degree of $L(\g C^n_P,T)$ as a polynomial in $T$). It is more convenient to consider its deficiency:
\medskip
{\bf Definition 0.1.11.} The analytic rank of $\g C^n_P$ at $\infty$ is
$\bar \kappa- \deg_TL(\g C^n_P,T)$. It is denoted by $r_\infty=r_\infty(n,P)=r_\infty(n,a_*)$.
\medskip
{\bf Remark 0.1.12. } $r_\infty$ is not invariant under the natural inclusion of the set of polynomials of degree $m$ to the set of polynomials of degree $m'$, where $m'>m$.
\medskip
{\bf 0.1.13. Varieties of points of a given analytic rank, and polynomials defining them.} We recall definitions of (0.0.2a), (0.0.4). We fix $q,n$, \ $m>i\ge0$, and $\g c=1$ or $\g c=\infty$. We consider the set of points $(a_0,\dots, a_m)\in \n F_q^{m+1}$ such that $r_\g c(n,a_*)\ge i$. It is obviously an algebraic set, i.e. the set of $\n F_q$-zeroes of some polynomials in $a_0,\dots, a_m$. These polynomials come from coefficients of $\det(I_{\bar \kappa} - \Cal M(P,n,\bar \kappa)_{fc}T)$ considered as polynomials in $t,T$. We denote these polynomials by $D_{**fc}$, resp. $H_{**fc}$ for $\g c=1$, resp. $\g c=\infty$ (we do not give an exact definition of $D_{**fc}$, $H_{**fc}$ because we do not need it. See (0.2.1a) for the definition of $H_{**}$ for the case that we need). For fixed $q,n,m,i$ their set depends on 2 parameters, because they are coefficients of polynomials in $t$, $T$.
\medskip
Hence, we can consider algebraic varieties denoted by $X_\g c(q,n,m,i)_{fc}$ --- the sets of zeroes in $\bar \n F_q^{m+1}$ of $D_{**fc}$, resp. $H_{**fc}$ for $\g c=1$, resp. $\g c=\infty$. We have: the set of points $(a_0,\dots, a_m)\in \n F_q^{m+1}$ such that $r_\g c(n,a_*)\ge i$ is $X_\g c(q,n,m,i)_{fc}(\n F_q)$ --- the set of $\n F_q$-points of $X_\g c(q,n,m,i)_{fc}$.
\medskip
Polynomials $H_{**fc}$ are homogeneous, hence we can consider $X_\infty(q,n,m,i)_{fc}\subset \n P^m(\bar \n F_q)$ as projective varieties.
\medskip
{\bf Remark 0.1.14.} For the case of integer $\frac{m+n}{q-1}$ we can consider the non-trivial parts $r_{\g c,nt}$ of $r_\g c$, corresponding to zeroes of $L_{nt}(\g C^n_P,T)$. The same consideration is applied to the varieties $X_\g c(q,n,m,i)_{fc}$: we can consider varieties $X_{\g c,nt}(q,n,m,i)_{fc}$.
\medskip
{\bf Remark 0.1.15.} \footnotemark \footnotetext{The authors are grateful to an anonymous reviewer who indicated them the subject of the present remark.} We have: only the sets $X_\g c(q,n,m,i)_{set, \ fc}=X_\g c(q,n,m,i)_{fc}(\n F_q)$ are canonically defined (as sets of twists having $r_\g c\ge i$), but not $X_\g c(q,n,m,i)_{fc}$ as varieties. Really, by definition, varieties $X_\g c(q,n,m,i)_{fc}$ are sets of $\bar \n F_q$-zeroes of polynomials $D_{**fc}$, resp. $H_{**fc}$ for $\g c=1$, resp. $\g c=\infty$. In its turn, $D_{**fc}$, $H_{**fc}$ come from the proof of [GL16], Theorem 3.3. In principle, it can happen that we get another formula for $L(\g C^n_P,T)$ giving other systems of polynomials $\g D_{*}$, $\g H_{*}$. The sets of $\n F_q$-zeroes of $\g D_{*}$, resp. $\g H_{*}$ are the same sets $X_\g c(q,n,m,i)_{fc}(\n F_q)$ for $\g c=1$, resp. $\g c=\infty$ (because they are sets of twists having $r_\g c\ge i$), but varieties of $\bar \n F_q$-zeroes of $\g D_{*}$, $\g H_{*}$ can be another.
\medskip
We neglect this phenomenon, because at the moment we do not know other natural systems of polynomials that define $X_\g c(q,n,m,i)_{fc}(\n F_q)$. See [GL16], page 121 for a possible invariant definition of the dimension of $X_\g c(q,n,m,i)_{fc}(\n F_q)$.
\medskip
{\bf 0.2. Subject of the present paper: $\g c=\infty$, $q=2$, $n=0$, field = $\n C$.}
\medskip
As we mentioned above, study of behavior of $L(\g C^n_P,T)$ at $T=\infty$ is much simpler than at $T=1$, so in the present paper we shall consider only this problem. Further, we can consider $a_0,\dots, a_m$ as abstract elements, not necessarily as elements of $\n F_q$. Particularly, they can belong to $\n C$. For simplicity, we shall consider the case of (0.1.8): $\frac{m+n}{q-1}$ is integer, $\kappa:=\frac{m+n}{q-1}-1$, we  consider the non-trivial part of the $L$-function, and we omit the subscript "$nt$". Hence, we define a matrix $\Cal M(P,n,\kappa)\in M_{\kappa\times \kappa}(\n Z[t])$  using the same formula (0.1.4a), but in characteristic 0. Namely, the $(\g i,\g j)$-th entry of $\Cal M(P,n,\kappa)$ is defined as follows:

$$\Cal M(P,n,\kappa)_{\g i,\g j}:=\sum_{l=0}^n (-1)^l\binom{n}{l}a_{\g jq-\g i-l}\ t^{n-l}\eqno{(0.2.0)}$$
We consider its characteristic polynomial:
$$Ch(\Cal M(a_*,n,\kappa),T):=\det(I_{ \kappa} - \Cal M(P,n, \kappa)T)\in\n Z[a_0,\dots, a_m][t][T]\eqno{(0.2.1)}$$

Varieties $X_\g c(q,n,m,i)\subset \n C^{m+1}$ are defined like the above $X_\g c(q,n,m,i)_{fc}$ (recall that we omit the subscript $nt$). For $\g c=\infty$ we use the same notation for its projectivization, i.e. $X_\infty(q,n,m,i)\subset \n P^m(\n C)$.
\medskip
More exactly, for $\g c=\infty$ the polynomials $$H_{**}=H_{i\g j,qn}(m)=H_{i\g j,qn}(a_0,\dots,a_m)\in\n Z[a_0,\dots, a_m]\eqno{(0.2.1a)}$$ are defined by the following formula (they are coefficients of $Ch(\Cal M(a_*,n,\kappa),T)$ at $t^*T^*$):

$$Ch(\Cal M(a_*,n,\kappa),T)=\sum_{\iota=0}^\kappa \sum_{\g j=0}^{n(\kappa-\iota)}H_{\iota\g j,qn}(m)t^\g jT^{\kappa-\iota}\eqno{(0.2.2)}$$

Hence, $X_\infty(q,n,m,i)$ is the set of zeroes of $H_{\iota\g j,qn}(m)$ for all $\iota\in [0,\dots, i-1]$ and all $\g j\in [0,\dots, n(\kappa-\iota)]$. We shall  need also the corresponding projective schemes ($\iota$ and $\g j$ run over the same set $\iota\in [0,\dots, i-1]$, $\g j\in [0,\dots, n(\kappa-\iota)]$):
$$X_{S,\infty}(q,n,m,i):=\hbox{Proj }(\n C[a_0,\dots, a_m]/ <H_{\iota\g j,qn}(m)>)\eqno{(0.2.2a)}$$
(here $<H_{\iota\g j,qn}(m)>$ is the ideal generated by all $H_{\iota\g j,qn}(m)$, and the subscript $S$ indicates that we consider a scheme instead of a set of points as earlier).
\medskip
{\bf Remark 0.2.3.} According (0.1.2.7), for $q=2$ there are no twists. An exact relation between $L$-functions of $\g C_P$, $\g C_{PQ^{q-1}}$ (here $Q\in \n F_q[\th]$) is given in [GL16], (5.6.1) --- these $L$-functions are equal up to finitely many Euler factors. Clearly this formula holds only in characteristic $p$, and in characteristic 0 we have a non-trivial theory for $q=2$.
\medskip
{\bf 0.2.3a.} As it was mentioned in (0.0.7), we can consider formally polynomials $H_{**}$ for $n=0$. For this case the matrix $\Cal M_{nt}(a_*,0,\kappa)$ does not depend on $t$, and hence the polynomials $H_{**}$ depend on only one parameter, not on two parameters. For $q=2$ we denote them by $D(m,i):=H_{i0,20}(m)$ (see also (1.2.1); recall that $D(m,i)\in \n C[a_0,\dots, a_m]$ ).
\medskip
Now we can formulate the subject of the present paper: study of $X_\infty(q,n,m,i)$, $X_{S,\infty}(q,n,m,i)$ for the case $q=2$, $n=0$ (see (0.0.7)). Varieties $X_\infty(2,0,m,i)$, resp. schemes $X_{S,\infty}(2,0,m,i)$ are denoted simply by $X(m,i)$, resp. $X_S(m,i)$. Namely, we define (here $\iota=0,\dots,i-1$; see also Definition 1.3):
$$X(m,i)=\{(a_0:...:a_m)\in \n P^m|\ D(m,\iota)(a_0,\dots,a_m)=0\} \ \ \hbox{ (set of points)}$$
$$X_{S}(m,i):=\hbox{Proj }(\n C[a_0,\dots,a_m]/ <D(m,\iota)>)\ \ \hbox{ (projective scheme)}$$

It turns out that for $q=2$ the case of any $n$ is conjecturally the same as the case $n=0$. Namely, we have
\medskip
{\bf Conjecture 0.2.4.} For $q=2$, $\forall i,\g j,m,n$ $\exists \ \varkappa$ such that $$(H_{i\g j,2n}(m))^\varkappa\in \ <D(m,0),\dots,D(m,i)>$$ --- the ideal generated by
$D(m,0),\dots,D(m,i)$ (recall that polynomials $H_{i\g j,2n}(m)$, $D(m,*)$ are over $\n C$, not over $\n F_q$).
\medskip
This conjecture holds for $i=0$ and some other cases, see [GL16], (II.6) and A1 for numerical data. See also Proposition 10.5 for a simple proof for $i=0$, $n=1$, and [ELS] for two particular cases.
\medskip
{\bf Remark 0.2.5.} 1. Conjecture 0.2.4 implies that $X_{\infty}(2,n,m,i)$ (which is equal to Supp $X_{S,\infty}(2,n,m,i)$ --- the support of $X_{S,\infty}(2,n,m,i)$ ) does not depend on $n$ and is equal to $X(m,i)$, although the schemes $X_{S,\infty}(2,n,m,i)$ are different. In particular, multiplicities of their irreducible components depend on $n$. See, for example, [GL16], Conjecture 9.7.8: both $X_{\infty}(2,0,m,3)$, $X_{\infty}(2,1,m,3)$ consist of 4 irreducible components (because they coincide), but multiplicities of these components in $X_{S,\infty}(2,0,m,3)$, $X_{S,\infty}(2,1,m,3)$ are different (see the last two lines of the table of [GL16], Conjecture 9.7.8).
\medskip
2. Finding of analogs of Conjecture 0.2.4 for $q>2$ is a subject of further research. [GL16], Theorem 8.6a indicates that such analogs exist.
\medskip
3. Analogs of Conjecture 0.2.4 for $q>2$ show that $H_{i\g j,qn}(m)$ are highly dependent. Since $H_{**fc}$ of (0.1.13) is the reduction of $H_{i\g j,qn}(m)$ and $D_{**fc}$ of (0.1.13) are linear combinations of $H_{**fc}$ we can expect that  $D_{**fc}$ are also dependent, and hence dimensions of $X_1(q,n,m,i)_{fc}$ maybe are higher than the na\"{\i}ve parameter count predicts.
\medskip
{\bf 0.3. Contents of the present paper.} We give an explicit description of irreducible components of $X(m,i)$, and their rational parameterizations, in combinatorial terms. Some results are conditional: they depend on conjectures based on computer calculations. Their proof is a subject of further research.

There are two constructions. Irreducible components of $X(m,i)$ are denoted by $C_{ij\g k}(m)\subset X(m,i)$, where $j\ge1$ is an invariant of an irreducible component (i.e. $j$ is canonically defined by the component), and $\g k$ is a label of this component (i.e. assignment of $\g k$ is arbitrary), see 2.4, 2.5 for details. Conjecturally, they form series: $i,\ j, \ \g k$ are fixed and $m$ grows. The minimal possible value of $m$ is $i+j$, and the corresponding component is called the minimal component (see 2.4.1a).
\medskip
The first construction (Theorem 5.6; Conjecture 5.10) describes the minimal components.
Namely, let $T$ be a finite rooted binary tree\footnotemark \footnotetext{This $T$ has nothing common with the argument of the $L$-function.} (see (4.1)) and $F$ a forest --- a disjoint union of such trees. We consider a weight function $w$ --- a function on the set of nodes of $F$ satisfying some properties ((4.4), (a), (b), (c) or (a$'$), (b$'$), (c$'$)). There are 3 groups $\Aut(F)$, $\g G(F)$ and the Galois group $\Gal(\n Q(exp(2\pi\sqrt{-1} /2^d))/\n Q)$ (here $d$ is a sufficiently large number) acting on the set of pairs $(F,w)$ where $F$ is a forest and $w$ is a weight function on it.
\medskip
The main results of the first construction are the following:
\medskip
{\bf Theorem 5.6.} Any pair $(F,w)$ defines an irreducible component of $X(m,i)$. If two pairs $(F_1,w_1)$ and $(F_2,w_2)$ belong to the same orbit of $\g G(F)\rtimes\Aut(F)$ then they define the same component. The action of the Galois group on the set of $(F,w)$ is the Galois action on the corresponding components.
\medskip
{\bf Conjecture 5.10.} These irreducible components of $X(m,i)$ are the minimal ones. All minimal irreducible components of $X(m,i)$ are obtained by the above construction, and if two pairs $(F_1,w_1)$ and $(F_2,w_2)$ define the same component then they belong to the same orbit of $\g G(F)\rtimes\Aut(F)$.
\medskip
This conjecture is supported by computer calculations (see Tables A2.2, A3).
\medskip
The second construction (Proposition 7.1; Conjecture 7.3) describes all elements of a series in terms of the minimal component of this series. Proposition 7.1 is conditional, it depends on Conjecture 0.2.4.
\medskip
It turns out that Conjecture 0.2.4 and polynomials $H_{i\g j,21}$ for $q=2$, $n=1$ play an important role in the second construction, as well as in the construction of the odd lift, see Corollary 10.4. This is wonderful, because a priori polynomials
$D(m,i)$ "do not know" that they come from the theory of Carlitz modules, and that there is some relation between them and $H_{i\g j,21}$.
\medskip
{\bf 0.4. Subjects of further research.} (1) It is desirable to continue computer calculations of [GL16], Section 6, for the case $q=3$, $n=1$. Are there $P\in \n F_3[\th]$ such that $r_1(P)>6$? Are there more $P\in \n F_3[\th]$ such that $r_1(P)=6$? What are the dimensions of $X_1(3,1,m,i)_{fc}$ for $i=1$ and 2 (for other $i$ the calculation can be too difficult), or, at least, what are their growths as $m\to\infty$?
\medskip
(2) To prove all conjectures of the present paper. Particularly, we should either prove Conjecture 0.2.4, or find an independent proof of Proposition 7.1 and Conjecture 7.3. Also, we should prove that Theorem 5.6 gives an exhaustive description of irreducible components, i.e. prove Conjecture 5.10.
\medskip
(3) To find more properties of varieties $X(m,i)$ and their irreducible components. For example, from one side, they are surjective images of powers of $\n P^1$. From another side, there are some tautological sheaves on them (see 1.5). We should find these sheaves in terms of $O(n)$ on $\n P^1$.
\medskip
(4) Generalize the results of the present paper to other $q$ and $n$ (i.e. to the case $q$, $n$ arbitrary, $T=\infty$, characteristic 0), first to the case $q=3$, $n=1$ --- the first non-trivial case of Problem 0.0.2. Most likely we shall have to use $q$-ary trees for their description, instead of binary trees. Here the situation is more complicated, because the statement of Conjecture 0.2.4 should be modified. Particularly, we should study varieties of zeroes of $H_{i\g j}$ and not only of $D(m,i)$.
\medskip
(5) Solve Problems 0.0.2, 0.0.3. This is the case $q$, $n$ arbitrary, $T=1$, characteristic $p$. We think that to pass from $T=\infty$ to $T=1$ is the most complicated part of the work.
\medskip
(6) Further, we can generalize the above results to the case of other Drinfeld modules (not necessarily the Carlitz modules).
\medskip
{\bf 0.5. Organization of the paper.} Section 1 contains an elementary self-contained definition of varieties $X(m,i)$. We formulate in Section 2 conjectures on degrees,
multiplicities and Jordan forms of their irreducible components, and we describe their series $C_{ij\g k}(m)$.
Section 4 contains definitions concerning the weighted rooted binary trees, forests and weights --- objects that will be used later. Section 5 is the main part of the paper. It gives the first construction --- a description of the minimal irreducible components of $X(m,i)$ in terms of weighted rooted binary forests (Theorem 5.6; Conjecture 5.10). Further, formulas for the degree (Section 5.11), intersection with the trace hyperplane (Section 5.12), multiplicity (Section 5.13), Jordan form (Section 5.14), fields of definition and Galois action on them (Section 5.15) are given.
Section 6 contains explicit examples of constructions of Section 5 for some types of forests. Section 7 contains the second construction - the construction of series corresponding to a fixed minimal irreducible component. Sections 8 and 9 contain another constructions of irreducible components which are called respectively the even and odd lifts. Section 10 gives a relation between the odd lift and Conjecture 0.2.4. The appendix to the present paper contains some tables. Table A1 justifies Conjecture 0.2.4. Tables A2, A3 describe irreducible components of $X(m,i)$ for $i\le 6$ and for some other cases. They justify Conjectures 5.10, 7.3.
\medskip
{\bf Some notations.}
\medskip
\settabs 5 \columns
\+$a_0,...,a_m$&Coordinates in $\n P^m$\cr
\+$A_1$, $A_2$, $A$&Matrices used in the proof of (5.4)\cr
\+$\g a_0,...,\g a_\g m$&Analogs of $a_0,...,a_m$ for the odd lift, see Sections 9 -- 10\cr
\+$\al$&The number (label) of a tree in a forest\cr
\+$\Cal A_0$, $\Cal A_1$&Matrices, see (9.1), (9.2)\cr
\+$\aleph$&A map, see above (5.8.11.1)\cr
\+$b_0,...,b_{\vartheta}$&Coordinates of a point in a projective space, see Section 7\cr
\+$\g b$&Entries of matrices $A_1$, $A_2$, $A$, see (5.4)\cr
\+$\be_{m,i,j}$; $\be_{m,i}$&Some inclusions, see (2.4.0); their union, see below (8.4)\cr
\+$c_1,...,c_j$&Numbers parameterizing points in $X(m,i)$, see (5.0.1)\cr
\+$C$; $C_{ij\g k}$&Element of $Irr(i,j)$, see (2.4.2); (2.5)\cr
\+$\bar C_{ij\g k}$; $C_{ij\g k}(m)$&Irreducible components of $X(m,i)$, see (2.5)\cr
\+$C(F)$&Irreducible components corresponding to a forest $F$, see (5.15)\cr
\+$\g c$&Element of $\n F_q[t]\cup\infty$ --- order of zero of $L$-function at this point\cr
\+& is studied. See (0.1.9), (0.1.10)\cr
\+$\g c_*=(\g c_0,...,\g c_m)$&Numbers used in Section 5, see (5.4.4)\cr
\+$\g C$; $\g C^{n}$&Carlitz module; its $n$-th tensor power\cr
\+$\g C^{n}_P$&$P$-twist of the $n$-th tensor power of a Carlitz module, \cr
\+&see below (0.0.1c) and (0.1.2.4)\cr
\+$\p$&Complete algebraically closed field of characteristic $p$, see beginning of (0.1)\cr
\+$\Cal C_2$, $\Cal C$&Elementary correcting matrix; their union, see below (5.4.3)\cr
\+$Ch(\g M(m))$&Characteristic polynomial of $\g M(m)$, see Section 1\cr
\+$CH$&Chow ring\cr
\+$d$&Depth of a tree (forest)\cr
\+$D_{**fc}$&See (0.1.13), (0.0.2a)\cr
\+$D(m,i)$&Polynomials defining $X(m,i)$, see (0.2.3a), (1.2)\cr
\+$D_1(u)$, $D_2(u)$&Descendant trees of a ramification node $u$ of a tree, see (5.11.1a), (4.6)\cr
\+$\g D_*$&Conjectural system of polynomials defining varieties $X_\g c(q,n,m,i)_{fc}$, \cr
\+&see Remark 0.1.15\cr
\+$e$, $e_1$, $e_n$&Element of a basis of $\g C$; $\g C^{n}$, see (0.1.1a,b) and the above lines\cr
\+$\g d(C)$, $\g d(C_{ij\g k})$&Coefficient describing degree of an irreducible component, see (2.4.2.1)\cr
\+$\de$; $\de_\al$&Degree of im $\vf(T,w)$; im $\vf(T_\al,w)$, see below (5.11.7); (5.11.8)\cr
\+$f_*$&Quantities of elements in some equivalence classes, see the end of (5.11)\cr
\+$F$&Forest (disjoint union of finite rooted binary trees)\cr
\+$\g F_{ij}$&Set of pairs: forest having $i$ nodes consisting of $j$ trees, and weight on it, \cr
\+&see above (5.9.1)\cr
\+$\bar \g F_{ij}$&Its quotient set by an action of a group, see (5.9), (5.10)\cr
\+$\vf$&Map parameterizing irreducible components of $X(m,i)$, see (5.0)\cr
\+$\vf_{ij}$&Map from the set of pairs \{forest, weight\} to the set of irreducible \cr
\+&components of $X(m,i)$, see (5.9.1)\cr
\+$\bar \vf_{ij}$ &Its quotient map\cr
\+$\tilde \vf_{ij}$&$\bar \vf_{ij}\circ\g p_{ij}$, see (5.12))\cr
\+$\phi_*$&See (5.14)\cr
\+$\Phi$&Subfield of $\p$, see (0.1.0)\cr
\+$fc$&Abbreviation for the finite characteristic case\cr
\+$g$&Element of $\Aut(\g C^n)$, see (0.1.2), (0.1.2.2)\cr
\+$G_k$&Automorphism group of the complete tree of depth $k$, see (4.3)\cr
\+$\g g$&Element of $\g G$\cr
\+$\g G$&A group acting on weights on a tree, forest, see (4.5)\cr
\+$\ga$&Depth of a contraction of a tree, see (5.11.2)\cr
\+$\Ga_d$&Galois group $\Gal(\n Q(\g z_d)/\n Q)$, where $\g z_d=\zeta_{2^d}$\cr
\+$h_*$&Lengths of some final branches of a tree, see (5.14)\cr
\+$h(u)$&Height of a node of a tree, see (5.14), the third line\cr
\+$H$&Trace hyperplane, see (3.1)\cr
\+$H_{**}$&See (0.0.6), (0.2.1a) \cr
\+$H_{**fc}$&See (0.0.4), (0.1.13) \cr
\+$H_{i\g j,qn}(m)$&Polynomials --- coefficients of a characteristic polynomial, see (0.2.2)\cr
\+$\g h(\vk)$&Quantity of nodes $u$ of a forest such that $h(u)\le\vk$, see (5.14.4.1)\cr
\+$\g H_*$&Conjectural system of polynomials defining varieties $X_\g c(q,n,m,i)_{fc}$, \cr
\+&see Remark 0.1.15\cr
\+$\eta_\be$&See (4.2.2)\cr
\+$\th$&Independent variable, see (0.0.1c) and the above lines, and Section 0.1\cr
\+$\vartheta$&Dimension of $\n P^{\vartheta}$, Section 7\cr
\+$i$&Quantity of equations defining $X(m,i)$ in $\n P^m$; \cr
\+&$=$ quantity of nodes in a forest\cr
\+$Irr(X(m,i))$&See (2.3)\cr
\+$Irr(i,j)$&Sets describing irreducible components of $X(m,i)$, see (2.4)\cr
\+$Irr(i)$&Their union, see below (8.4)\cr
\+$j$&Quantity of trees in a forest; $=$ a constant in (2.4)\cr
\+$\g j$&Subscript for $H_{i\g j,qn}(m)$, see (0.2.2)\cr
\+$k$&Depth of a node in a tree, forest\cr
\+$k$&Depth of a complementary forest, see (5.15)\cr
\+$K$&Multiindex, see (5.0.2)\cr
\+$\g k$&The number (label) of an element in $Irr(i,j)$\cr
\+$\bar \kappa$, resp. $\kappa$&Size of the matrix $\g M$, resp. $\g M_{nt}$, see (0.1.4), (0.1.8)\cr
\+$\vk$&Occasional integer parameters in some parts of the paper\cr
\+$\Cal K(m,i)$&Inclusion of sets of irreducible components, see (2.4.1), (2.4a)\cr
\+$l$&Quantity of final nodes in a tree (forest)\cr
\+$L, \ L_\al, \ L_\g m$&Class of hyperplane in Chow ring of $\n P^*$, see (5.11.8); above (10.6)\cr
\+$(\g l_0:\g l_1)$&Point in $\n P^1$, see Section 9\cr
\+$L_e$, $L_o$&Even; odd lift of irreducible components, see (8.4); (9.10)\cr
\+$\g L_e$&Even lift on sets parameterizing irreducible components, see (8.5)\cr
\+$(\la_0:...:\la_m)$&Coordinates of a point in $\im \vf$, see (5.0.1)\cr
\+$\la_*(c_1,...,c_j)$&Polynomials defining the map $\vf$, see (5.5), (5.5D)\cr
\+$\La_1(m)$&Linear subspace, see (8.2)\cr
\+$m$&Dimension of the ambient space $\n P^m$. See $X(m,i)$\cr
\+$M$&Anderson t-motive, see (0.1.1)\cr
\+$\g m$&Analog of $m$ for odd lift, see Sections 9 -- 10\cr
\+$\g M(m)(*,...,*)$&Matrix (principal object of the paper), see Section 1.\cr
\+$\Cal M(a_*,n,\bar \kappa)_{fc}$&Matrix defining $L$-function, see (0.1.4a)\cr
\+$\Cal M_{nt}(a_*,n,\kappa)_{fc}$&Non-trivial submatrix of $\Cal M(a_*,n,\bar \kappa)$, see (0.1.8)\cr
\+$\Cal M(a_*,n,\kappa)$&The same matrix with entries in $\n Z$, see (0.2.0)\cr
\+$\mu(C)$, $\mu(C_{ij\g k})$&Multiplicity of an irreducible component, see (2.4.2.2)\cr
\+$n$&Exponent of a Carlitz module; dimension of an Anderson t-motive $M$\cr
\+$N$&Maximal degree, see (5.0.3)\cr
\+$\g n$&Depth of a ramification node in a tree, see (6.3)\cr
\+$\g n$&Occasional integer parameters in some parts of the paper\cr
\+$\g N$&See (5.11.8)\cr
\+$\nu$&Map used in the second construction, see Section 7\cr
\+$\xi$&Element of $(\n P^1)^j$, see (5.0.1)\cr
\+$O$&Order of non-ramification nodes in a tree, forest, see below (5.4.2)\cr
\+$\om$&An element of $\p^*$, see above (0.1.1c)\cr
\+$\om$&Linear map of projective spaces, see: below (5.11.7); and (5.11.8c)\cr
\+$p$&Prime, characteristic of the base field. Introduction\cr
\+$P$&Initial polynomial, see (0.0.1c)\cr
\+$\g p$&See (5.11.8a)\cr
\+$\g p_{ij}$&Projection, see (5.9.2)\cr
\+$\g P(i)$&The set of partitions of $i$\cr
\+$\g P(\g Y)$&Polynomial, see (5.5.1.1)\cr
\+$\pi_{m,i}$, $\pi(i,j)$&Maps describing the Jordan form, see (2.3), (2.4.3)\cr
\+$\tilde \pi$&A partition, see (5.14.4.1.1)\cr
\+$\Cal P$&Polynomial, see A2.1.1\cr
\+$pr$&Projection map, see (10.3)\cr
\+$\psi$&Segre map, see (5.11.8b)\cr
\+$q$&Power of $p$\cr
\+$q_1(u)$, $q_2(u)$&Quantities of nodes in $D_1(u)$, $D_2(u)$, see (5.13)\cr
\+$\n Q(d)$&$\n Q(exp(2\pi\sqrt{-1} /2^d))$, see (4.5.4)\cr
\+$\g q$ &Prime, see A2.1\cr
\+$r$&Introduction: Rank of an Anderson t-motive\cr
\+$r$&Sections 4--6: Root of a tree\cr
\+$R_u$&Subtree formed by $u$ and all its descendants, see (5.14)\cr
\+$r_1$, $r_\g c$, $r_\infty$&Introduction: Analytic rank, see (0.0.2), (0.1.10), (0.1.11)\cr
\+$r_1$, $r_2$&Sections 4, 5: Roots of $D_1(r)$, $D_2(r)$, see (4.6)\cr
\+$\rho$&Constant, see Lemma 5.11.1.1; lines below (5.12.1.4)\cr
\+$\varrho_*$&See (5.14)\cr
\+$s$&Any Segre map, see (5.11.8), (5.11.9)\cr
\+$S_\infty$&See above (0.1.3B.1)\cr
\+$\g s(T)$, $\g s(F)$&Partition associated to a tree, forest, see (5.14.4)\cr
\+$\g s^{-1}(\varsigma)_1$&Subset of trees in $\g s^{-1}(\varsigma)$, see below (5.14.5.1)\cr
\+$\sigma$&Symmetric polynomial, see (5.4), (5.5), (6) \cr
\+$\varsigma$&Partition, see (5.14)\cr
\+$t$&Generator of the Anderson ring $\p\{t,\tau\}$\cr
(Usually, this generator is denoted by $T$. We have to use the notation $t$ in order to avoid confusion with the argument of $L$-function.)
\+$\g t$&Regular map, see (5.0.2)\cr
\+$T$&Introduction: Argument of the $L$-function\cr
\+$T$&Sections 1 - 10: Finite rooted binary tree\cr
\+$\tau$&Introduction: Generator of the Anderson ring $\p\{t,\tau\}$\cr
\+$\tau$&Sections 4--7: Map on sets of weights, see (4.6)\cr
\+$u$&Non-ramification node of a tree, forest\cr
\+$u'$&The right neighbor of a non-final non-ramification node $u$ of a tree, forest\cr
\+$U$&Independent variable in a characteristic polynomial, see (1.2)\cr
\+$\Upsilon_{ij}$&Constant in a binomial formula, see (2.7)\cr
\+$\g u$&See (0.1.3B.2)\cr
\+$v(x)$&A column vector, see (5.1)\cr
\+$w$; $w_a$; $w_\mu$&Weight on a tree, forest; its additive, multiplicative form\cr
\+$W(T)$, $W(F)$&Set of weights on a tree, forest\cr
\+$X(m,i)$&Projective varieties --- main objects of the paper, see (0.0.7), (0.2.3a), (1.3B)\cr
\+$X_S(m,i)$&Schemes whose support is $X(m,i)$, see (0.2.3a), (1.3A)\cr
\+$X_{S,\infty}(q,n,m,i)$&Schemes, see (0.2.2a)\cr
\+$X_{*_1}(q,n,m,i)_{*_2}$&Varieties (here $*_1=\g c$ or $(nt,\g c)$ where $\g c =1$ or $\g c=\infty$,  \cr
\+&$*_2=\emptyset$ or "fc"), see (0.0.2a) -- (0.0.6), (0.1.13), (0.2)\cr
\+$\chi$&Veronese map, see below (5.11.7)\cr
\+$\chi_u$&A vector, see proof of (5.6)\cr
\+$Y_*$&Correspondence on $Irr(*,*)$, see (3.2)\cr
\+$Z$&Irreducible component of $X(m,i)$\cr
\+$Z_i$&Map, see (5.12.1)\cr
\+$\g z_d$&$=\zeta_{2^d}=exp(2\pi\sqrt{-1} /2^d)$\cr
\+$\Cal Z$&See (5.14.3)\cr
\medskip
{\bf 1. Definitions.} Let $m, \ i$ satisfy $m>i\ge1$. We give here explicit definitions of $\Cal M$, $D(m,i)$ and $X(m,i)$ for the case $q=2$, $n=0$, at $\g c=\infty$. We shall consider the non-trivial part of $\Cal M$ and $L$ (see 0.1.8), hence $\kappa=m-1$.
\medskip
Let $a_*=(a_0, ...,a_m)$ be any objects, $a_\vk=0$ for $\vk\not\in \{0,\dots,m\}$. The matrix $\Cal M(a_*,0,m-1)$ is denoted by $\g M(m)(a_0,...,a_m)$, i.e. it is a $(m-1)\times (m-1)$-matrix whose $(\al, \be)$-th entry is equal to $a_{2\be-\al}$
(this is $\g M(P,m)$ of [GL16], Section 9). If it is clear what objects $(a_0, ...,a_m)$ are kept in mind, we write $\g M(m)$ instead of $\g M(m)(a_0,...,a_m)$.
\medskip
{\bf Example.} For $m=6,7$ \ \ $\g M(m)$ are the following (their structure is slightly different for odd and even $m$):
\medskip
\newpage
\centerline{$\g M(6)(a_0,...,a_6)$  \ \ \ \ \ \ \ \ \ \ \ \ \ \ \ \ \ \ \ \ \ \ \ \ $\g M(7)(b_0,...,b_7)$ }
\medskip
$$\left(\matrix a_1&a_3&a_5&0&0\\ a_0&a_2&a_4&a_6&0 \\ 0&a_1&a_3&a_5&0 \\ 0&a_0&a_2&a_4&a_6\\ 0&0&a_1&a_3&a_5\endmatrix \right)
\ \ \ \ \ \   \left(\matrix b_1&b_3&b_5&b_7&0&0\\ b_0&b_2&b_4&b_6&0& 0 \\ 0&b_1&b_3&b_5&b_7&0 \\ 0&b_0&b_2&b_4&b_6&0\\ 0&0&b_1&b_3&b_5&b_7\\ 0&0&b_0&b_2&b_4&b_6 \endmatrix \right) \eqno{(1.1)}$$
\medskip
$\g M(m)$ is a truncated Hurwitz matrix, i.e. the Hurwitz matrix (see [H], p. 273 - 274) without the last row and the last column. Also, $\g M(m)$ is a permutation of rows and columns of the Sylvester matrix of two polynomials $P_{[1]}$, resp. $P_{[0]}$ whose coefficients are $a_i$ with odd $i$, resp. even $i$.
\medskip
Let $Ch(\g M(m))$ be the $(-1)^{m-1}\cdot$ characteristic polynomial of $\g M(m)$:
$$Ch(\g M(m))=|\g M(m)-U\cdot I_{m-1}|$$
\medskip
{\bf Definition 1.2.} $D(m,i)\in \n Z[a_0,\dots, a_m]$ are coefficients at $U^i$ of $Ch(\g M(m))$ considered as a polynomial in $U$:
$$Ch(\g M(m))=D(m,0)+D(m,1)\ U+D(m,2)\ U^2+\dots + D(m,m-2)\ U^{m-2}+(-U)^{m-1}\eqno{(1.2.1)}$$
They are homogeneous polynomials of degree $m-1-i$. Hence, $D(m,i)=H_{i0,0}(m)$ where $H_{i\g j,n}(m)$ are from (0.2.2).
\medskip
Particularly, we have: $D(m,0)=|\g M(m)|$, $D(m,m-2)=(-1)^m tr (\g M(m))=(-1)^m(a_1+a_2+...+a_{m-1})$ (the trace hyperplane).
\medskip
{\bf Example 1.2.2.} For $m=4$ we have
\medskip
(a) $D(4,0)=-a_0a_3^2-a_1^2a_4+a_1a_2a_3$;
\medskip
(b) $D(4,1)=a_0a_3-a_1a_2-a_1a_3+a_1a_4-a_2a_3$;
\medskip
(c) $D(4,2)=a_1+a_2+a_3$.
\medskip
{\bf Definition 1.3. A.} $X_S(m,i)\subset \n P^m(\n C)$ (subscript $S$ means scheme) is a projective scheme of the first $i$ polynomials $D(m,i')$, $i'=0,1,...,i-1$:
$$X_S(m,i):= \Proj \n C[a_0,\dots,a_m]/<D(m,0),\dots,D(m,i-1)>\eqno{(1.3.A.1)}$$
where (as in (0.2.4)) $<D(m,0),\dots,D(m,i-1)>$ is the ideal generated by $D(m,0),\dots,D(m,i-1)$.
\medskip
{\bf B.} $X(m,i)$ is the support of $X_S(m,i)$, i.e. the set of $(a_0:...:a_m)\in \n P^m(\n C)$ such that $\forall \ i' =0,\dots, i-1$ we have $D(m,i')(a_0,...,a_m)=0$.
\medskip
{\bf Example 1.4.} For $m=4$ \ $X(4,1)$ is a cubic threefold (see 1.2.2, (a)): $$X(4,1)=\{(x_0:...:x_4)\in \n P^4(\n C)\ |\ -x_0x_3^2-x_1^2x_4+x_1x_2x_3=0\}$$ $X(4,2)$ (resp. $X(4,3)$) is a surface (resp. a curve) in $\n P^4(\n C)$. See [GL16], 9.15 for more details on these varieties.
\medskip
{\bf 1.5.} There exist some locally free sheaves on desingularization of $X(m,i)$. Really, let $(a_0,...,a_m)$ be a $m+1$-uple of numbers such that $(a_0:...:a_m)\in X(m,i)$. We can consider a $(m-1)\times (m-1)$-matrix $\g M(m)(a_0,...,a_m)$ as a linear transformation of the affine space $\n C^{m-1}$. Kernels and images of powers of this transformation give us vector spaces associated to the point $(a_0:...:a_m)\in X(m,i)$. These vector spaces associated to any point of $X(m,i)$ define sheaves on $X(m,i)$. Their restrictions to irreducible components of $X(m,i)$ are locally free on their open subvarieties of non-singular points. After desingularization of the irreducible components, these sheaves become locally free. Conjecturally, all irreducible components of $X(m,i)$ are rational varieties, moreover, they are surjective images of powers of $\n P^1$. There exists a research problem to  express these sheaves in terms of $O(*)$ on $\n P^1$.
\medskip
{\bf 1.6.} The purpose of the present paper is to give a description of irreducible components of $X(m,i)$. In principle, there exists a description of these components as sets of zeroes of some homogeneous polynomials (such description exists for all algebraic varieties in $\n P^m$). Really, since they are (conjecturally) rational varieties, for all of them we give only a parametric description. This is because a description of algebraic varieties as sets of zeroes of polynomials usually is not convenient. For example, even the image of the simplest regular map $\chi: \n P^1 \to \n P^\de$ (the Veronese map) defined by the formula $\chi (c:{c'})=(c^\de: c^{\de-1}{c'}:...: c{c'}^{\de-1}:{c'}^\de)$ has no canonical description as a set of zeroes of polynomials.
\medskip
{\bf 1.7. Terminology on a field of definition.} An affirmation "A variety $X$ is defined over a field $F$" can have different meanings. In the present paper, we shall consider only the following meaning. Our $X$ are always fixed subsets of $\n P^m(\n C)$. For $t\in \n P^m(\n C)$ and $\s\in \Aut(\n C)$ the point $\s(t)\in \n P^m(\n C)$ is well-defined. In the present paper, an affirmation "A variety $X$ is defined over a field $F$" (where $X\subset \n P^m(\n C)$ ) will mean the following:
$\forall \ t\in X, \ \ \forall \ \s\in \Aut(\n C/F)\hbox{ we have } \s(t)\in X$.
\medskip
Particularly, this $F$ is not necessarily the minimal field of definition (which is the subfield of $\n C$ corresponding to the subgroup of all $\s\in \Aut(\n C)$ such that $\s(X)=X$).
\medskip
{\bf Example 1.7.1.} Let $\g z\in \n C$ be any constant. According the present terminology, the $\n P^1$ in $\n P^2$ given by the formula $x_0+\g zx_1+x_2=0$ (where $(x_0:x_1:x_2)$ are homogeneous coordinates in $\n P^2$) is defined over $\n Q(\g z)$, and this is its minimal field of definition (according other terminologies, $\n P^1$ is clearly defined over $\n Q$).
\medskip
Analogically, for $\s\in \Aut(\n C)$ the set $\s(X)$ is the set of $\s(t)$ for all $t\in X$.
\medskip
{\bf 2. Conjectural description of irreducible components of $X(m,i)$.}
\medskip
The below conjectures generalize and correct the ones of [GL16], (9.7), see A2 for details. They come from explicit calculations (see Tables A2.2, A3) and results and conjectures of Sections 5, 7.
\medskip
{\bf Conjecture 2.1.} Formula (1.3.A.1) represents $X_S(m,i)$ as a complete intersection.
\medskip
By definition of the complete intersection, this means that all irreducible components of $X(m,i)$ have codimension $i$ \ --- \ the quantity of polynomials

$D(m,0),\dots,D(m,i-1)$.
\medskip
Later we shall assume the truth of Conjecture 2.1 in statements of all other conjectures of Sections 2, 3.
\medskip
Let $\g z_d:=\zeta_{2^d}=exp(2\pi\sqrt{-1} /2^d)$, i.e. $\g z_1=-1$, $\g z_2=\sqrt{-1}$. In the present paper, number $i$ never means $\sqrt{-1}$, but always the quantity of equations defining $X(m,i)$ in $\n P^m$ = the quantity of nodes in a forest.
\medskip
{\bf Conjecture 2.2.} All irreducible components of $X(m,i)$ are defined over $\n Q(\g z_{(i-2)/2})$ for even $i$, $(\n Q(\g z_{(i-1)/2})\cap \n R)\cdot \n Q(\g z_{(i-3)/2})$ for odd $i$.
\medskip
Really, Section 5.15 contains an exact (conjectural) formula for the field of definition of any irreducible component of $X(m,i)$. Conjecture 2.2 is an immediate corollary of results of 5.15.
\medskip
{\bf 2.3. Jordan form, partition, multiplicity.} Let $Irr(X(m,i))$ be the set of irreducible components of $X(m,i)$, and let us denote by $\g P(i)$ the set of partitions of $i$.
Let $(a_0:...:a_m)$ be a generic point of an irreducible component of $X(m,i)$. The Jordan form of its matrix $\g M(m)$ has exactly $i$ zeroes on its diagonal. Really, by definition of $X(m,i)$, it has $\ge i$ zeroes on its diagonal. If for a generic point $(a_0:...:a_m)$ the Jordan form of $\g M(m)$ has $\ge i+1$ zeroes on its diagonal, then this irreducible component is contained in $X(m,i+1)$ which contradicts to Conjecture 2.1.
\medskip
{\bf 2.3.1. Partition.} Let us recall a definition of a partition associated to the 0-part of a Jordan form of a matrix. A Jordan form is a block diagonal matrix whose blocks are Jordan blocks. We consider all its 0-blocks (i.e. blocks having 0 on diagonal). Let there be $\vk$ such blocks of sizes $\g n_1,\dots, \g n_\vk$, and let $\g n_\bullet$ be the sum of these numbers. A partition of $\g n_\bullet$:
$$\g n_\bullet=\g n_1+...+\g n_\vk$$
depends only on the initial matrix and not on a choice of its Jordan form.
\medskip
For the above generic point $(a_0:...:a_m)$ we have $\g n_\bullet=i$, hence the above construction defines a partition of $i$.
\medskip
{\bf Conjecture 2.3.1.1.} Let $Z\subset X(m,i)$ be an irreducible component. There exists a Zariski open subset $Z_0\subset Z$ such that the above partition of $i$ is the same for all points $(a_0:...:a_m)\in Z_0$.
\medskip
See Remark 5.14.1.1 showing that this conjecture can be deduced from Conjectures 2.1, 5.8, 7.1, 7.3 of the present paper.
\medskip
This partition is called the partition corresponding to an irreducible component. This construction defines a map $\pi_{m,i}: Irr(X(m,i)) \to \g P(i)$. See examples in Section 5.14.
\medskip
{\bf 2.3.2. Multiplicity.} Let $Z$ be an irreducible component of $X(m,i)$. By its multiplicity $\mu(Z)$ we mean its multiplicity in $X_S(m,i)$.
\medskip
{\bf Remark 2.3.3.} Let us give a model example of the notion of multiplicity for non-experts in algebraic geometry. Let us consider a scheme

Spec $\n C[X,Y]/<Y^2-X^2+X^3, Y>$. We have: $Y^2-X^2+X^3=0$ is a singular cubic and $Y=0$ is a straight line passing through its singularity. Hence, the scheme has two points in its support: (0,0) (singularity) and (1,0). Their multiplicities in this scheme are 2 and 1 respectively.
\medskip
For another example, let us consider $X(4,1)$ from Example 1.4 (see [GL16], 9.15 for more details and explicit equations of $X(4,i)$ and their irreducible components). The set of its singularities $\Sing (X(4,1))$ is a plane (the reader should imagine the direct product of the above singular cubic, and a plane). The set $\{D(4,1)=0\}$ is a 3-dimensional cone in $\n P^4$ containing $\Sing (X(4,1))$. Hence, $X_S(4,2)$ has degree $3\cdot2=6$, and the plane $\Sing (X(4,1))$ has multiplicity 2 in it. It is the first of two irreducible components of $X_S(4,2)$. It is denoted by $C_{221}$, $m=4$ in A2.2. The second irreducible component of $X_S(4,2)$ has degree 4 and multiplicity 1, it is denoted by $C_{211}$, $m=4$ in A2.2.
\medskip
See also A2.1 for a method of explicit calculation of multiplicities, in terms of multiplicities of roots of some polynomials.
\medskip
{\bf Main conjecture 2.4.} For any $i\ge1$ and any $j=1,\dots,i$ there exist sets $Irr(i,j)$, and for any $m\ge i+j$ there exist canonical (see below for the meaning of this word) injective maps
$$\be_{m,i,j}: Irr(i,j)\hookrightarrow Irr(X(m,i))\eqno{(2.4.0)}$$ such that the following holds:
\medskip
{\bf 2.4.1.} $Irr(X(m,i))$ is a disjoint union of all im $\be_{m,i,j}$ (here and below $\sqcup$ means the disjoint union):
$$Irr(X(m,i))=\bigsqcup_{j=1}^{\hbox{min }(i,m-i)} \be_{m,i,j}(Irr(i,j))$$
Particularly, for $m\ge 2i$ the set $Irr(X(m,i))$ does not depend on $m$ (really, if $m\ge 2i$ then $\hbox{min }(i,m-i)=i$ and $Irr(X(m,i))=\bigsqcup_{j=1}^{i} Irr(i,j)$ - because $\be_{m,i,j}$ are inclusions). The word "canonical" means that if we fix an irreducible component $Z$ of $X(m,i)$ belonging to im $\be_{m,i,j}$ for some $j$, then the number $j$ and the irreducible component $$\be_{m+1,i,j}({\be_{m,i,j}}^{-1}(Z))\eqno{(2.4.1.1)}$$ of $X(m+1,i)$ are canonically defined by $Z$, and do not depend on any other choices (recall that $Z$ in (2.4.1.1) is considered as an element of a set $Irr(X(m,i))$, but not a component itself). Equivalently, we have a canonical inclusion denoted by $\Cal K(m,i)$, see 2.4a below: $$Irr(X(m,i))\hookrightarrow Irr(X(m+1,i)): \ \ Z\mapsto \be_{m+1,i,j}({\be_{m,i,j}}^{-1}(Z))$$ which is an isomorphism if $m\ge 2i$. See (2.5.2) below for an example. See also Picture 1.
\medskip
{\bf 2.4.1a.} For $m=i+j$ the irreducible components $\be_{i+j,i,j}(Irr(i,j))$ of $X(i+j,i)$ are called the minimal (irreducible) components. They are exactly elements of $Irr(X(m,i))- Irr(X(m-1,i))$ (difference of sets is with respect to the above inclusion).
\medskip
{\bf 2.4.2.} Let $C\in Irr(i,j)$. There exist integer positive numbers $\g d(C)$, $\mu(C)$ such that $\forall \ m\ge i+j$ we have:
$$\deg \be_{m,i,j}(C)=\g d(C)\binom{m-i}{j}\eqno{(2.4.2.1)}$$
(particularly, the degree of the minimal component corresponding to $C$ is $\g d(C)$ );
$$\hbox{$\forall \ m\ge i+j$ we have: }\mu(\be_{m,i,j}(C))=\mu(C)\eqno{(2.4.2.2)}$$

{\bf 2.4.3.} The partition of the Jordan form of $\be_{m,i,j}(C)$ does not depend on $m$, i.e. there exists a map $\pi(i,j): Irr(i,j) \to \g P(i)$ such that for any $m\ge i+j$ we have a commutative diagram
$$\matrix Irr(i,j)&&\overset{\be_{m,i,j}}\to{\to}&&Irr(X(m,i))\\
\\ &\pi(i,j)\searrow&&\swarrow \pi_{m,i}\\
\\ &&\g P(i)
\endmatrix$$
\medskip
{\bf 2.4a. Alternative formulation of 2.4.1.} Let us fix $i$. For any $m>i$ there exists a canonical inclusion $\Cal K_{m,i}: Irr(X(m,i)) \hookrightarrow Irr(X(m+1,i))$ which is an isomorphism if $m\ge 2i$.
\medskip
Relation with the above notations: $Irr(i,j)$ is defined as $Irr(X(i+j,i))-\Cal K_{i+j-1,i}(Irr(X(i+j-1,i))$ (we assume $Irr(X(i,i))=\emptyset$, hence $Irr(i,1)=Irr(X(i+1,i))$ ). Further, let $Z$ be an irreducible component of $X(m,i)$ belonging to im $\be_{m,i,j}$ for some $j$. We have $\Cal K_{m,i}(Z)$ is $\be_{m+1,i,j}({\be_{m,i,j}}^{-1}(Z))$ (see 2.4.1.1) in the above notations.
\medskip
Picture 1 illustrates 2.4.1 in terms of the inclusions $\Cal K_{m,i}$. The horizontal lines represent sets $Irr(X(m,i))$. Their subdivision as a disjoint union of $Irr(i,j)$ is shown by points subdividing the whole line by segments. Elements corresponding to minimal components are red. Maps $\Cal K_{m,i}$ are drawn as vertical hooked arrows. (2.4.1) is interpreted as follows: the partition of a Jordan form of a component depends only on its position in a column in Picture 1, but not on a horizontal line (i.e. it does not depend on $m$).
\medskip
{\bf 2.5. }We use the following notations. We order elements of $Irr(i,j)$ by some manner; the $\g k$-th element of $Irr(i,j)$ is denoted by $C_{ij\g k}$ (i.e. $\g k$ is not necessarily a number, but a label of an element of $Irr(i,j)$). The component $\be_{m,i,j}(C_{ij\g k})$ of $X(m,i)$ is denoted by $C_{ij\g k}(m)$. The minimal component corresponding to $C_{ij\g k}$ is denoted by $\bar C_{ij\g k}$; we have $\bar C_{ij\g k}=C_{ij\g k}(i+j)$.
\medskip
Let now $i,j,\g k$ be fixed, and let $m\ge i+j$ vary. A component $C_{ij\g k}(m)$ is an irreducible component of $X(m,i)$. The set of all these components (as $m$ varies) is called a series generated by $C_{ij\g k}$ (or by $\bar C_{ij\g k}$). To get the image of a series in Picture 1, we choose a point on any red segment (minimal component), and we consider points vertically down from it.
\medskip
A construction of minimal components is given in Theorem 5.6. Conjecture 5.10 affirms that this construction gives us all minimal components. Proposition 7.1, Conjecture 7.3 show how to construct all components of a series starting from its minimal component.
\medskip
{\bf Remark 2.5.1.} All irreducible components of $X(i+1,i)$ are the minimal ones, because the only possible value of $j$ for them is $j=1$, hence $m=i+j$. See Picture 1.
\medskip
{\bf Example 2.5.2.} For $i\le6$, resp. for $j\ge i-3$ sets $Irr(i,j)$, numbers $\g d(C_{ij\g k})$, $\mu(C_{ij\g k})$ and Jordan forms are given in Table A2.2, resp. A3. Namely, elements of $Irr(i,j)$ correspond to lines of Tables A2.2, A3 with these $i, \ j$.
\medskip
Let us illustrate Picture 1 and Table A2.2 for $i=5$. We have:
\medskip
The set $Irr(5,1)$ consists of 6 elements (lines $C_{511} - C_{516}$ of the Table A2.2);
\medskip
The set $Irr(5,2)$ consists of 5 elements (lines $C_{521} - C_{525}$ of the Table A2.2);
\medskip
etc.,
\medskip
The set $Irr(5,5)$ consists of 1 element (line $C_{551}$).
\medskip
Further, (2.4.1) or 2.4a imply that $Irr(X(6,5))$ consists of 6 elements (lines $C_{511} - C_{516}$ of the Table A2.2), all these components are the minimal ones;
\medskip
$Irr(X(7,5))$ consists of 11 elements (lines $C_{511} - C_{525}$ of the Table A2.2), components corresponding to the lines $C_{521} - C_{525}$ are the minimal ones;
\medskip
etc.,
\medskip
Finally, $Irr(X(10,5))$ consists of 16 elements (all lines of the Table A2.2 having $i=5$), the only minimal component of $Irr(X(10,5))$ is the component corresponding to the line $C_{551}$. For $m>10$ the set $Irr(X(m,5))$ consists of the same 16 elements, all of them are not minimal.
\medskip
{\bf Conjecture 2.6.} All irreducible components of $X(m,i)$ are rational varieties.
\medskip
This conjecture follows from Theorem 5.6, Conjectures 5.10, 7.3. Let us give an idea of the proof of this fact. First, formulas (5.5.1), (5.4.4) give us rational parametrization of varieties $\im \vf(F,w)$ (minimal irreducible components of $X(m,i)$). Results of Section 5.11 show that this rational parametrization is really a rational parametrization in variables $c_*^{2^{\ga_*}}$ where $\ga_*$ are some coefficients, see Proposition 5.11.4. Moreover, this construction gives us a rational 1 -- 1 correspondence between $(\n P^1)^j$ (a rational variety) and the corresponding minimal irreducible component of $X(m,i)$. Further, Conjecture 5.10 affirms that all minimal irreducible component of $X(m,i)$ are obtained using the construction of 5.4, 5.6, hence they are all rational. Finally, Conjecture 7.3 affirms that any irreducible component of $X(m,i)$ is a surjective image of a minimal irreducible component times a projective space. The map  $\nu_{i+j,\vartheta}$ from Conjecture 7.3 is an inclusion (see its definition in the beginning of Section 7), hence any irreducible component of $X(m,i)$ is a rational variety.

We do not give a complete proof of $\{(5.6), (5.10), (7.3)\} \implies \{(2.6)\}$: this is meaningless without having a proof of (5.10), (7.3).
\medskip
{\bf 2.7. Formula for $\sum_{\g k\in Irr(i,j)}\g d(C_{ij\g k})\mu(C_{ij\g k})$.}
\medskip
Most results of the present paper either are proved rigorously or have a serious justification. Exception is the formula (5.13.2) for $\mu(C)$: it is purely experimental. Its evidence comes not only from computer calculations, but also from the below formula (2.8). Let us deduce (2.8) from the above conjectures.

Since $X_S(m,i)$ is (conjecturally) a complete intersection and deg $D(m,i')=m-i'-1$ we get that deg $X_S(m,i)=(m-1)(m-2)\cdot...\cdot(m-i)$.
According (2.1) and (2.4.2), for any $m$ and $i$ the same number is equal to $$\sum_{j=1}^i\ [ \ \sum_{\g k\in Irr(i,j)}\g d(C_{ij\g k})\mu(C_{ij\g k})\ ] \ \binom{m-i}{j}$$

Let $i$ be fixed. If constants $\Upsilon_{ij}$ $\forall m$ satisfy $$\sum_{j=0}^i\ \ \Upsilon_{ij} \ \binom{m-i}{j}=(m-1)(m-2)\cdot...\cdot(m-i)$$
then they are defined uniquely. A well-known combinatorial formula usually written as
$$\binom{\g n}{\g l}=\binom{\varkappa}{0}\binom{\g n-\varkappa}{\g l}+\binom{\varkappa}{1}\binom{\g n-\varkappa}{\g l-1}+ \binom{\varkappa}{2}\binom{\g n-\varkappa}{\g l-2}+... +\binom{\varkappa}{\g l}\binom{\g n-\varkappa}{0}$$
gives us immediately that Conjectures 2.1, 2.4 imply $\forall i,j$ $$\sum_{\g k\in Irr(i,j)}\g d(C_{ij\g k})\mu(C_{ij\g k})=\binom{i-1}{i-j}\cdot i!\eqno{(2.8)}$$
\medskip
We see that formula (2.8) holds for all cases covered by Tables A2.2, A3, which gives evidence for Conjectures 2.1, 2.4, 5.8, 5.10. Example (a): $i=6$, $j=1$, Table A2.2 (this is the simplest case where there are two non-equivalent weights on a tree, see A2 for details). We have: $Irr(6,1)$ consists of 12 elements, (2.8) holds: $$32\cdot1+8\cdot2+16\cdot3+16\cdot4+16\cdot5+8\cdot6+4\cdot8+8\cdot10+8\cdot10+$$ $$+4\cdot10+8\cdot15+4\cdot20=720 \hbox{ and }\binom55\cdot6!=720$$ (b): $j=i-3$, Table A3. We have: $Irr(i,i-3)$ consists of 6 elements, (2.8) holds:
$$8(i-3)\cdot(i!)/24+8(i-3)(i-4)\cdot(i!)/12+8\binom{i-3}{3}\cdot(i!)/8+2(i-3)\cdot(i!)/12+$$ $$+4(i-3)\cdot(i!)/8+2(i-3)(i-4)\cdot(i!)/6=\binom{i-1}{3}\cdot i!$$

{\bf Remark.} Numbers $\g d(C_{ij\g k})$, resp. $\mu(C_{ij\g k})$ are given by (5.11.11), resp. (5.13.2), (5.13.3), and the set $Irr(i,j)$ is defined by Conjecture 5.10. Recall that these formulas are conjectural. There is a problem: to prove (2.8) using (5.11.11), (5.13.2), (5.13.3) as definitions of $\g d(C_{ij\g k})$, $\mu(C_{ij\g k})$, and Conjecture 5.10 as definition of $Irr(i,j)$. We do not think that this is easy, especially because there is no uniform description of the sets of forests and weights entering in 5.10.
\medskip
\newpage
{\bf 3. Intersection with the trace hyperplane.}\footnotemark \footnotetext{Contents of the present section, as well as Subsection 4.6, are used only for results of Subsection 5.12, and hence can be skipped while the first reading.} We shall look in this section what happens if we pass from $X(m,m-2)$ to $X(m,m-1)$. Namely, $X(m,m-2)$ is a surface in $\n P^m$ and $X(m,m-1)$ is its hyperplane section (recall that we assume truth of Conjecture 2.1). Really, we have $$X(m,m-1)=X(m,m-2) \cap \hbox{ the trace hyperplane }H:=\{a_1+...+a_{m-1}=0\}\eqno{(3.1)}$$ (because $D(m,m-2)=\pm(a_1+a_2+...+a_{m-1})$).

We let $i=m-2$. We shall see that formula (3.1) and Conjectures 2.1, 2.4 imply some relations between the sets $Irr(i,j)$ and numbers $\g d(C_{ij\g k})$, $\mu(C_{ij\g k})$ for $j=1,2$. Hence, we assume in this section that these conjectures hold. We have (identifying elements of $Irr(*,*)$ with their $\be_{*,*,*}$-images) $Irr(X(i+2,i))=Irr(i,1)\sqcup Irr(i,2)$ and $Irr(X(i+2,i+1))=Irr(i+1,1)$. Formula (3.1) implies that we have a correspondence ( = binary relation = multi-valued function) $$Y_i: Irr(i,1)\sqcup Irr(i,2) \to Irr(i+1,1)\eqno{(3.2)}$$ --- an irreducible component $Z \subset X(i+2,i)$ goes to $Z\cap H$ which is a union of irreducible components of $X(i+2,i+1)$.
\medskip
{\bf Conjecture 3.3.} The converse correspondence $Y_i^{-1}: Irr(i+1,1)\to Irr(i,1)\sqcup Irr(i,2)$ is a function, i.e. any irreducible component of $X(m,m-1)$ is contained only in one irreducible component of $X(m,m-2)$, or, the same, no curve --- intersection of irreducible components of $X(m,m-2)$ --- is contained in the trace hyperplane $H$.
\medskip
{\bf Remark.} Irreducible components of $X(m,m-2)$ are surfaces in $\n P^m$, so for large $m$ we can expect that their intersection is empty. Really, at least for $m=4$ it is not the case: [GL16], 9.15 shows that for $m=4$ the intersection of irreducible components of $X(4,2)$ is a curve, and Conjecture 3.3 for $m=4$ is true. The authors do not know what are intersections of irreducible components of $X(m,m-2)$ for $m>4$; this is a subject of further research.
\medskip
Conjecture 3.3 can be deduced from Conjectures 2.1, 5.10, 7.3 using explicit description of irreducible components given in Sections 5, 7. We do not need this deduction. See also Remark 5.12.2.
\medskip
Intersection with a hyperplane preserves both degree and multiplicity. For $C_{i1\g k}\in Irr(i,1)$ we have deg $C_{i1\g k}(i+2)=2\g d(C_{i1\g k})$, for $C_{i2\g k}\in Irr(i,2)$ we have deg $C_{i2\g k}(i+2)=\g d(C_{i2\g k})$, hence we get
\medskip
{\bf Proposition 3.4.} If Conjectures 2.1, 2.4 hold then
$$\hbox{ For $C_{i1\g k}\in Irr(i,1)$ we have }\sum_{\tilde C\in Y_i(C_{i1\g k})} \g d(\tilde C)=2\g d(C_{i1\g k})\eqno{(3.4.1)}$$
$$\hbox{ For $C_{i2\g k}\in Irr(i,2)$ we have }\sum_{\tilde C\in Y_i(C_{i2\g k})} \g d(\tilde C)=\g d(C_{i2\g k})\eqno{(3.4.2)}$$
$$\hbox{ For $C_{ij\g k}\in Irr(i,1)\sqcup Irr(i,2)$, $\tilde C\in Y_i(C_{ij\g k})$ we have }\mu(\tilde C)=\mu(C_{ij\g k}).\eqno{(3.4.3)}$$
$$(j=1, \ 2)\eqno{\square}$$
\medskip
{\bf 4. Rooted binary trees.} Minimal irreducible components of $X(m,i)$ are defined in terms of rooted binary trees, their unions and weights on them.
\medskip
{\bf 4.1.} A rooted tree is a tree having a fixed node (the root). The root defines an orientation on a tree (at all below pictures, the orientation is from the left to the right). A rooted tree is called binary if any node (including the root) has $\le 2$ right neighbors. Later in the present paper the term "tree" always (without any exception) will mean a finite rooted binary tree. The quantity of nodes of a tree is denoted by $i$.
\medskip
{\bf Example 4.1.1.} Here $i=23$:
$${\underset{23}\to{\bullet}}^{\nearrow^{{\underset{8}\to{\diamond}}^ {\nearrow^{\underset{6}\to{\circ}-{\underset{5}\to{\diamond}}<^{\circ2-\ast1}_{\circ 4-\ast3}}}
_{\searrow_{\ast7}}}}_ {\searrow_{\overset{22}\to{\circ}-{\underset{21}\to{\diamond}}^{\nearrow^{\circ13-\diamond12<^{\circ10-\ast9}_{\ast11}}}_{\searrow_{\underset{20}
\to{\circ}-{{\underset{19}\to{\diamond}}<^{\ast14}_{\overset{18}\to{\circ}-{\underset{17}\to{\diamond}}<^{\ast15}_{\ast16}}}  }    }     }}\eqno{(4.1.2)}$$ The nodes are numbered. All edges (arrows) are equivalent, they are drawn by different symbols only by typesetting reasons. Here
$\ast$ are final nodes, $\diamond$ are ramification nodes, $\bullet$ is the root node (it can be either a non-ramification node or a ramification node, as in the present example).
\medskip
{\bf Example 4.1.3.} All trees having $i\le5$, up to isomorphism (see Table A2.2, lines $i\le5$, $j=1$, last column):
\medskip
\+$i=1$&$i=2$&$i=3$&$i=4$&$i=5$\cr
\medskip
\+$\bullet$&$\bullet-\circ$&$\bullet-\circ-\circ$  &$\bullet-\circ-\circ-\circ$\ &$\bullet-\circ-\circ-\circ-\circ$  \cr
\medskip
\+&&$\bullet<^\circ_\circ$ &$\bullet-\diamond<^\circ_\circ $ & $\bullet-\circ-\diamond<^\circ_\circ$\cr
\medskip
\+&&&$\bullet<^{\circ-\circ}_\circ$ &$\bullet-\diamond<^{\circ-\circ}_\circ$\cr
\medskip
\+&&&&$\bullet<^{\circ-\circ-\circ}_\circ$\cr
\medskip
\+&&&&$\bullet<^{\circ-\circ}_{\circ-\circ}$\cr
\medskip
\+&&&&$\bullet<^{\diamond<^\circ_\circ}_{\circ}$\cr
\medskip
{\bf Remark 4.1.4.} The quantity of trees with $i$ nodes is $a_{we}(i+1)$ where $a_{we}(l)$ is a Wedderburn-Etherington number (see [W2]) --- the quantity of ramification types of trees with $l$ final nodes. It is given by a recurrent formula $$a_{we}(1)=1, \ \ \ a_{we}(\g n)=\sum_{i=1}^{\frac{\g n-1}2}a_{we}(i)a_{we}(\g n-i)\hbox{ for $\g n$ odd, } $$ $$ a_{we}(\g n)=\sum_{i=1}^{\frac{\g n-2}2}a_{we}(i)a_{we}(\g n-i)+\binom{a_{we}(\g n/2)+1}{2}\hbox{ for $\g n$ even}\eqno{(4.1.4.1)}$$  Apparently these formulas do not appear in the present theory.
\medskip
{\bf 4.2a. Definition.} The depth of a node is the distance between this node and the root (the length of any edge is 1; the depth of the root is 0), and the depth of a tree is the maximal depth of its nodes.
\medskip
{\bf 4.2b. Definition.} A simple tree is a tree without ramification nodes, its length ( = depth +1) is the quantity of nodes.
\medskip
A complete tree of depth $k$ is a tree such that all its final nodes have depth $k$ and all its non-final nodes are ramification nodes. Example for $k=3$: $\circ^{\nearrow^{\circ^{\nearrow^{\circ<^\circ_\circ}}_{\searrow_{\circ<^\circ_\circ}}    }}_{\searrow_{\circ^{\nearrow^{\circ<^\circ_\circ}}_{\searrow_{\circ<^\circ_\circ}} }}$
\medskip
A forest $F$ is a disjoint unordered union of trees, it is denoted $$F=\bigsqcup_{\al=1}^j T_\al\eqno{(4.2.1)}$$
\medskip
{\bf Remark.} Later we associate to $F$ some minimal irreducible components of $X(i+j,i)$ where $j$ is from (4.2.1) and $i$ is the quantity of nodes of $F$.
\medskip
Also, we can group isomorphic constituent trees as follows: $$F=\bigsqcup_{\be=1}^\de \eta_\be T_\be\eqno{(4.2.2)}$$ i.e. $F$ is a disjoint union of $\eta_1$ copies of $T_1$, $\eta_2$ copies of $T_2$, ... , $\eta_\de$ copies of $T_\de$, and $T_{\be_1}\ne T_{\be_2}$.
\medskip
{\bf 4.2.3.} By default, the quantity of nodes of a tree $T$ is denoted as above by $i$, the quantity of final nodes is denoted by $l$ (hence the quantity of ramification nodes is $l-1$), and the depth of $T$ is denoted by $d$. For a forest $F$ we use by default notations of (4.2.1) (i.e. $\al$ is the number (label) of a constituent tree in $F$, $j$ is the quantity of trees in $F$), numbers $i$, $l$, $d$ for $T_\al$ are denoted by $i_\al$, $l_\al$, $d_\al$ respectively. Numbers $i$, $l$, $d$ for $F$ mean the quantity of nodes, the quantity of final nodes and the depth respectively (the depth of a forest is, by definition, the maximal depth of its trees). Hence, we have $i=\sum_{\al=1}^j i_\al$, \ \ \ $l=\sum_{\al=1}^j l_\al$, \ \ \ $d=\max_{\al=1}^j d_\al$.
\medskip
A list of all forests having $i\le6$ is given in the last column of Table A2.2, their $i$ is indicated in the first column of the table, their $j$ is given in the second column.
\medskip
{\bf 4.3. Automorphisms of trees and forests.} Notation: Let $G$ be a group. By default, the semidirect product $G^n\rtimes S_n$ is with respect to the action of the symmetric group $S_n$ on $G^n$ by permutation of factors.
\medskip
An automorphism of a tree (a forest) is a permutation of nodes preserving the root(s) and the edges. We denote the automorphism group of the complete tree of depth $k$ by $G_k$. It satisfies a recurrent formula\footnotemark \footnotetext{Proof: eliminating the root we get two complete trees of depth $k-1$. An automorphism can either preserve them, or interchange. A subgroup of $G_k$ preserving them is $G_{k-1}^2$; interchanging gives us a semidirect product by $\n Z/2$.}: $G_1=\n Z/2\n Z$, $G_k=G_{k-1}^2\rtimes \n Z/2\n Z$. We have $\#(G_k)=2^{2^k-1}$. Clearly for any tree $T$ we have $Aut(T)=G_{\vk_1}\times G_{\vk_2}\times ...\times G_{\vk_\ga}$ for some $\vk_1,\vk_2,...,\vk_\ga$ and for a forest $F=\sqcup_{\be=1}^\de \eta_\be T_\be$ we have\footnotemark \footnotetext{This follows from a general fact. Let $X$ be a connected set with some structure and $X^{\sqcup n}$ the disjoint union of $n$ copies of $X$. Then $\Aut(X^{\sqcup n})=(\Aut(X))^n\rtimes S_n$.} $$Aut(F)=(Aut(T_1)^{\eta_1}\rtimes S_{\eta_1})\times (Aut(T_2)^{\eta_2}\rtimes S_{\eta_2})\times ... \times (Aut(T_\de)^{\eta_\de}\rtimes S_{\eta_\de})$$

Examples: $G_2$ is the dihedral group of order 8. The automorphism group of the tree (4.1.2) is $G_1\times G_1=(\n Z/2)^2$: the first factor interchanges branches (1 - 2 - 5) and (3 - 4 - 5), the second factor interchanges branches (15 - 17) and (16 - 17).
\medskip
{\bf 4.4. Weighted trees.} A weight on a tree is a function $w$ on its nodes
satisfying the below conditions. There are two equivalent descriptions of a weight: multiplicative and additive. Let $T$ be a tree. A multiplicative weight function $$w_\mu: \hbox{ \{nodes of }T\} \to \ <\g z_d>\ \subset \n C^*$$ where $<\g z_d>$ is the multiplicative group generated by $\g z_d$, satisfies conditions:
\medskip
(a) If $y$ is the left neighbor of $x$ then $w_\mu(y)=w_\mu(x)^2$;
\medskip
(b) If $y$ is a ramification node, $x_1$ and $x_2$ are its right neighbors then $$w_\mu(x_2)=-w_\mu(x_1)$$
\medskip
(c) $w_\mu$(root)=1.
\medskip
The corresponding additive weight $w_a$ is defined as follows. Let $x$ be a node of depth $k$. We have: $w_a(x)\in \n Z/2^k$ such that $w_\mu(x)=\g z_k^{w_a(x)}$. The above (a) - (c) become
\medskip
(a$'$) If $x$ is a right neighbor of $y$, where $y$ of depth $k$, then either $w_a(x)=w_a(y)$ or $w_a(x)=w_a(y)+2^{k}$ (here $\forall \ \vk$ $\ \ \n Z/2^\vk$ is identified with $\{0,1,\dots, 2^\vk-1\}$);
\medskip
(b$'$) If $y$ is a ramification node of depth $k$, $x_1$ and $x_2$ are its right neighbors then $w_a(x_2)\ne w_a(x_1)$, i.e. $w_a(x_2)= w_a(x_1)\pm2^k$;
\medskip
(c$'$) $w_a$(root)=0.
\medskip
A weight of a forest $F=\sqcup_{\al=1}^j T_\al$ is a set of weights $(w_1,...,w_j)$ of $T_1,...,T_j$. The set of weights of a tree $T$, resp. a forest $F$ is denoted by $W(T)$, resp. $W(F)$.
\medskip
{\bf 4.4.1.} For any $T$ we have $\#(W(T))=2^{i-l}$. Really, $T$ has $l-1$ ramification nodes. For each ramification node $x_\vk$ $(\vk=1,...,l-1)$ we choose and fix one of its right neighbors. We denote it by $x_\vk(1)$. Let $R$ be the set \{root node, $x_1(1),x_2(1),...,x_{l-1}(1)$\}. To choose an element $w$ of $W(T)$, we choose the values of $w(x)$ for all nodes $x\in T$ from the left to the right. For any $x\in T-R$ we have 2 choices, for any $x\in R$ we have 1 choice. So, $\#(W(T))=2^{\#(T-R)}=2^{i-l}$.
\medskip
Clearly $\#(W(F))=2^{i-l}$ as well ($i$, $l$ for $F$).
\medskip
{\bf 4.5. Action of groups on $W(T)$, $W(F)$.} There are 3 group actions on $W(F)$. The first action is the obvious action of $Aut(F)$ on $W(F)$.

 Let $F$ be a fixed forest. To define the second action, we let $\g G_\al=\g G(T_\al):=\n Z/2^{d_\al}$ and $\g G=\g G(F):=\prod_{\al=1}^j \g G_\al=\prod_{\al=1}^j\n Z/2^{d_\al}$. The second action is an action of the group $\g G$. It is defined as follows. Let $w_a$ be a weight written additively. Namely, if $x$ is a node of $T_\al\subset F$ of depth $k$ then $w_a(x)\in \n Z/2^k$. Let $\g g=(\g g_1,...,\g g_j)\in \g G$ where $\g g_\al\in \g G_\al$. We define $$(\g g(w_a))(x):=\bar\g g_\al+w_a(x)\eqno{(4.5.1)}$$ where $\bar\g g_\al:=\ve_{d_\al,k}(\g g_\al)$ is the image of $\g g_\al$ in the epimorphism $\ve_{d_\al,k}: \n Z/2^{d_\al}\to\n Z/2^{k}$.
\medskip
{\bf (4.5.2)} Formula (4.5.1) shows that this action of $\g G$ on $W(F)$ comes from the actions of $\g G_\al$ on $W(T_\al)$ defined by (4.5.1).
\medskip
{\bf Example 4.5.3.} Let $(F,w)$ be the following weighted forest (additive weights are indicated at nodes; elements of $\n Z/2^\vk$ are identified with $0, \dots, 2^\vk-1$):

\settabs 20 \columns
\medskip
\+&&&$^{\bullet0^{\nearrow^{\circ0<^{\circ0-\circ0}_{\circ2-\circ6-\circ6}}}_{\searrow_{\circ1<^{\circ3-\circ3-\circ11}_{\circ1-\circ5}  }}}$& &&&&$^{\bullet0-\circ1^{\nearrow^{\circ3}}_{\searrow_{\circ1-\circ5  }}}$
\cr

\medskip
We have $j=2$, $d_1=4$, $d_2=3$, $\g G=\n Z/16 \times \n Z/8$. We choose $\g g=(13,6)\in \g G$. The weight $\g g(w)$ is the following:

\medskip
\+&&&$^{\bullet0^{\nearrow^{\circ1<^{\circ1-\circ5}_{\circ3-\circ3-\circ3}}}_{\searrow_{\circ0<^{\circ0-\circ0-\circ8}_{\circ2-\circ2}  }}}$& &&&&$^{\bullet0-\circ1^{\nearrow^{\circ1}}_{\searrow_{\circ3-\circ3  }}}$
\cr

\medskip
We simply add 13, resp. 6 to weights of all nodes of $T_1$, resp. $T_2$, and take residues of obtained sums modulo $2^k$, where $k$ is the depth of a node.
\medskip

{\bf (4.5.4)} We denote $\n Q(d):=\n Q(\g z_d)$, $\Ga_d:=\Gal(\n Q(d)/\n Q)=(\n Z/2^d)^*$.
The third (multiplicative) action is the Galois action of $\Ga_d$ on values of $w_\mu(x)$. It is easy to see that in the above (additive) notations, it is defined by the formula $$(\ga(w_a))(x):=\bar\ga\cdot w_a(x)$$ where $\ga\in \Ga_d=(\n Z/2^d)^*\subset \n Z/2^d$ and $\bar\ga\in \n Z/2^k$ is as above.
\medskip
Clearly for one tree $T$ the action of $\Aut(T)$ commutes with the action of $\g G(T)$ on $W(T)$. For a forest $F$ such that some $\eta_\be$ of (4.2.2) are $>1$ the action of $\Aut(F)$ on $\g G(F)$ is not trivial.\footnotemark \footnotetext {Namely, for any $\be$ such that $\eta_\be>1$ we have a surjection $\Aut(F) \twoheadrightarrow S_{\eta_\be}$ (permutation of components). Element $\vk$ of $\Aut(F)$ acts on $\g g=(\g g_1,...,\g g_j)\in \g G$ permuting elements $\g g_\al$ according the images of $\vk$ in $S_{\eta_\be}$, for all $\be$.} We denote the corresponding semidirect product by $\g G(F)\rtimes\Aut(F)$. This group acts on $W(F)$. Further, the actions of $\Aut(F)$ and of $\Ga_d$ on $W(F)$ commute, and the actions of $\g G(F)$ and $\Ga_d$ give an action on $W(F)$ of the group $\g G(F)\rtimes\Ga_d$, the action of $\Ga_d$ on $\g G(F)$ in this semidirect product is the natural one (namely, $\ga\in \Ga_d=(\n Z/2^d)^*$ acts on $\g g=(\g g_1,...,\g g_j)\in \g G$ by the formula $\ga(\g g)=(\bar\ga\g g_1,...,\bar\ga\g g_j)$ where for $\al$-th term of this expression we have $\bar\ga\in (\n Z/2^{d_\al})^*$ as above).

\medskip
{\bf 4.6. Elimination of the root.}\footnotemark \footnotetext {Contents of the present subsection are used only for results of Subsection 5.12.} Let $T$ be a tree and $r$ its root. There are two possibilities:
\medskip
(a) $r$ is not a ramification node;
\medskip
(b) $r$ is a ramification node.
\medskip
In the case (a), we denote by $\bar T$ a tree obtained from $T$ by elimination of the root. We have a map $\tau: W(T) \to W(\bar T)$ defined as follows. Let $\bar r$ be the right neighbor of $r$ (i.e. $\bar r$ is the root of $\bar T$), and $w\in W(T)$ written additively. We have $w(\bar r)=0$ or $w(\bar r)=1$, and for any $u\in \bar T$ we have $w(u)\equiv w(\bar r)$ mod 2. We let $\forall \ u \in \bar T$ $$\tau(w)(u):=\frac{w(u)-w(\bar r)}2$$ Clearly $\forall \ w\in W(T)$, \ $\forall \ \g g\in \g G(T) \ \exists \ \bar \g g\in \g G(\bar T)$ (depending on $w$) such that $\tau(\g g(w))=\bar \g g(\tau(w))$. Warning: the map $\g g\mapsto \bar \g g$ (where $w$ is fixed) is not a homomorphism $\g G(T)\to \g G(\bar T)$. For any $\bar w\in W(\bar T)$ we have: $\tau^{-1}(\bar w)$ consists of two elements which belong to one $\g G(T)$-orbit. Particularly, $\tau$ is an isomorphism on the set of $\g G$-orbits.
\medskip
In the case (b), eliminating the root we get two trees. We denote them by $D_1(r)$, $D_2(r)$. Their disjoint union $D_1(r)\sqcup D_2(r)$ is a forest having $j=2$. Let $\al=1,2$. Two maps $\tau_\al: W(T) \to W(D_\al(r))$ are defined by the similar manner. Let $r_\al$ be the root of $D_\al(r)$, and $u\in D_\al(r)$. We let $$\tau_\al(w)(u):=\frac{w(u)-w(r_\al)}2\eqno{(4.6.1)}$$ Maps $\tau_1$, $\tau_2$ define a map $\tau: W(T) \to W(D_1(r)\sqcup D_2(r))$. For $(w_1,w_2)\in W(D_1(r)\sqcup D_2(r))$ we have $\tau^{-1}(w_1,w_2)$ consists of 2 elements, because $\forall \ w\in W(T)$ we have $w(r_1)\ne w(r_2)$.
\medskip
Clearly $\forall w\in W(T)$, \ $\forall \ \g g\in \g G(T), \ \forall \ \al=1,2 \ \exists \ \g g_\al\in \g G(D_\al(r))$ such that $$\tau_\al(\g g(w))= \g g_\al(\tau_\al(w))\eqno{(4.6.2)}$$

{\bf (4.6.3).} The converse usually is not true: if $\g g_\al\in \g G(D_\al(r))$ ($\al=1,2)$ then it is few likely that $\exists \ \g g\in \g G(T)$ such that (4.6.2) holds for both $\al=1,2$. Moreover, for a given $(w_1,w_2)\in W(D_1(r)\sqcup D_2(r))$ two elements of $\tau^{-1}(w_1,w_2)$ usually do not belong to one $\Aut(T)\times\g G(T)$-orbit.
\medskip
Example: Let $(T,w)$ be ${\overset{0}\to{\circ}<^{\overset{0}\to{\circ}-\overset{0}\to{\circ}- \overset{0}\to{\circ}- \overset{0}\to{\circ}}_{\underset{1}\to{\circ}-\underset{1}\to{\circ}- \underset{1}\to{\circ}}}$ (numbers at nodes indicate the additive weight), hence $D_1(r)=\circ-\circ-\circ-\circ$, $D_2(r)=\circ-\circ-\circ$, and $\tau_1(w)$, $\tau_2(w)$ are the 0-weights (i.e. an additive weights equal to 0 on any node). Let $(\g g_1, \g g_2)\in \g G(D_1(r)\sqcup D_2(r))$ be $(0,2)$. Two elements of $\tau^{-1}(\g g_1(w_1),\g g_2(w_2))$ are the following:
\medskip
${\overset{0}\to{\circ}<^{\overset{0}\to{\circ}-\overset{0}\to{\circ}- \overset{0}\to{\circ}- \overset{0}\to{\circ}}_{\underset{1}\to{\circ}-\underset{1}\to{\circ}- \underset{5}\to{\circ}}}$ and ${\overset{0}\to{\circ}<^{\overset{1}\to{\circ}-\overset{1}\to{\circ}- \overset{1}\to{\circ}- \overset{1}\to{\circ}}_{\underset{0}\to{\circ}-\underset{0}\to{\circ}- \underset{4}\to{\circ}}}$ \ \ \
Obviously there is no $\g g\in \g G(T)$ such that (4.6.2) for them holds. Further, $\Aut(T)$ is trivial, and these two weights are not in one $\g G$-orbit.
\medskip
{\bf 5. The first construction: a weighted forest defines a minimal irreducible component.} Let $F$ be a forest, and let $w\in W(F)$ be a weight of $F$. We let $m=i+j$ (see (4.2.3) for notations). We associate to the pair $(F,w)$ a map $$\vf=\vf(F,w): (\n P^1)^j \to \n P^m\eqno{(5.0)}$$ Later we shall see that its image is contained in $X(m,i)$. Since $(\n P^1)^j$ is irreducible, Conjecture 2.1 and the below Conjecture 5.8 imply that for any $(F,w)$ the image of $\vf(F,w)$ is an irreducible component of $X(m,i)$. Some explicit examples are given in (5.5.2 A, B), (6.1), (6.2.1), (6.2.3), (6.5.1 --- 6.5.3), (7.5.2), (7.5.3).
\medskip
The definition of $\vf(F,w)$ is given in (5.5.D), (5.5.1), numbers $\g c_*$ used in (5.5.1) are defined in (5.4.4). Nevertheless, it is more convenient to define $\vf(F,w)$ not by direct formulas (5.5.1), (5.4.4), but as a solution to some linear equations (5.2) - (5.3), see Section 5.1.1. This definition follows.
\medskip
{\bf 5.0.1.} Let $\xi=\{(c_1:c'_1), (c_2:c'_2), \dots, (c_j:c'_j)\}\in (\n P^1)^j$.  We denote $\vf(\xi)$ by $(\la_0:...:\la_m)\in \n P^m$. Numbers $\la_0,\dots,\la_m$ are polynomials in $c_1,c'_1, c_2,c'_2, \dots, c_j,c'_j$.
\medskip
Let $\n A^1:=\n P^1-\{(1:0)\}$ be the affine part of $\n P^1$ (the set of finite points) and $(\n A^1)^j\subset (\n P^1)^j$ their product. Any regular map $\g t: (\n A^1)^j\to \n P^m$ can be uniquely prolonged to a map from $(\n P^1)^j$ to $\n P^m$. Really, let $\g t$ be given by the formula $$\g t(c_1,\dots, c_j)=(\sum_K z_{K0}c^K: \sum_K z_{K1}c^K: ... :\sum_K z_{Km}c^K)\eqno{(5.0.2)}$$ where $K=(\vk_1,\dots,\vk_j)$ is a multiindex, $c^K=c_1^{\vk_1}\cdot c_2^{\vk_2} \cdot \dots\cdot c_j^{\vk_j}$, and $z_{K*}$ are coefficients. Let $N$ be the maximal degree of $K$ entering in (5.0.2). The prolongation of $\g t$ to $(\n P^1)^j$ is given by the formula
$$\g t(c_1:c'_1,\dots, c_j:c'_j)=(\sum_K z_{K0}c^K{c'}^{N-K}: \sum_K z_{K1}c^K{c'}^{N-K}: ... :\sum_K z_{Km}c^K{c'}^{N-K})\eqno{(5.0.3)}$$
where $N-K$ is the multiindex $(N-\vk_1,\dots,N-\vk_j)$ and ${c'}^{N-K}={c'_1}^{N-\vk_1}\cdot {c'_2}^{N-\vk_2} \cdot \dots\cdot {c'_j}^{N-\vk_j}$.
\medskip
Hence, in order to simplify formulas, we can consider the case of the affine part of each $\n P^1$ in $(\n P^1)^j$, i.e. $\forall \ \al=1,...,j$ we let $c'_\al=1$.
\medskip
We need a notation: for any $x$ we denote by $v(x)$ the following column vector of size $m-1\times1$: $$v(x):=(1,x,x^2,x^3,\dots,x^{m-2})^{tr}\eqno{(5.1)}$$
\medskip
We use notations of (4.2.1), (4.2.3), and we use the multiplicative form $w_\mu$ of $w$.
\medskip
{\bf 5.1.1.} The system of equations defining $\la_0,\dots,\la_m$ consists of:
\medskip
{\bf A.} One matrix equation for any final node $u$ of $F$. Let $k=k(u)$ be the depth of $u$ (see 4.2.a), and $\al$ the number of the tree containing $u$: $u\in T_\al\subset F$. The corresponding equation ($\la_0,\dots,\la_m$ are unknowns, $c_1,\dots, c_j$ parameters) is the following (here 0 is the 0-matrix column of size $m-1\times1$):
$$\g M(m)(\la_0,...,\la_m)\cdot v(w_\mu(u)c_\al^{2^{d_\al-k}})=0\eqno{(5.2)}$$
\medskip
{\bf B.} One matrix equation for any non-final, non-ramification node $u$ of $F$. Let $k$ and $\al$ be as above, and $u'$ the right neighbor of $u$. The corresponding equation is the following: $\exists$ a scalar $\be_u\in \n C$ depending on $u$ such that
$$\g M(m)(\la_0,...,\la_m)\cdot v(w_\mu(u)c_\al^{2^{d_\al-k}})=\be_u \cdot v(w_\mu(u')c_\al^{2^{d_\al-k-1}})\eqno{(5.3)}$$ holds.
\medskip
{\bf Proposition 5.4.} For any $(F,w)$, for any generic $(c_1,\dots, c_j)\in \n C^j$ equations (5.2), (5.3) have 1 non-zero solution, up to a multiplication by a scalar, i.e. one solution in $\n P^m$ (generic means that $(c_1,\dots, c_j)$ does not belong to a finite union of hypersurfaces).
\medskip
{\bf Proof.} For any number $\g b$ we define matrices
\medskip
$A_1(\g b):=\left(\matrix 0&1&0&\g b&0&\g b^2&0&\g b^3&\dots & \g b^{{\frac m2}-1}&0
\\1&0&\g b&0&\g b^2&0&\g b^3&0&\dots & 0& \g b^{\frac m2}\endmatrix \right)$ for even $m$,
\medskip
$A_1(\g b):=\left(\matrix 0&1&0&\g b&0&\g b^2&0&\g b^3&\dots & 0& \g b^{\frac {m-1}2}
\\1&0&\g b&0&\g b^2&0&\g b^3&0&\dots & \g b^{\frac {m-1}2} & 0\endmatrix \right)$ for odd $m$
\medskip
(for any $m$ it is a $2\times (m+1)$-matrix) and
\medskip
$A_2(\g b):=\left(\matrix 1&-\g b&\g b^2&-\g b^3&\g b^4&&\dots&&(-1)^m\g b^m\endmatrix \right)$ (a $1\times (m+1)$-matrix). Further, we denote a column matrix $(\la_0,...,\la_m)^{tr}$ by $\la_*$.
\medskip
(5.2) is equivalent to $$A_1(w_\mu(u)c_\al^{2^{d_\al-k}})\cdot \la_*=0\eqno{(5.4.1)}$$ and (5.3) is equivalent to $$A_2(w_\mu(u') c_\al^{2^{d_\al-k-1}})\cdot \la_*=0\eqno{(5.4.2)}$$ To prove these affirmations, it is sufficient to write down explicitly (5.2), (5.3), (5.4.1), (5.4.2), and we get a proof immediately. For example, (5.4.2) is equivalent to (5.3), because $w_\mu(u)c_\al^{2^{d_\al-k}}=(w_\mu(u')c_\al^{2^{d_\al-k-1}})^2$; $\be_u$ of (5.3) is equal to $\la_1+\la_3\g b^2+\la_5\g b^4+...$ where $\g b=w_\mu(u')c_\al^{2^{d_\al-k-1}}$.
\medskip
{\bf 5.4.2a.} Recall that $F$ is a forest. Let us make a list of all non-ramification nodes of $F$: $u_1,\dots,u_{i-l+j}$ and fix this order $O$ (see (4.2.3) for $l$). We form a $m\times (m+1)$-matrix $A=A(F,w,c_*)$ as follows. $A$ is a block matrix having one block column and $i-l+j$ block lines. For $\vk=1, \dots, i-l+j$ the $A_{\vk,1}$-th block entry of $A$ is a matrix $A_1(w_\mu(u_\vk)c_\al^{2^{d_\al-k_\vk}})$ (resp. $A_2(w_\mu(u'_\vk)c_\al^{2^{d_\al-k_\vk-1}})$) if $u_\vk$ is a final (resp. non-final, non-ramification) node of $T_\al\subset F$, where as above $u'_\vk$ is the right neighbor of a non-final node $u_\vk$, and $k_\vk$ is the depth of $u_\vk$. So, $A$ is a disjoint union of matrices $A_1(w_\mu(u_\vk)c_\al^{2^{d_\al-k_\vk}})$, $A_2(w_\mu(u'_\vk)c_\al^{2^{d_\al-k_\vk-1}})$ arranged vertically, according the order $O$, for all non-ramification nodes of $F$. See (6.1.1) for an example. (5.4.1), (5.4.2) are equivalent to $$A\la_*=0\eqno{(5.4.3)}$$
It is sufficient to show that for a generic $(c_1,\dots, c_j)$ the matrix $A$ is of the maximal rank $m$. We denote by $V_{\vk,\la}(x_1,\dots,x_{\vk})$ the Vandermonde matrix of size $\vk\times\la$ (its $(\mu,\nu)$-th entry is $x_{\mu}^{\nu-1}$). Let $\Cal C_2=\Cal C_2(\g b)$ (an elementary correcting matrix) be $\left(\matrix \sqrt{\g b}&1 \\ -\sqrt{\g b}&1\endmatrix \right)$, so we have $\Cal C_2(\g b)\cdot A_1(\g b)=V_{2,m+1}(\sqrt{\g b},-\sqrt{\g b})$. Let $\Cal C$ (a correcting matrix) be a block diagonal matrix consisting of $\Cal C_2(w_\mu(u)c_\al^{2^{d_\al-k}} )$ on block diagonal positions corresponding to final nodes $u$ of $F$ and of 1-s on the diagonal positions corresponding to non-final, non-ramification nodes. All blocks are ordered by the order $O$.
\medskip
{\bf 5.4.4.} Let us consider the sequence $u_1,\dots,u_{i-l+j}$ from 5.4.2a. Now we replace any final (resp. non-final, non-ramification) node $u$ by two numbers $\pm \sqrt{w_\mu(u)c_\al^{2^{d_\al-k}}}$ (resp. by one number $-w_\mu(u')c_\al^{2^{d_\al-k-1}}$, where $u'$ is the right neighbor). We denote the obtained sequence by $\g c_*=(\g c_1,\dots,\g c_m)$. See examples in 5.5.2.
\medskip
In these notations we have
$$\Cal C\cdot A(F,w,c_*)=V_{m,m+1}(\g c_*)\eqno{(5.4.5)}$$
Hence, (5.4.3) is equivalent to $$V_{m,m+1}(\g c_*)\cdot \la_*=0\eqno{(5.4.6)}$$

For any $(x_1,\dots,x_m)$ the $\vk$-th maximal minor (the determinant of the matrix obtained by elimination of its $(m+1-\vk)$-th column, $\vk=0,\dots,m$) of $V_{m,m+1}(x_1,\dots,x_m)$ is
$$\sigma_{m-\vk}(x_1,\dots,x_m)\cdot \prod_{m\ge\iota>\gimel\ge1}(x_\iota-x_\gimel)$$ where $\sigma$ is the symmetric polynomial (we let $\sigma_0=1$).
\medskip
{\bf 5.4.7.} For a generic $(c_1,\dots, c_j)$ the numbers $(\g c_1,\dots,\g c_m)$ are different. Really, first we can consider $c_\al\ne0$, hence $\sqrt{w_\mu(u)c_\al^{2^{d_\al-k}}}\ne-\sqrt{w_\mu(u)c_\al^{2^{d_\al-k}}}$. Further, let $u_1\in T_{\al_1}$, $u_2\in T_{\al_2}$ be two different nodes of depths $k_1$, $k_2$. If $\al_1\ne\al_2$ then equality $\g c_{*1}=\g c_{*2}$ (here $*1$, $*2$ are indices coming from $\al_1$, $\al_2$) can occur only if $c_{\al_1}$, $c_{\al_2}$ satisfy a non-trivial relation. For example, if $u_1, \ u_2$ are both non-final, non-ramification nodes then this relation is $-w_\mu(u'_1)c_{\al_1}^{2^{d_{\al_1}-k_1-1}}=-w_\mu(u'_2)c_{\al_2}^{2^{d_{\al_2}-k_2-1}}$; for the case when $u_1$ and/or $u_2$ is (are) final node(s) we have a similar relation. There are finitely many such relations, because there is one relation for a pair of nodes. If $\al_1=\al_2$ and $k_1\ne k_2$ then again equality $\g c_{*1}=\g c_{*2}$ can occur only if $c_{\al_1}$, $c_{\al_2}$ satisfy a non-trivial relation, because $c_{\al_1}$ enters in formulas for $\g c_*$ with different degrees. Finally, if $\al_1=\al_2$ and $k_1=k_2$ then coefficients $\pm \sqrt{w_\mu(u)}$, $-w_\mu(u')$ of $\g c_*$ for $u_1$, $u_2$ are always different: this follows immediately from the properties of $w_\mu$ (it is necessary to treat separately cases of final and non-final nodes $u_1$, $u_2$).
\medskip
So, if $(c_1,\dots, c_j)$ do not satisfy a finite number of polynomial relations, then all $(\g c_1,\dots,\g c_m)$ are different. Hence, for this case $\la_*\in \n P^m$ is unique. $\square$
\medskip
{\bf 5.5.} If $(\g c_1,\dots,\g c_m)$ are different then $$\la_\vk=(-1)^{m-\vk}\sigma_{m-\vk}(\g c_*), \ \ \vk=0,\dots,m\eqno{(5.5.1)}$$ Obviously $\sigma_{m-\vk}(\g c_*)$ are polynomials in $(c_1,\dots, c_j)$, we denote them by $\lambda_\vk (c_1,\dots, c_j)$. We can consider $\lambda_\vk (c_1,\dots, c_j)$ for any $c_1,\dots, c_j$, not necessarily $c_1,\dots, c_j$ satisfying the property that all $(\g c_1,\dots,\g c_m)$ are different.
\medskip
{\bf Remark 5.5.1.1.} (a) It is convenient (see (5.11.4)) to consider a polynomial $\g P(\g Y)=\g P_{F,w}(\g Y):=\prod_{\vk=1}^m (\g Y-\g c_\vk)=\sum_{\vk=0}^m\la_\vk \g Y^\vk$ where $\g Y$ is an abstract variable.\footnotemark \footnotetext{This polynomial has nothing common with polynomials defining $X(m,i)$ and/or its irreducible components, as algebraic varieties. $\g P(\g Y)$ will be used for calculation of the degrees of irreducible components of $X(m,i)$ in Section 5.11.} Here $\la_\vk=\la_\vk(c_1,\dots,c_j)$, i.e., strictly speaking, $\g P_{F,w}(\g Y)\in \n Q(\g z_d)[c_1,\dots,c_j,\g Y]$.
\medskip
(b) Polynomials $\la_\vk$ do not depend on the order $O$. Really, is we change $O$, then we get the corresponding change of the set $\g c_1,\dots, \g c_j$. Since $\la_\vk$ are symmetric polynomials in $\g c_1,\dots, \g c_j$, they do not depend on their order.
\medskip
{\bf Definition 5.5.D.} For any $\xi=\{(c_1:1), (c_2:1), \dots, (c_j:1)\}\in (\n P^1)^j$ we define $\vf(F,w)(\xi)\in \n P^m(\n C)$ as $(\lambda_0 (c_1,\dots, c_j):...:\lambda_m (c_1,\dots, c_j))$.
\medskip
We extend this definition to $(\n P^1)^j$ from $(\n A^1)^j$ by homogenization of polynomials $\lambda_\vk (c_1,\dots, c_j)$, see (5.0.3).
\medskip
{\bf Example 5.5.2. (A)} Let $F=T$ be the following tree: $\circ<^{\circ-\circ}_{\circ-\circ}$, $w_\mu=1$ on the root and on the nodes of the lower branch, $w_\mu=-1,\g z_2$ on the nodes of the upper branch. Let the order $O$ be the following: lower branch from the left to the right, upper branch from the left to the right. We have (here $j=1$, $c=c_1$, $c'=c'_1$, $i=5$, $m=6$):
$$\g c_*=(-c,\ \sqrt{c},\ -\sqrt{c},\ -\g z_2c,\ \sqrt{\g z_2c},\ -\sqrt{\g z_2c})$$ (the first 3 values correspond to the lower branch, the second 3 values correspond to the upper branch; from the left to the right), and
$$\vf(F,w)(c:1)=(-c^4:(-1+\g z_2)c^3:(1-\g z_2)c^3+\g z_2c^2:$$
$$:-2\g z_2c^2:(-1-\g z_2)c+\g z_2c^2:(1+\g z_2)c:1)\eqno{(5.5.2.1)}$$
(map on the affine part, i.e. $c'=1$);
$$\vf(F,w)(c:c')=(-c^4:(-1+\g z_2)c^3c':(1-\g z_2)c^3c'+\g z_2c^2{c'}^2:$$
$$:-2\g z_2c^2{c'}^2:(-1-\g z_2)c{c'}^3+\g z_2c^2{c'}^2:(1+\g z_2)c{c'}^3:{c'}^4)\eqno{(5.5.2.2)}$$
(map on whole $\n P^1$). This is $C_{515}$ from Table A2.2 (because of its degree which is 4, and the tree).
\medskip
{\bf (B)} Let $F$ be the following (nodes are numbered):
$\overset{6}\to{\circ}<^{\overset{3}\to{\circ}<^{\circ1}_{\circ2} } _{\underset{5}\to{\circ}-\underset{4}\to{\circ}}$ \ \ $\overset{11}\to{\circ}-\overset{10}\to{\circ}<^{\overset{7}\to{\circ} } _{\underset{9}\to{\circ}-\underset{8}\to{\circ}}$
\medskip
We have for it: $i=11$, \ $j=2$, \ $m=13$, \ $d_1=2$, \ $d_2=3$. We consider the following multiplicative weight function $w$; we indicate also the depth function $k$:
\medskip
\settabs 14 \columns
\+Node &&1&2&3&4&5&6&7&8&9&10&11\cr
\medskip
\+Weight $w$&&$-\g z_2$&$\g z_2$&$-1$&1&1&1&$-1$&1&1&1&1\cr
\medskip
\+Depth $k$&&2&2&1&2&1&0&2&3&2&1&0\cr
\medskip
Let the order $O$ of non-ramification nodes be the following: 1, 2, 4, 5, 7, 8, 9, 11. The sequence $\g c$ corresponding to this order is the following:
\settabs 6 \columns
\medskip
\+Elements of $\g c_*$ &&$\pm \sqrt{-\g z_2c_1}$&$\pm \sqrt{\g z_2c_1}$&$\pm \sqrt{c_1}$&$-c_1$\cr
\medskip
\+Corresponding nodes&&1&2&4&5\cr
\medskip
\+Elements of $\g c_*$ &&$\pm \sqrt{-c_2^2}$&$\pm \sqrt{c_2}$&$-c_2$ &$-c_2^2$\cr
\medskip
\+Corresponding nodes&&7&8&9&11\cr
\medskip
The polynomial $\g P(\g Y)(c_1,c_2)$ is equal to
$$(\g Y^4+c_1^2)(\g Y^2-c_1)(\g Y^2+c_2^2)(\g Y^2-c_2)(\g Y+c_1)(\g Y+c_2)(\g Y+c_2^2)$$
Hence, the corresponding irreducible component of $X(13,11)$ is defined over $\n Q$ (see 5.15). In non-homogeneous coordinates we have
$$\vf(F,w)((c_1:1), (c_2:1))=(c_1^4c_2^6:c_1^3c_2^6+c_1^4c_2^5+c_1^4c_2^4:-c_1^4c_2^5+c_1^4c_2^4+c_1^4c_2^3-c_1^3c_2^6+$$ $$c_1^3c_2^5+c_1^3c_2^4: ...: c_1c_2^2+c_2^3+c_1c_2+c_2^2-c_1-c_2:c_1+c_2+c_2^2:1)$$ In homogeneous coordinates we have
$$\vf(F,w)((c_1:c'_1), (c_2:c'_2))=(c_1^4c_2^6:c_1^3c'_1c_2^6+c_1^4c_2^5c'_2+c_1^4c_2^4{c'_2}^2:-c_1^4c_2^5c'_2  +c_1^4c_2^4{c'_2}^2+$$ $$c_1^4c_2^3{c'_2}^3-c_1^3c'_1c_2^6+c_1^3c'_1c_2^5c'_2+c_1^3c'_1c_2^4{c'_2}^2: ...: c_1{c'_1}^3c_2^2{c'_2}^4+{c'_1}^4c_2^3{c'_2}^3+$$ $$c_1{c'_1}^3c_2{c'_2}^5+{c'_1}^4c_2^2{c'_2}^4-c_1{c'_1}^3{c'_2}^6-{c'_1}^4c_2{c'_2}^5:c_1{c'_1}^3{c'_2}^6 +{c'_1}^4c_2{c'_2}^5+{c'_1}^4c_2^2{c'_2}^4:{c'_1}^4{c'_2}^6)$$
See (6.1), (6.2.1), (6.2.3), (6.5.1 --- 6.5.3), (7.5.2), (7.5.3) for other examples.
\medskip
{\bf Theorem 5.6.} $\forall \ F,w$ we have: im $\vf(F,w) \subset X(m,i)$.
\medskip
{\bf Proof.} We fix a generic $c_*=(c_1,...,c_j)\in (\n A^1)^j$, and let $\la_*$ be defined by (5.5.1) = (5.5.D). First, we need
\medskip
{\bf Lemma 5.6.1.} Let $u$ be any ramification node, $u'_1$, $u'_2$ its right neighbors, $\al,k$ as in (5.1.1.A), $c_\al\ne0$. Then $\exists$ scalars $\be_{1,u}, \ \be_{2,u}\in \n C$ such that
$$\g M(m)(\la_0,...,\la_m)\cdot v(w_\mu(u)c_\al^{2^{d_\al-k}})=\be_{1,u} \cdot v(w_\mu(u'_1)c_\al^{2^{d_\al-k-1}})+\be_{2,u} \cdot v(w_\mu(u'_2)c_\al^{2^{d_\al-k-1}})\eqno{(5.6.2)}$$
\medskip
{\bf Proof.} We let $\g b=w_\mu(u'_1)c_\al^{2^{d_\al-k-1}}$, hence $w_\mu(u'_2)c_\al^{2^{d_\al-k-1}}=-\g b$ and $w_\mu(u)c_\al^{2^{d_\al-k}}=\g b^2$. We have: if numbers $y_1,y_2$ are roots to the system

$$\matrix \la_1+\la_3\g b^2+\la_5\g b^4+...=y_1+y_2\\ \\ \la_0+\la_2\g b^2+\la_4\g b^4+...=\g b y_1-\g b y_2\endmatrix\eqno{(5.6.3)}$$ then values $\be_{1,u}=y_1$, $\be_{2,u}=y_2$ satisfy (5.6.2). System (5.6.3) has determinant $-2\g b\ne0$. $\square$
\medskip
According (5.4.7), for generic $c_*$ numbers $w_\mu(u)c_\al^{2^{d_\al-k}}$ are different; it is easy to see that this is true for ramification nodes $u$ as well. We consider these $c_*$. For any node $u$ of $F$ we denote the column vector $v(w_\mu(u)c_\al^{2^{d_\al-k}})$ by $\chi_u$. They are linearly independent, there are $i$ such vectors (because $i$ is the quantity of nodes in the forest; we have one column vector for a node. The matrix formed by these column vectors is the transpose of a Vandermonde matrix. If all numbers $w_\mu(u)c_\al^{2^{d_\al-k}}$ are different then the rank of this matrix is $i$ (the quantity of columns; this is because the determinant of a square Vandermonde matrix with different second-column entries is non-zero), hence these vectors are linearly independent).
\medskip
We order vectors $\chi_u$ according the value of $k(u)$ (see (5.1.1 A) for the definition of $k(u)$): if $k(u_1)<k(u_2)$ then the vector $\chi_{u_1}$ is preceding to $\chi_{u_2}$; if $k(u_1)=k(u_2)$ then the ordering of $\chi_{u_1}$, $\chi_{u_2}$ is arbitrary. Let $\psi_1,\dots, \psi_{m-i-1}$ be other vectors such that $\{\chi_u$, $\psi_*\}$ form a basis of $\n C^{m-1}$.
\medskip
We consider vectors $\{\chi_u$, $\psi_*\}$ as subdivided in two blocks: $\{\chi_1, \dots, \chi_i\}$ is the first block denoted by $\chi_*$, and $\psi_1,\dots, \psi_{m-i-1}$ is the second block denoted by $\psi_*$. Let us consider any linear operator on $\n C^{m-1}$. This block subdivision of the basis defines a $2\times2$ block structure on the matrix of this linear operator in this basis.
\medskip
Let this linear operator be the linear operator $\g M$ on $\n C^{m-1}$ defined by the matrix $\g M(m)(\la_0,...,\la_m)$ in the initial basis of $\n C^{m-1}$. Hence, the matrix of $\g M$ in the basis $\{\chi_*$, $\psi_*\}$ is a $2\times2$ block matrix. Formulas (5.2), (5.3), (5.6.2) show that this matrix is block triangular, namely its (2,1)-block is 0. This is because the action of $\g M$ on any $\chi_u$ is a linear combination of other $\chi_*$ (formulas (5.2), (5.3), (5.6.2)) and not of $\psi_*$.
\medskip
{\bf 5.6.4.} Further, the (1,1)-block of this matrix (corresponding to the action of $\g M$ on $\{\chi_*\}$ ) is strictly lower triangular (i.e. its diagonal entries are 0).
\medskip
\{ This is because the image of any $\chi_u$ under the action of $\g M$ is a linear combination of some $\chi_{\tilde u}$ where:
\medskip
For a non-final, non-ramification node $u$ we have: $\tilde u$ is $u'$ --- the right neighbor of $u$ (see (5.3));
\medskip
For a ramification node $u$ we have:
$\tilde u$ consists of $u'_1$, $u'_2$ --- the right neighbors of $u$ (see (5.6.2));
\medskip
$\g M(\chi_u)=0$ for a final node $u$, see (5.2). Since $k(u'), \ k(u'_1), \ k(u'_2)>k(u)$ we have that $\chi_{u'}, \  \chi_{u'_1}, \chi_{u'_2}$ are to the right from the vector $\chi_u$ in the ordering of the basis $\{\chi_*$, $\psi_*\}$. This implies strict lower-triangularity. \}
\medskip
Finally, for any $2\times2$ block matrix such that its (2,1)-block is 0 and its (1,1)-block of size $i\times i$ is strictly lower triangular we have: $i$ lower coefficients of its characteristic polynomial (in degrees $0, \dots, i-1$) are 0. This implies the statement of the theorem for a generic $c_*$. Since $X(m,i)$ is closed, we get the desired. $\square$
\medskip
{\bf Proposition 5.7.} $\forall \ g\in Aut(F)$, $\g g\in \g G$ we have: if $w_2=g(\g g(w_1))$ then $\im \vf(F,w_1)=\im \vf(F,w_2)$ (as subsets of $X(m,i)\subset \n P^m(\n C)$ ).
\medskip
{\bf Proof.} It is sufficient to prove this proposition independently for the cases $w_2=g(w_1)$, $w_2=\g g(w_1)$. Let us consider the case $w_2=g(w_1)$. In this case we have: for fixed $c_1,\dots, c_j$ the set of $\g c_1, \dots, \g c_m$ of $\{F,g(w_1), c_1,\dots, c_j\}$ is a permutation (induced by an automorphism $g$) of the set of $\g c_1, \dots, \g c_m$ of $\{F,w_1, c_1,\dots, c_j\}$, hence $\forall \ \vk$ the polynomials $\la_\vk(c_1,\dots, c_j)$ for $w_1$ and $w_2$ coincide.

Let $w_2=\g g(w_1)$, where $\g g=(\g g_1,...,\g g_j)\in \g G$ from (4.5). Let $u$ be a non-final, non-ramification node, resp. a final node of $F$ and $\g c(u,w,c_1,\dots, c_j)$ be an element, resp. two elements of the set $\g c_1, \dots, \g c_m$ corresponding to the node $u$ and the weight $w$ according (5.4.4). Definition (5.4.4) and the multiplicative form of (4.5.1) (formula for the action of $\g g$ on weights) imply immediately that for the above $w_1, \ w_2, \ \g g$ we have
$$\g c(u,w_2,c_1,\dots, c_j)=\g c(u,w_1,\g z_{d_1}^{\g g_1}\cdot c_1,\dots, \g z_{d_j}^{\g g_j}\cdot c_j)$$
hence three maps: an isomorphism $\{c_\al\mapsto \g z_{d_\al}^{\g g_\al}\cdot c_\al\}$, and $\vf(F,w_1)$, and $\vf(F,w_2)$, form a commutative triangle. This means that $\im \vf(F,w_1)=\im \vf(F,w_2)$. $\square$
\medskip
{\bf Conjecture 5.8.} $\forall \ F,w$ dim im $\vf(F,w)=m-i$.
\medskip
{\bf Remark 5.8.1.} The above conjecture is not surprising, because $j=m-i$, the source of $\vf$ is $(\n P^1)^j$ and $\vf$ "cannot have a fibre of dimension $>0$". This can be easily checked for any fixed $F$, $w$ (see, for example, Example 6.2.3), but a proof for all $F$, $w$ is too complicated.
\medskip
{\bf Corollary 5.9.} Conjectures 2.1, 5.8 imply that $\forall \ F,w$ we have: im $\vf(F,w)$ is an irreducible component of $X(m,i)$.\footnotemark \footnotetext{Really, let $\vf: X \to Y$ be a map from an irreducible variety $X$ of dimension $n$ to a variety $Y$ such that all irreducible components of $Y$ have the same dimension $n$, and such that $\dim \im \vf=n$. Then $\im \vf$ is an irreducible component of $Y$.} $\square$
\medskip
For fixed $i,j$ we denote by $\g F_{ij}$ the set of pairs $(F,w)$ for all $F$ whose $i,j$ are from (4.2.3), i.e $F$ consists of $j$ trees and has $i$ nodes. We have:
\medskip
{\bf 5.9.1.} Let us assume truth of Conjectures 2.1, 5.8. Under this assumption, Corollary 5.9 implies that there exists a map $\g F_{ij}\to Irr (X(m,i))$: $(F,w)\mapsto \im \vf(F,w)$. We denote it by $\vf_{ij}$.
\medskip
{\bf 5.9.2.} Let $\bar \g F_{ij}$ be the quotient set of $\g F_{ij}$ by the equivalence relation of Proposition 5.7, i.e. $(F,w_1)$ is equivalent $(F,w_2)$ iff $\exists \ g\in Aut(F)$, $\g g\in \g G$ such that $w_2=g(\g g(w_1))$. The natural projection $\g F_{ij}\to\bar \g F_{ij}$ will be denoted by $\g p_{ij}$. Proposition 5.7 implies (here and below we assume the truth of Conjectures 2.1, 5.8) that $\vf_{ij}$ factors throw $\bar \g F_{ij}$.
\medskip
Let us recall some definitions of Conjecture 2.4. There exist (abstract) sets $Irr(i,j)$ and injective maps (2.4.0)
$$\be_{m,i,j}: Irr(i,j)\hookrightarrow Irr(X(m,i))$$
For $m=i+j$ the elements $\be_{m,i,j}(Irr(i,j))$ are called the minimal irreducible components of $X(m,i)$, see 2.4.1a. We identify $Irr(i,j)$ with a subset of $Irr(X(m,i))$ via the inclusion $\be_{m,i,j}$.
\medskip
{\bf Conjecture 5.10.} We assume truth of Conjectures 2.1, 2.4, 5.8. In notations of (2.4), (2.5), im $\vf_{ij}=Irr(i,j)\subset Irr(X(m,i))$. Hence, there exists a map $\bar \vf_{ij}$ making the following diagram commutative:
$$\matrix \g F_{ij}&\overset{\vf_{ij}}\to\to&Irr(X(m,i))\\ \\ \downarrow\g p_{ij}&&\uparrow\be_{m,i,j}\\ \\ \bar \g F_{ij}&\overset{\bar \vf_{ij}}\to\to& Irr(i,j) \endmatrix \eqno{(5.10.1)}$$ Moreover, we conjecture that $\bar \vf_{ij}: \bar \g F_{ij}\to Irr(i,j)$ is an isomorphism.
\medskip
We see that Conjecture 5.10 describes the set $Irr(i,j)$ in combinatorial terms: it is the quotient set of $\g F_{ij}$ (a set of pairs \{forest, weight\}) by the equivalence relation of Proposition 5.7 (see 5.9.2). Since for any $m$ the set $Irr(X(m,i))$ is a union of some $Irr(i,j)$, we get that Conjectures 5.10, 2.4.1 give a complete description of $Irr(X(m,i))$ in combinatorial terms.
\medskip
{\bf Remark 5.10.2.} (1) Particularly, $\forall \ F,w$ we have: im $\vf(F,w)$ is a minimal component of $X(m,i)$ (because its dimension is $j$ which is $m-i$).
\medskip
(2) Clearly Conjecture 5.10 has meaning if Conjecture 2.4 is true (if not then the sets $Irr(i,j)$ have no meaning).
\medskip
(3) The origin of Conjecture 5.10 is purely experimental. We proved (Theorem 5.6) that im $\vf(F,w) \subset X(m,i)$. Subsequent conjectures indicate that $\vf(F,w)$ gives us irreducible components of $X(m,i)$. Computer calculations (Table A2.2) show that $\vf(F,w)$ gives us:
\medskip
(a) Only minimal irreducible components of $X(m,i)$ (i.e. not non-minimal components);
\medskip
(b) All minimal irreducible components of $X(m,i)$;
\medskip
(c) Any minimal irreducible component of $X(m,i)$ is obtained only once (up to the action of $Aut(F)$, $\g G$ on $\g F_{ij}$, see (5.7)).
\medskip
This is exactly the contents of Conjecture 5.10.
\medskip
(4) To prove Conjecture 5.10 we need to show that all generalized eigenvectors\footnotemark \footnotetext{Recall that a vector $v$ is a generalized eigenvector of a linear transformation $\vf$ with generalized eigenvalue $\la$ if $\exists \ \vk$ such that $(\vf-\la \cdot \Id)^\vk(v)=0$.} of $\g M(m)(\la_0,...,\la_m)$ with generalized eigenvalues 0 (here $(\la_0,...,\la_m)$ belongs to a minimal component of $X(m,i)$ ) have the form described in (5.1), for $x=$ one of $\g c_\al$. Also, we can use (2.8): namely, we must prove that for fixed $i,\ j$, for all possible $F$ having the given $i,\ j$, and for all their $w$, for the corresponding $\g d(C_{ij\g k})$, $\mu(C_{ij\g k})$ we have: (2.8) holds. This will mean that there is no other components. Taking into consideration that apparently there is no good description of $F$ and $w$ (see below), it will not be easy to realize this idea.
\medskip
{\bf 5.11. Finding of the degree of} $\im \vf(F,w)$. First, we consider the case of one tree: $F=T$, i.e. $j=1$, hence $m=i+j=i+1$, and we have one parameter $c=c_1$. For any $\be=1,\dots,m$ we have: $\g c_\be$ is a monomial in $\sqrt{c}$ (see 5.5.2 as an example) and $\deg \g c_\be\in \frac12 \n Z$ is well-defined. Clearly deg $\im \vf(T,w)$ is less than or equal to $\sum_{\be=1}^m \deg \g c_\be$ (really, $\la_0(c)$ is a monomial in $c$ of degree $\sum_{\be=1}^m \deg \g c_\be$, and other $\la_\vk(c)$ are polynomials in $c$ of degree $\le\sum_{\be=1}^m \deg \g c_\be$. Hence, in homogeneous coordinates $(c:c')$ the map $\vf(F,w)$ is given by polynomials of degree $\sum_{\be=1}^m \deg \g c_\be$. Hence, the degree of its image is $\le\sum_{\be=1}^m \deg \g c_\be$).
\medskip
{\bf Proposition 5.11.1.} For $F=T$ one tree we have $\sum_{\be=1}^m \deg \g c_\be=2^d$.
\medskip
{\bf Proof.} The contribution of a node of depth $k$ to $\sum_{\be=1}^m \deg \g c_\be$ is:
\medskip
$2^{d-k}$ \ \ ( = deg $\left(\sqrt{w_\mu(u)c^{2^{d-k}}}\right)+$ deg $ \left(-\sqrt{w_\mu(u)c^{2^{d-k}}}\right)$ ) for a final node (see (5.4.4), line 3);
\medskip
$2^{d-k-1}$ \ \ ( = deg $w_\mu(u')c^{2^{d-k-1}}$ ) for a non-ramification non-final node (see (5.4.4), line 3);
\medskip
0 for a ramification node (because no number $\g c_*$ is attached to ramification nodes, see (5.4.4)).
\medskip
The sum of these numbers over all nodes of $T$ is $2^d$, this is easily proved by induction by the quantity of branches of $T$. $\square$
\medskip
For further applications (proof of Proposition 5.12.1) we shall need the exact value of $\la_0$ (see (5.5.1)) for any forest $F$.
\medskip
{\bf 5.11.1a.} We need a definition. Let us consider a ramification node $u$. All its descendants form two trees denoted by $D_1(u)$, $D_2(u)$. A ramification node is called a final ramification node if one or two of these trees is a simple tree. These simple trees are called the final branches. Example: tree (4.1.2). We have that for it the
nodes 5, 8, 12, 17, 19 are the final ramification nodes, branches $2-1$, \ $4-3$, \ $7$, \ $10-9$, \ $11$, \ $14$, \ $15$, \ $16$ are the final branches (of lengths 2, 2, 1, 2, 1, 1, 1, 1 respectively).
\medskip
{\bf Lemma 5.11.1.1.} $\la_0=(-1)^j \prod_{\al=1}^j c_\al^{2^{d_\al}}$.
\medskip
{\bf Proof.} The above proposition means that $\la_0=\rho \prod_{\al=1}^j c_\al^{2^{d_\al}}$ where a constant $\rho=\rho(F,w)$ is defined as follows (this follows immediately from (5.4.4), (5.5.1)): $\rho=(-1)^m \prod_u (-w(u)) \prod_{u'} (-w(u'))$ where the first product runs over all final nodes of $F$, and the second product runs over all nodes $u'$ which are the right neighbors of non-final non-ramified nodes of $F$ (the factor $(-1)^m$ comes from (5.5.1)). Hence, $\rho=(-1)^m [(-1)^l \prod_u w(u)] [(-1)^{i+j} \prod_{u'} w(u')]$ (because the quantity of non-final non-ramified nodes of $F$ is $\sum_\al (i_\al -l_\al-(l_\al-1))=i-2l+j$). We have $m=i+j$, hence
$\rho=(-1)^l \prod_u w(u) \prod_{u'} w(u')$. Further, for any tree $T$ and any its weight $w$ we have $\rho(T,w)=-1$. This is proved by induction by the $l$. Really, this is obvious for $l=1$ --- a simple tree. Let us make the induction step. Namely, let $T$ be any tree, $u$ a final ramification node of $T$, $D_1(u) = \{u_1, \dots, u_\vk\}$ a final branch, and $T':= T-D_1(u)$. Hence, $l(T')=l-1$. We denote by $w'$ the restriction of $w$ from $T$ to $T'$, and by $u'_1$ the descendent of $u$ in $T'$. We have $$\frac{\rho(T,w)}{\rho(T',w')}=\frac{(-1)^l(\prod_{\al=2}^\vk w(u_\al))\cdot w(u_\vk)}{(-1)^{l-1}w(u'_1)}=1$$ because $\prod_{\al=2}^\vk w(u_\al))\cdot w(u_\vk)=w(u_1)$ which is $-w(u'_1)$. Hence we get the desired. The proof that for a forest $F$ we have $\rho(F,w)=(-1)^j$ follows immediately. $\square$
\medskip
Let us return to the case of one tree $F=T$; recall that for this case we have $m=i+1$. Let us evaluate the degree of the covering $\vf(T,w): \ \n P^1 \to \im \vf(T,w)\subset \n P^m(\n C)$. We need the following
\medskip
{\bf Definition 5.11.2.} Let $T'$ be a tree. A tree $T$ obtained from $T'$ by replacing some final nodes of the maximal depth of $T'$ by complete trees (of some depths, probably different for different nodes) is called an extension of $T'$, and $T'$ is called a contraction of $T$. Number $\ga:=d(T)-d(T')$ is called the depth of the contraction.
\medskip
{\bf Example 5.11.3.} \ \ \ $T: \ \ \circ^{\nearrow^{\circ^{\nearrow^{\circ<^\circ_\circ}}_{\searrow_{\circ<^\circ_\circ}}    }} _{\searrow _{\circ-\circ^{\nearrow^{\circ<^\circ_\circ}}_{\searrow_{\circ<^\circ_\circ}} }}$, \ \ \ $T': \ \ \circ^{\nearrow^{\circ^{\nearrow^{\circ}}_{\searrow_{\circ}}    }} _{\searrow_{_{\circ -\circ }}}$, \ \ \ $\ga=2$.
\medskip
{\bf Proposition 5.11.4.} Let $T$ have a contraction $T'$ of depth $\ga$. Then $\la_*$ from (5.5.1) are polynomials in $c^{2^\ga}$.
\medskip
{\bf Proof.} Non-ramified nodes of $T$ either have depth $< d-\ga$, or --- for some $\ga'$ satisfying $0\le \ga'\le \ga$ --- are final nodes of complete subtrees of depth $\ga'$ that replace some final nodes of depth $d-\ga$ of $T'$. If $u$ is a non-ramified non-final node of $T$ of depth $< d-\ga$ then the factor $\g Y-w(u)c^{2^{d-k-1}}$ of $\g P(\g Y)$ (Remark 5.5.1.1) corresponding to $u$ contains $c$ in a power which is a multiple of $2^\ga$. If $u$ is a final node of $T$ of depth $< d-\ga$ then the same holds for the factor $(\g Y-\sqrt{w(u)c^{2^{d-k}}})(\g Y+\sqrt{w(u)c^{2^{d-k}}})$ of $\g P(\g Y)$ corresponding to $u$. Let now $\tilde u$ be a final node of depth $d-\ga$ of $T'$ which is replaced by a complete tree $CT_{\tilde u}$ of depth $\ga'$ in $T$. We have ($u$ runs over the set of final nodes of $CT_{\tilde u}$)
$$\prod_{u} (\g Y-\sqrt{w(u)c^{2^{\ga-\ga'}}}\ )\ (\g Y+\sqrt{w(u)c^{2^{\ga-\ga'}}}\ )=\g Y^{2^{\ga'+1}}-w(\tilde u)c^{2^{\ga}}\eqno{(5.11.4.1)}$$
(this formula holds, because the set $\{w(u)\}$, $u$ runs over the set of final nodes of $CT_{\tilde u}$, coincides with the set of roots $w(\tilde u)^{2^{-\ga'}}$. Hence, the set of the roots to \{left hand side of (5.11.4.1) $=0$, $\g Y$ is an unknown\} is the set $w(\tilde u)^{2^{-\ga'-1}}c^{2^{\ga-\ga'-1}}$. Obviously the set of the roots to \{right hand side of (5.11.4.1) $=0$, $\g Y$ is an unknown\} is the same set).
\medskip
We see that $c$ enters in all factors of $\g P(\g Y)$ with a degree which is a multiple of $2^\ga$, hence the proposition. $\square$
\medskip
Let us recall definitions from algebraic geometry. Let $\vf: \n P^1 \to \n P^m$ be a non-constant map. It defines a surjection $\n P^1 \twoheadrightarrow $ im $\vf$. Its degree (the dimension as a vector space over $\n C$ of $O(\vf^{-1}(t))$, where $t$ a generic point of im $\vf$) is called the degree of the covering $\vf$. We have also the degree of the curve im $\vf$ (considered as a subvariety of $\n P^m$). Let $\vf$ be defined in coordinates by the formulas $\vf(x:y)=(P_0(x,y):\dots P_m(x,y))$ where $P_0,\dots, P_m$ are homogeneous polynomials of degree $\vk$. We have
$$\hbox{ (degree of the covering $\vf)\ \cdot $ (degree of im }\vf)=\vk$$
We apply this formula to the case $\vf=\vf(T,w)$.
\medskip
{\bf Corollary 5.11.5.} If $T$ has a contraction of depth $\ga$ then the degree of the covering $\vf(T,w)$ is $\ge 2^{\ga}$, and hence deg im $\vf(T,w)$ is $\le 2^{d-\ga}$.
\medskip
{\bf 5.11.5a.} Let $T$ be fixed. Its contraction $T'$ such that the depth of $T'$ is the minimal possible (and hence the depth of contraction is the maximal possible) is called the minimal contraction of $T$.
\medskip
{\bf Conjecture 5.11.6.} Let $\ga$ be the depth of the minimal contraction of $T$. Then the degree of the covering $\vf$ is $2^{\ga}$, and hence deg im $\vf(T,w)$ is $2^{d-\ga}$.
\medskip
Let us justify this conjecture. We denote $c^{2^\ga}$ by $\bar c$. Proposition 5.11.4 means that $\la_*$ are polynomials in $\bar c$. We can expect that the map $\bar c \mapsto (\la_0:...:\la_m)$ is an inclusion at a generic point. This is exactly 5.11.6.
\medskip
{\bf Corollary 5.11.7.}  Irreducible components of degree 1 (having $j=1$) correspond to complete trees\footnotemark \footnotetext{Really, degree 1 means that $d=\ga$, i.e. the contracted tree has depth 0. This means that the contracted tree is the one-node-tree, hence $T$ is a complete tree.}. See (6.4), (11.3) and A7 for explicit examples.
\medskip
Below, until the end of 5.11, we shall assume truth of (5.11.6). We denote the (conjectural) degree of im $\vf(T,w)$ by $\de=2^{d-\ga}$. Let $\g T: \n P^1\to \n P^1$ be the map defined by $\g T(c:{c'})=(c^{2^{\ga}}:{c'}^{2^{\ga}})$ and let $\chi: \n P^1\to \n P^{\de}$ be the Veronese map (symmetric product). Recall that it is defined by the formula $$\chi (c:{c'})=(c^\de: c^{\de-1}{c'}:...: c{c'}^{\de-1}:{c'}^\de)$$

It is well-known that there exists a linear rational map $\omega: \n P^\de\to \n P^{i+1}$ such that $$\vf(T,w)=\omega\circ \chi\circ\g T\eqno{(5.11.7a)}$$ Really, let  $\vf(T,w)(c:{c'})=(\la_0(c,{c'}):...:\la_{i+1}(c,{c'}))$ where for $\g j=0,\dots, i+1$ we have $\la_\g j(c,{c'})=\sum_{K=0}^{\de} z_{K\g j}c^{2^{\ga}K}{c'}^{2^{\ga}(\de-K)}$ (compare with (5.0.3)). (5.11.7a) holds if we define
$$\om(\tilde\la_0:\tilde\la_1:...:\tilde\la_\de)=(\sum_{K=0}^{\de} z_{K0}\tilde\la_{\de-K}:\sum_{K=0}^{\de} z_{K1}\tilde\la_{\de-K}: ... :\sum_{K=0}^{\de} z_{K,i+1}\tilde\la_{\de-K})$$
\medskip
{\bf 5.11.8.} Let us consider the case of arbitrary $F$. The above notations will bear a subscript $\al$ for a tree $T_\al$. Particularly, we have maps $\omega_\al\circ \chi_\al\circ\g T_\al: \n P^1 \to \n P^{i_\al+1}$ who are equal to $\vf(T_\al,w_\al)$. We consider the Segre map
$$\psi: \n P^{i_1+1}\times \n P^{i_2+1} \times ... \times \n P^{i_j+1} \to \n P^{\g N-1} \hbox{ where }\g N=\prod_\al(i_\al+2)$$ Coordinates in $\n P^{\g N-1}$ are multiindices $\vk_1,\dots,\vk_j$ where $\vk_\al\in [0,\dots,i_\al+1]$. Further, we consider a linear rational map $\g w: \n P^{\g N-1} \to \n P^m$ (recall that $m=i+j=\sum_{\al=0}^j (i_\al+1)$ ) defined by a $\g N\times(m+1)$-matrix whose $\{(\vk_1,\dots,\vk_j),\vk\}$-th entry (here $(\vk_1,\dots,\vk_j)$ is the number of the coordinate of $\n P^{\g N-1}$, and $\vk\in[0,\dots,m]$ is the number of the coordinate of $\n P^m$) is 1 if $\vk=\vk_1+...+\vk_j$ and 0 otherwise.
\medskip
Further, we consider a set-theoretic product\footnotemark \footnotetext{Notation aleph is used here because of lack of symbols. It has nothing common with cardinality of sets.}
\medskip
$\aleph:=(\omega_1\circ\chi_1\circ\g T_1) \times (\omega_2\circ\chi_2\circ\g T_2) \times ... (\omega_j\circ\times \chi_j\circ\g T_j): \n P^1\times \n P^1 \times ... \times \n P^1 \to \n P^{i_1+1}\times \n P^{i_2+1} \times ... \times \n P^{i_j+1}$,
\medskip
We have: $$\hbox{ the map $\g w\circ\psi\circ\aleph: \n P^1\times \n P^1 \times ... \times \n P^1 \to \n P^m$ is $\vf(F,w)$ of (5.0)} \eqno{(5.11.8.1)}$$
This follows immediately from the formula for symmetric polynomials of groups of variables. Namely, let $\tilde x_{\al}= \{\tilde x_{\al1}, \tilde x_{\al2},\dots, \tilde x_{\al,i_\al}\}$ be a $\al$-th set of variables and $\tilde x=\tilde x_{1}\cup \tilde x_{2}\cup ... \cup\tilde x_{j}$ their union. Then $\si_k(\tilde x)=\sum_K \si_{\vk_1}(\tilde x_1)\cdot \si_{\vk_2}(\tilde x_2)\cdot ... \cdot  \si_{\vk_j}(\tilde x_j) $ where $K=(\vk_1,\dots,\vk_j)$ is a multiindex satisfying $\vk_1+...+\vk_j=k$. (5.11.8.1) follows from this formula and from the formula (5.5.1) (see (5.4.4) for its notations).
\medskip
Let us calculate deg im $\vf(F,w)$. For any Segre map $$s: \n P^{n_1}\times \n P^{n_2}\times ...\times \n P^{n_j}\to \n P^N$$ we denote by $L_\al$ the class of $\n P^{n_1}\times ...\times \n P^{n_{\al-1}}\times H_\al\times \n P^{n_{\al+1}} ...\times \n P^{n_j}$ in $CH(\n P^{n_1}\times ...\times \n P^{n_j})$, where $H_\al$ is a hyperplane in $\n P^{n_\al}$. Let $L$ be the class of the hyperplane of $\n P^N$. The Chow ring of $\n P^{n_1}\times \n P^{n_2}\times ...\times \n P^{n_j}$ is $\n Z$-generated by $L_1^{\varkappa_1}\cdot L_2^{\varkappa_2}\cdot ...\cdot L_j^{\varkappa_j}$, where $\varkappa_\al\le n_\al$ (because the Chow ring of $\n P^\vk$ is $\n Z[L]/L^{\vk+1}$). Let $CH(s)_*$ be the direct image map of Chow groups. We have:
$$CH(s)_*(L_1^{\varkappa_1}\cdot L_2^{\varkappa_2}\cdot ...\cdot L_j^{\varkappa_j})=\frac{(\sum_{\al=1}^j (n_\al-\varkappa_\al))!} {\prod_{\al=1}^j (n_\al-\varkappa_\al)!}\ \ L^{N-(\sum_{\al=1}^j (n_\al-\varkappa_\al))}\eqno{(5.11.9)}$$
\{ Let us recall the proof of (5.11.9). First, we have $s^{-1}(L)=\sum_{\al=1}^j L_\al$. Since $s^{-1}:CH(\n P^N)\to CH(\n P^{n_1}\times \n P^{n_2}\times ...\times \n P^{n_j})$ is a ring homomorphism, $\forall \ \vk$ we have $s^{-1}(L^\vk)=(\sum_{\al=1}^j L_\al)^\vk$. We take $\vk=\sum_{\al=1}^j n_\al$. Since $L_\al^{n_\al+1}=0$ in $CH(\n P^{n_\al})$, we get that the only non-zero term in $(\sum_{\al=1}^j L_\al)^\vk$ for this $\vk$ is $\frac{(\sum_{\al=1}^j n_\al)!} {\prod_{\al=1}^j (n_\al)!} \prod_{\al=1}^j L_\al^{n_\al}$, i.e. $s^{-1}(L^{\sum_{\al=1}^j L_\al})$ is the class of $\frac{(\sum_{\al=1}^j n_\al)!} {\prod_{\al=1}^j (n_\al)!}$ points in $\n P^{n_1}\times \n P^{n_2}\times ...\times \n P^{n_j}$. Hence, the same multinomial coefficient is the degree of the Segre map. Taking into consideration that a representative of $L_1^{\varkappa_1}\cdot L_2^{\varkappa_2}\cdot ...\cdot L_j^{\varkappa_j}$ in $CH(\n P^{n_1}\times \n P^{n_2}\times ...\times \n P^{n_j})$ is the product $\n P^{n_1-\varkappa_1}\times \n P^{n_2-\varkappa_2}\times ...\times \n P^{n_j-\varkappa_j}$ and applying the formula of the degree of the Segre map to this product, we get immediately (5.11.9) \}.
\medskip
For our situation we have $n_\al=i_\al+1$, $\varkappa_\al=i_\al$, i.e. all $n_\al-\varkappa_\al$ are 1, hence the multinomial coefficient of (5.11.9) is $j!$, hence the degree of the map $\om\circ\psi\circ\aleph$ is $\prod_{\al=1}^j \de_\al \cdot j!$ (recall that the degree of a map from a variety to another irreducible variety of the same dimension is the quantity of elements in the preimage of a generic point).
\medskip
Isomorphisms between $(T_\al,w_\al)$ and $(T_{\al'},w_{\al'})$ imply that $\om\circ\psi\circ\aleph$ is a finitely-sheeted covering, hence deg im $\vf(F,w)$ is a divisor of $\prod_{\al=1}^j \de_\al \cdot j!$ Namely, let us consider an equivalence relation on the set $T_1,\dots,T_j$: two elements $T_\al$, $T_{\al'}$ are equivalent iff there exists an isomorphism $\bh: T_{\al}\to T_{\al'}$ and $\g g\in \g G_\al$ such that $w_\al=\g g(w_{\al'}\circ\bh)$ (it is easy to check that this is really an equivalence relation).
\medskip
Let $j'$ be the quantity of classes in the above equivalence relation, and $f_1,\dots,f_{j'}$ quantities of elements in these equivalence classes (therefore, $f_1+...+f_{j'}=j$). We get that the degree of the covering $\om\circ\psi\circ\aleph: (\n P^1)^j \to \im \om\circ\psi\circ\aleph$ is $\ge f_1!\cdot...\cdot f_{j'}!$
\medskip
{\bf Conjecture 5.11.10.} This degree is equal to $f_1!\cdot...\cdot f_{j'}!$, hence (taking into consideration Conjecture 5.11.6) we get (notations of (5.10.1)): $$\g d(\bar \vf_{ij}\circ\g p_{ij}(F,w))=\deg \im \vf(F,w)=\frac{2^{\sum_{\al=1}^j d_\al} \cdot j!}{2^{\sum_{\al=1}^j \ga_\al}\cdot f_1!\cdot...\cdot f_{j'}!}\eqno{(5.11.11)}$$
(Recall the notations: $\bar \vf_{ij}\circ\g p_{ij}(F,w)$ is an (abstract) element of $Irr(i,j)$, see (5.10.1). For $C\in Irr(i,j)$ (i.e. $C$ is an element) the number $\g d(C)$ is the degree of $\be_{i+j,i,j}(C)$ --- the minimal irreducible component corresponding to $C$, see (2.4.2.1), case $m=i+j$).
\medskip
All entries of Table A2.2 satisfy this conjecture. See (6.1.4), (6.2.2), (6.3.3), (6.4.1) for explicit examples.
\medskip
{\bf 5.12.} Recall that in Section 3 we considered relations between $X(m,m-2)$ and its hyperplane section $X(m,m-1)$, see formula (3.1). Now we show that the correspondence $Y_i$ of (3.2) has a natural interpretation in terms of weighted forests. Roughly speaking, the inverse map $Y_i^{-1}$ corresponds to elimination of the root, see (4.6). This fact will be used for a deduction of the multiplication formula (5.13.3) from (5.13.2).
\medskip
Let us give details. We denote epimorphisms $\bar \vf_{ij}\circ\g p_{ij}:  \g F_{ij}\to Irr(i,j)$ (see (5.10.1)) by $\tilde \vf_{ij}$. Further, we define a map $\g Z_i: \g F_{i+1,1}\to \g F_{i1}\sqcup \g F_{i2}$ as follows. Let $r$ be the root of $T$. We use notations of (4.6).
\medskip
(a) $r$ is a non-ramification node. We have $\g Z_i(T,w)=(\bar T, \tau (w))\in \g F_{i1}$.
\medskip
(b) $r$ is a ramification node. We have $\g Z_i(T,w)=(D_1(r) \sqcup D_2(r), \tau(w))\in \g F_{i2}$.
\medskip
{\bf Proposition 5.12.1.} Conjectures 0.2.4, 3.3, 5.10 imply that the diagram
$$\matrix \g Z_i: &\g F_{i+1,1}&\to &\g F_{i1}&\sqcup &\g F_{i2}\\ \\ &\tilde \vf_{i+1,1}\downarrow&&\tilde \vf_{i1}\downarrow&&\tilde \vf_{i2}\downarrow\\ \\ Y_i^{-1}: & Irr(i+1,1)&\to &Irr(i,1)&\sqcup& Irr(i,2)\endmatrix$$
is commutative.
\medskip
{\bf Remark.} Without assuming truth of Conjecture 3.3, we can understand commutativity of the above diagram as follows: for any $(T,w)\in \g F_{i+1,1}$ we have: $\tilde \vf_{i+1,1}(T,w)\in Y_i(\tilde \vf_{i\al}(\g Z_i(T,w)))$, where $\al =1$ if $\g Z_i(T,w)\in \g F_{i1}$, and $\al =2$ if $\g Z_i(T,w)\in \g F_{i2}$.
\medskip
{\bf Proof.} Let us consider first the case "$r$ is a non-ramification node". Namely, let $(T_1,w_1)\in \g F_{i1}$, i.e. $T_1$ is a tree having $i$ nodes. Let $d$ be its depth. We denote by $(T,w)\in \g F_{i+1,1}$ a tree and weight such that $(T_1,w_1)=(\bar T, \tau (w))$, i.e. the root of $T$ is a non-ramification node, $T_1$ is obtained from $T$ by elimination of the root, and $w_1$ is the restriction of $w$ to $T_1$ ($w$ is written multiplicatively). We denote $\tilde \vf_{i1}(T_1,w_1)$ by $C_{i1\g k}$ (notations of Section 2) where $\g k$ is the number (label) of the element. Hence, im $\vf(T_1,w_1)=\bar C_{i1\g k}$. We must find $Y_i(C_{i1\g k})$. First, we construct an irreducible component $C_{i1\g k}(i+2)$ of $X(i+2,i)$ that belongs to the series of $C_{i1\g k}$. This construction (depending on Conjecture 0.2.4) is given in Proposition 7.1 (its proof does not depend on the present considerations, i.e. there is no vicious circle). Let $(c:1)\in \n P^1$ be a parameter; we let $\vf(T_1,w_1)(c)=(\la_0:...:\la_{i+1})$. According (7.1.1), $C_{i1\g k}(i+2)$ is a surface whose parametric equations are $$\varpi_\varkappa=\la_{\varkappa-1}b_1+\la_\varkappa b_0, \ \ \varkappa=0,\dots,i+2\eqno{(5.12.1.1)}$$ where $\varpi_0,\dots, \varpi_{i+2}$ are projective coordinates of a point of $C_{i1\g k}(i+2) \subset \n P^{i+2}$ corresponding to parameters $(c:1)\in \n P^1$, $(b_0:b_1)\in \n P^1$ and $\la_{-1}=\la_{i+2}=0$. Further, by definition of $Y_i(C_{i1\g k})$, the set $Y_i(\tilde \vf_{i1}(T_1,w_1))$ is the set of irreducible components of the intersection of $C_{i1\g k}(i+2)$ with the trace hyperplane $$D(i+2,i)=0, \ \ \hbox {i.e. } \ \sum_{\varkappa=1}^{i+1}\varpi_\varkappa=0\eqno{(5.12.1.2)}$$ Since $\bar C_{i1\g k}\subset X(i+1,i)$ we get that $\sum_{\varkappa=1}^{i}\la_\varkappa=0$ (because $X(i+1,i)$ is contained in the hyperplane $D(i+1,i-1)=0$). Hence, (5.12.1.2) becomes $\la_0b_1+\la_{i+1}b_0=0$, i.e. (because $\la_{i+1}=1$, $\la_{0}=-c^{2^d}$, see Lemma 5.11.1.1) in order to satisfy (5.12.1.2), numbers $b_1,b_0$ must be $b_0= c^{2^d}$, $b_1=1$. Hence, (5.12.1.1) becomes $$\varpi_\varkappa=\la_{\varkappa-1}+\la_\varkappa c^{2^d}, \ \ \varkappa=0,\dots,i+2\eqno{(5.12.1.3)}$$

{\bf 5.12.1.3a.} We need more notations (they will be used in formulas (5.12.1.4), (5.12.1.5) and below, throughout the whole Section 5.12). Let $\g c_*=(\g c_1,\dots,\g c_m)$ be a sequence from (5.4.4). The pair of numbers of this sequence $\pm \sqrt{w_\mu(u)c_\al^{2^{d_\al-k}}}$ (reps. the number $-w_\mu(u')c_\al^{2^{d_\al-k-1}}$) corresponding to a final (resp. to a non-ramification non-final) node $u$ is denoted by $\g c(u,F,w,c_*)$.
\medskip
Equations (5.12.1.3) coincide with the equations defining im $\vf(T,w)$. Really, let $u_1$ be a node of $T_1$ and $u$ the corresponding non-root node of $T$. We have $$\g c(u_1,T_1,w_1,c)=\g c(u,T,w,c)\eqno{(5.12.1.4)}$$ Let $r$ be the root of $T$. We have $$\g c(r,T,w,c)=-c^{2^d}\eqno{(5.12.1.5)}$$ Clearly $\rho(T,w)=\rho(T_1,w_1)$ ($\rho$ from the proof of 5.11.1.1). Further,
we have ($\g P$ from 5.5.1.1) $\g P_{T_1,w_1}(\g Y)=\sum_{\varkappa=0}^{i+1}\la_\varkappa \g Y^\varkappa$. (5.12.1.3) implies that $\sum_{\varkappa=0}^{i+2} \varpi_\varkappa \g Y^\varkappa=\g P_{T_1,w_1}(\g Y)\cdot (\g Y+c^{2^d})$. Formulas (5.12.1.4), (5.12.1.5) imply that it is equal to $\g P_{T,w}(\g Y)$. This means exactly that $Y_i(C_{i1\g k})=\vf_{i+1,i}(T,w)$, which proves the proposition for the case "$r$ is a non-ramification node".
\medskip
Now let us consider the case "$r$ is a ramification node". Let $(F,w)\in \g F_{i2}$. We have $F=T_1\sqcup T_2$, $w=(w_1,w_2)$ where $w_\al$ is a weight of $T_\al$ (here and below $\al=1,2$). We have $i_1+i_2=i$. Let $(c_\al:1)\in \n P^1$ be parameters. We denote $\g P_\al(\g Y)=\g P_{T_\al,w_\al}(\g Y)=\sum_{\varkappa=0}^{i_\al+1}\la_{\al,\varkappa}\g Y^\varkappa \in \n Z[c_\al,\g Y]$, $\la_{\al,\varkappa}\in \n Z[c_\al]$. We have $\g P_{(F,w)}(\g Y)=\g P_{1}(\g Y)\g P_{2}(\g Y)$. We denote $\g P_{(F,w)}(\g Y)=\sum_{\varkappa=0}^{i_1+i_2+2}\varpi_\varkappa \g Y^\varkappa$. We have $(\varpi_0:...:\varpi_{i_1+i_2+2})\in \n P^{i_1+i_2+2}$ are parametric equations of im $\vf(F,w)$, where $c_1,c_2$ are parameters. By definition of $Y_i$, we have: $Y_i(\tilde \vf_{i2}(F,w))$ is the set of irreducible components of its intersection with the trace hyperplane $$H:=\{\sum_{\varkappa=1}^{i_1+i_2+1}\varpi_\varkappa=0\}\eqno{(5.12.1.6)}$$ Since for $\al=1,2$ we have $\sum_{\varkappa=1}^{i_\al}\la_{\al,\varkappa}=0$ and $\la_{\al,0}=-c_\al^{2^{d_\al}}$, (5.12.1.6) becomes $c_1^{2^{d_1}}=-c_2^{2^{d_2}}$.
\medskip
Let $d_1\le d_2$. We get: if a point of $\im \vf(F,w)$ belongs to an irreducible component that belongs to $Y_i(\tilde \vf_{i2}(F,w))$ then $$c_1=\g z_{d_1+1}^\de c_2^{2^{d_2-d_1}}\eqno{(5.12.1.7)}$$ where $\de\in \n Z/2^{d_1+1}$ is odd. Different values of $\de$ correspond to different irreducible components of $H\cap \im \vf(F,w)$ (probably some of them can coincide).
\medskip
We denote by $T$ a tree whose root $r$ is ramified and $D_\al(r)=T_\al$. Let $r_\al$ be the root of $T_\al$, i.e. $r_1,r_2$ are descendants of $r$ in $T$. Let $w(\de)$ be a weight of $T$ such that $$w(\de)_a(r_1)=1\hbox{ and } \tau(w(\de))=((\frac{\de-1}2)w_1,w_2)\eqno{(5.12.1.8)}$$ where $\frac{\de-1}2\in \n Z/2^{d_1}$ and the action of $\frac{\de-1}2$ on $w_1$ is the action of $\g G$ on $W(T_1)$ defined by (4.5.1).
\medskip
Let $u$ be a non-ramification node of $F$. We can consider $u$ as a node of $T$ as well.
\medskip
{\bf Lemma 5.12.1.9.} Let $c_1$ satisfies (5.12.1.7). Then $\g c(u,F,w,(c_1,c_2))=\g c(u,T,w(\de),c_2)$ (notations of (5.12.1.3a)).
\medskip
{\bf Proof.} For $u\in T_2$ this is clear. Let us consider the case of $u\in T_1$ is a final node of depth $k$. We have $\g c(u,F,w,(c_1,c_2))=\pm \sqrt{w_1(u)c_1^{2^{d_1-k}}}$ and $\g c(u,T,w(\de),c_2)=\pm \sqrt{w(\de)(u)c_2^{2^{(d_2+1)-(k+1)}}}$, because the depth of $T$ is $d_2+1$, the depth of $u$ as a node of $T$ is $k+1$. Further, (4.6.1) and (5.12.1.6) imply that $\frac{w(\de)_a(u)-1}2=\frac{\de-1}2+{w_1}_a(u)$ (equality in additive form, in $\n Z/2^k$), i.e. $w(\de)(u)=\g z_{k+1}^\de w_1(u)$ (multiplicative form). Substituting this value and the value of $c_1$ from (5.12.1.7) to these formulas we get the desired.

For $u\in T_1$ a non-ramification non-final node the proof is the same (we use its right neighbor $u'$ instead of $u$ in the above calculation). $\square$
\medskip
Proposition 5.12.1 follows immediately from this lemma, because $\tilde \vf_{i+1,1}(T,w(1))\in Y_i(\tilde \vf_{i2}(F,w))$ and all elements of $\g F_{i+1,1}$ are of the form $(T,w(1))$ for some $(F,w)\in \g F_{i2}$. $\square$
\medskip
{\bf Remark 5.12.2.} Conjecture 5.10 gives evidence that the converse correspondence $Y_i^{-1}: Irr(i+1,1)\to Irr(i,1)\sqcup Irr(i,2)$ is a function (considered as a particular case of a correspondence). Really, according Conjecture 5.10, a fiber of $\tilde \vf_{i+1,1}$ is a $\Aut(T)\times \g G$-orbit. Description of $\g Z_i$ shows that $\g Z_i$ of this orbit is contained in an orbit of the analogous group for $(\bar T, \tau(w))$ or $(D_1(r) \sqcup D_2(r), \tau(w))$. Hence, its $\tilde \vf_{i1}$ or $\tilde \vf_{i2}$-image is one point.

This is not a proof of the fact that Conjecture 5.10 implies Conjecture 3.3, because Proposition 5.12 depends on Conjecture 3.3. To get rid of a vicious circle, we need more work; this is a subject of further research.
\medskip
{\bf Remark 5.12.3.} (4.6.3) shows why $Y_i$ is not 1 -- 1, i.e. $Y_i^{-1}$ is not injective. We have $\g G(F)=\g G(D_1(r))\times \g G(D_2(r))$. $(\g g_1,\g g_2)$ of (4.6.3) belongs to $\g G(F)$, and $\tilde \vf_{i2}(F,\tau(w))=\tilde \vf_{i2}(F,(\g g_1,\g g_2)(\tau(w)))$. The fiber $\tau^{-1}(\tau w)$ consists of two elements, and 4 elements $\tau^{-1}(\tau w)$, $\tau^{-1}((\g g_1,\g g_2)(\tau w))$ usually are not $\g G(T)$-equivalent.
\medskip
Example for $i=5$, see Table A2.2: $Y_5^{-1}(C_{618a})=Y_5^{-1}(C_{618b})=C_{522}$. This is the only example for $i\le5$ of non-injectivity of $Y_i^{-1}$. We have $C_{618a}$, $C_{618b}$ are defined over $\n Q[\sqrt{-1}]$, they are $\n Q$-conjugate (i.e. $C_{618b}=\s(C_{618a})$ where $\s$ is the complex conjugation. This follows easily from the results of Section 5.15).
\medskip
{\bf Remark 5.12.4.} The subject of the present subsection is to find a relation between $Irr(i,1)\sqcup Irr(i,2)$ and $Irr(i+1,1)$ based on the fact that $X(m,m-1)$ is a hyperplane section of $X(m,m-2)$. Analogously, the intersection of $X(i+j,i)$ and hypersurface $D(i+j,i)=0$ gives us a correspondence from $\sqcup_{\vk=1}^j Irr(i,\vk)$ to $\sqcup_{\vk=1}^{j-1} Irr(i+1,\vk)$. An analog of (5.12.1) for this situation is a subject of future research.
\medskip
{\bf 5.13. Multiplicity formula.} Let $u$ be a ramification node of $T$. We denote by $q_\al(u)$ the quantity of nodes of $D_\al(u)$ ($\al=1$, 2).
\medskip
{\bf Conjecture 5.13.1.} The multiplicity $\mu(\im \vf(F,w))=\mu(F)$ of $\im \vf(F,w)$ depends only of $F$. It is given by the formulas:
$$\mu(F)=\frac{i!}{i_1!\cdot i_2!\cdot...\cdot i_j!}\prod_{\al=1}^j \mu(T_\al)\eqno{(5.13.2)}$$ $$\mu(T)=\prod_u\binom{q_1(u)+q_2(u)}{q_1(u)}\eqno{(5.13.3)}$$ where $u$ runs over the set of all ramification nodes of $T$ (the empty product is 1).
\medskip
At it was mentioned above, we have no theoretical justification of (5.13.2), it is purely experimental, based on computer calculations, see Table A2.2. Explicit examples are given in (6.1.4), (6.2.2), (6.3.3), (6.4.2), (6.4.3).
\medskip
{\bf Example 5.13.4.} Let us calculate $\mu(F)$ for $F$ of (5.5.2 (B)). We have $F=T_1\cup T_2$, ramification nodes of $T_1$ are nodes $u_3$, $u_6$, we have for them $q_1(u_3)=q_2(u_3)=1$, $q_1(u_6)=2, \ q_2(u_6)=3$, hence $\mu(T_1)=\binom21\binom52=20$. Ramification node of $T_2$ is the node $u_{10}$, we have for it $q_1(u_{10})=1, \ q_2(u_{10})=2$, hence $\mu(T_2)=\binom31=3$. We have $i=13$, $i_1=6$, $i_2=5$, hence (5.13.2) gives us $\mu(F)=60\cdot\binom{13}5$.
\medskip
{\bf Deduction of (5.13.3) from (5.13.2) and Conjectures 0.2.4, 5.10.} The below text is only an idea of this deduction. It should be considered as an argument supporting Conjecture 5.13.1.

We use induction: let (5.13.3) hold for all trees having $\le i$ nodes. Let $T$ be a tree having $i+1$ nodes such that its root $r$ is a ramification node, and $w\in W(T)$. We have $\tilde \vf_{i+1,1}(T,w)\in Irr(i+1,1)$. We consider $Y_i^{-1}(\tilde \vf_{i+1,1}(T,w))\in Irr(i,2)$. By 5.12.1, $Y_i^{-1}(\tilde \vf_{i+1,1}(T,w))=\tilde \vf_{i2}(D_1(r)\sqcup D_2(r),\tau(w))$. By induction supposition and by (5.13.2), $$\mu(\im \vf(D_1(r)\sqcup D_2(r),\tau(w)))=\frac{i!}{\#(D_1(r))!\#(D_2(r))!}\cdot \mu(\im \vf(D_1(r),\tau_1(w)))\cdot$$ $$\cdot \mu(\im \vf(D_2(r),\tau_2(w)))$$ Because of (3.4.3), we deduce that (5.13.3) holds for $\im \vf(T,w)$.
\medskip
For the case of $T$ having $i+1$ nodes such that its root $r$ is a non-ramification node, (5.12.1) and (3.4.3) show immediately that truth of (5.13.3) for $\bar T$ implies truth of (5.13.3) for $T$. $\square$
\medskip
{\bf 5.14. Jordan form.} We need more definitions. Let $T$ be a tree, and let $u$ be a node of $T$. We denote by $R_u$ the subtree of $T$ formed by $u$ and all its descendants; $u$ is the root of $R_u$. The depth of $R_u$ is called the height of $u$, it is denoted by $h(u)$.

Example: we consider the tree (4.1.2). At the present picture, a number at a node indicates its height:
$${\underset{7}\to{\bullet}}^{\nearrow^{{\underset{4}\to{\diamond}}^ {\nearrow^{\underset{3}\to{\circ}-{\underset{2}\to{\diamond}}<^{\circ1-\ast0}_{\circ 1-\ast0}}}
_{\searrow_{\ast0}}}}_ {\searrow_{\overset{6}\to{\circ}-{\underset{5}\to{\diamond}}^{\nearrow^{\circ3-\diamond2<^{\circ1-\ast0}_{\ast0}}}_{\searrow_{\underset{4}
\to{\circ}-{{\underset{3}\to{\diamond}}<^{\ast0}_{\overset{2}\to{\circ}-{\underset{1}\to{\diamond}}<^{\ast0}_{\ast0}}}  }    }     }}\eqno{(5.14.0)}$$

Finally, like in the proof of Theorem 5.6, we can consider a matrix $\g M(m)(\la_0,...,\la_m)$ as an operator $\g M$ on $\n C^{m-1}$.
\medskip
{\bf Lemma 5.14.1.} Conjectures 2.1, 5.8 imply: for any $F,w$, for generic $c_1,\dots, c_j$, for $(\la_0,...,\la_m)=\vf(F,w)(c_1,\dots, c_j)$, for any $\varkappa$ we have dim Ker $\g M^\varkappa=$ \{the quantity of nodes $u$ of $F$ such that $h(u)\le \varkappa-1$\}.
\medskip
{\bf Proof.} Let $u$ be a node such that $h(u)<\varkappa$. We have $\g M^\varkappa(\chi_u)=0$ (see the argument (5.6.4) of the proof of Theorem 5.6). For generic $c_1,\dots, c_j$ vectors $\chi_u$ are linearly independent, hence for generic $c_1,\dots, c_j$ we have dim Ker $\g M^\varkappa\ge$ \{the quantity of nodes $u$ of $F$ such that $h(u)\le \varkappa-1$\}. Let us prove that we have an equality. Numbers $\be_u$ from (5.3), $\be_{1,u}$, $\be_{2,u}$ from (5.6.2) are polynomials in $c_1,\dots, c_j$ and hence are $\ne0$ for generic $c_1,\dots, c_j$. This means that if $u$ is a node of height $\varkappa$ then $\g M^\varkappa(\chi_u)\ne0$. Hence, dim (Ker $\g M^\varkappa/ \Ker \g M^{\varkappa-1})\ge$ \{the quantity of nodes $u$ of $F$ such that $h(u)=\varkappa-1$\}. This means that if for some $\varkappa_0$ dim Ker $\g M^{\varkappa_0}/ \Ker \g M^{\varkappa_0-1}>$ \{the quantity of nodes $u$ of $F$ such that $h(u)=\varkappa-1$\} then dim Ker $\g M^{d+1}=i'$ for some $i'>i$, and hence im $\vf(F,w)\subset X(m,i')$. This contradicts to Conjectures 2.1, 5.8 (because according Conjecture 2.1 we have dim $X(m,i')=m-i'$ and according Conjecture 5.8 we have dim im $\vf(F,w)=m-i$. So, im $\vf(F,w)\subset X(m,i')$ cannot happen, because $m-i>m-i'$). $\square$
\medskip
{\bf Remark 5.14.1.1.} The above proof shows that for generic $c_1,\dots, c_j$ we have: dim Ker $\g M^\varkappa$ does not depend on $c_1,\dots, c_j$, i.e. there exists a Zariski open subset of the set of $c_1,\dots, c_j\in (\n P^1)^j$ where dim Ker $\g M^\varkappa$ is constant. Since the set of numbers dim Ker $\g M^\varkappa$ defines a partition (see below), we get that Lemma 5.14.1 supports Conjecture 2.3.1.1.
\medskip
Let $M$ be any 0-Jordan matrix (i.e. a union of Jordan blocks of eigenvalue 0) of size $n$. It defines a partition of $n$ denoted by $\varsigma$, see (2.3.1) for details. There is a relation between the rank of powers of $M$ and the conjugate\footnotemark \footnotetext{Partition obtained by reflection of its Young diagram, see [W1].} of $\varsigma$:
\medskip
{\bf 5.14.2.} If $n=p_1+...+p_k$ is the conjugate of $\varsigma$ then the rank of $M^\vk$ is $n-p_1-...-p_\vk$, i.e. dim Ker $M^\vk$ is $p_1+...+p_\vk$ (here and below, by default, for any partition we assume $p_1\ge...\ge p_k$).
\medskip
{\bf 5.14.3.} We shall need a notion of a partition of a partition. Let $\varsigma: \{i=d_1+...+d_l\}$ be a partition of $i$. A partition of $\varsigma$ is a representation of the set $\{1,...,l\}$ as a disjoint union: $\{1,...,l\}=\sqcup_{\g v=1}^j \Cal Z_\g v$. We denote $\Cal Z_\g v=\{z_{\g v1},\dots,z_{\g v,l_\g v}\}$, we let $i_\g v:=d_{z_{\g v1}}+...+d_{z_{\g v,l_\g v}}$. Hence, we get partitions (here maybe $i_1\ge...\ge i_j$ does not hold):
$$\bar \varsigma: \{i=i_1+...+i_j\};\eqno{(5.14.3.1)}$$ (a "quotient partition" of $\varsigma$) and
$$\varsigma_\g v: \{i_\g v=d_{z_{\g v1}}+...+d_{z_{\g v,l_\g v}}\}\eqno{(5.14.3.2)}$$ (subpartitions of $\varsigma$).

We do not distinguish different unions $\{1,...,l\}=\sqcup_{\g v=1}^j \Cal Z_\g v$ if they give us isomorphic (up to order) $\bar \varsigma$, $\varsigma_\g v$. Namely, let $\{1,...,l\}=\sqcup_{\g v=1}^{j'} \Cal Z'_\g v$ be another representation of the set $\{1,...,l\}$ as a disjoint union. These two unions are equivalent iff $j=j'$ and after some permutation of indices $\g v$ we can get: $\forall$ $\g v$ the sequences $d_{z_{\g v1}}, ... ,d_{z_{\g v,l_\g v}}$ and $d_{z'_{\g v1}}, ... ,d_{z'_{\g v,l_\g v}}$ coincide up to permutation of their terms (particularly, $l_\g v=l'_\g v$). For example, if $\varsigma : \{i=d_1+d_2+d_3\}$ and $d_2=d_3$ then unions $\{1,2,3\}=\{1,2\}\sqcup\{3\}$ and $\{1,2,3\}=\{1,3\}\sqcup\{2\}$ give us the same partition of $\varsigma$.
\medskip
Two extreme partitions corresponding to the cases $j=1$ and $j=l$ are called the upper and lower trivial partitions respectively.
\medskip
(5.14.3.3) We need also the opposite construction. If partitions $\varsigma_\g v$ from (5.14.3.2) are given then $\varsigma$ is called their union.
\medskip
{\bf 5.14.4. From a forest to a partition.} Here we describe two methods to find $\pi_{m,i}(\im \vf(F,w))$ --- the partition of the Jordan form of a matrix $\g M(m)(\la_0,...,\la_m)$ where $(\la_0,...,\la_m)\in \im \vf(F,w)$. \ \footnotemark \footnotetext{Recall that $\pi_{m,i}$ (see 2.3) is a map from the set of irreducible components of $X(m,i)$ to the set of partitions of $i$. We have $\im \vf(F,w)$ is an irreducible component of $X(m,i)$ and $\pi_{m,i}(\im \vf(F,w))$ is the corresponding partition of $i$.} For brevity, we denote the partition $\pi_{m,i}(\im \vf(F,w))$ by $\g s(F)$ (it does not depend on $w$, see below). The next subsection 5.14.5 contains the opposite construction --- from a partition to a forest. This construction is the inverse with respect to the second method.
\medskip
{\bf 5.14.4.1. First method.} (5.14.2) shows that $\g s(F)$ depends only on dim Ker $M^\vk$, while Lemma 5.14.1 shows that dim Ker $M^\vk$ depends only on the quantities of nodes of $F$ of given heights. Let us define by $\g h(\vk)$ the quantity of nodes $u$ of $F$ such that $h(u)\le\vk$. Numbers $\g h(\vk)$ define a partition $$\tilde \pi(F):=i=\g h(0) + [\g h(1)-\g h(0)] +  [\g h(2)-\g h(1)] +... + [\g h(d)-\g h(d-1)]\eqno{(5.14.4.1.1)}$$
Here $\g h(\vk)-\g h(\vk-1)$ is the quantity of nodes $u$ of $F$ such that $\g h(u)=\vk$. Since the left neighbor $u_l$ of a node $u$ satisfies $\g h(u_l)\ge\g h(u)+1$ we get that for any tree we have $$\g h(0) \ge \g h(1)-\g h(0)\ge \g h(2)-\g h(1)\ge ... \ge \g h(d)-\g h(d-1)\eqno{(5.14.4.1.2)}$$ Lemma 5.14.1 and 5.14.2 give us that $\g s(F)$ is the conjugate to $\tilde \pi(F)$. Particularly, it does not depend on $w$.
\medskip
{\bf Example:} For the tree (4.1.2) we have (we use heights of nodes of the tree (4.1.2) given at the picture (5.14.0)):
$$\g h(0)=8, \ \g h(1)=12, \ \g h(2)=15, \ \g h(3)=18, \ \g h(4)=20, \ \g h(5)=21, \ \g h(6)=22, \ \g h(7)=23$$
$$\tilde \pi(F): \ 23=8+4+3+3+2+1+1+1$$ $$\hbox {$\g s(F)$ is $23=8+5+4+2+1+1+1+1$} \eqno{(5.14.4.1.3)}$$

{\bf 5.14.4.2. Second method.} It is equivalent to the first method. As a first step, we consider a case $j=1$, i.e. $F=T$ is a tree. The corresponding partition $\varsigma=\g s(T)$ is described as follows. We denote $T_{(1)}:=T$. First, let $\g y_1$ be the minimal length of the final branches (see (5.11.1a)) of $T_{(1)}$. Further, let $\phi_1$ be the quantity of the final ramification nodes having one or two final branches of length
$\g y_1$.
We write the end of $\varsigma$ as follows: $i=...+...+ ... \ \ \ ...+\underset{\phi_1\hbox{ times }}\to {\underbrace{\g y_1+\g y_1+...+\g y_1}}$.
\medskip
Now, we cut off $\g y_1$ nodes from the end of all final branches of $T_{(1)}$. Particularly, if $u$ is a final ramification node such that both $D_1(u)$, $D_2(u)$ are simple trees of length $\g y_1$ then after this cutting $u$ will become a final node. If $u$ is a final ramification node such that only one of $D_1(u)$, $D_2(u)$ is a simple tree of length $\g y_1$ then after this cutting $u$ will become a simple (non-ramification, non-final) node.
\medskip
We denote the obtained tree by $T_{(2)}$, and we repeat the procedure: we denote by $\g y_2$ the minimal length of the final branches of $T_{(2)}$, we denote by $\phi_2$ the quantity of the final ramification nodes of $T_{(2)}$ having one or two final branches of length $\g y_2$. We write the end of $\varsigma$ as follows:

$i=...+...+ ... \ \ \ ...+\underset{\phi_2\hbox{ times }}\to {\underbrace{(\g y_1+\g y_2)+(\g y_1+\g y_2)+...+(\g y_1+\g y_2)}}+\underset{\phi_1\hbox{ times }}\to {\underbrace{\g y_1+\g y_1+...+\g y_1}}$.
\medskip
Now, we cut off $\g y_2$ nodes from the end of all final branches of $T_{(2)}$, etc.; the end of $\varsigma$ after the third step is written as follows:

$i=...+...+ ... \ \ \ ...+\underset{\phi_3\hbox{ times }}\to {\underbrace{(\g y_1+\g y_2+\g y_3)+(\g y_1+\g y_2+\g y_3)+...+(\g y_1+\g y_2+\g y_3)}}+$

$+\underset{\phi_2\hbox{ times }}\to {\underbrace{(\g y_1+\g y_2)+(\g y_1+\g y_2)+...+(\g y_1+\g y_2)}}
+\underset{\phi_1\hbox{ times }}\to {\underbrace{\g y_1+\g y_1+...+\g y_1}}$.
\medskip
{\bf 5.14.4.2.0.} The value of $\g y_{max}$ is the quantity of nodes in $T_{(max)}$ which is a simple tree, hence the value of $\phi_{max}$ is always 1.

We shall have $\g s(T)= \{i=\varrho_1+...+\varrho_l\}$, where $\varrho_1=d-1=\g y_1+\g y_2+...+\g y_{max}$, $\varrho_l=\g y_1$ and $l$ = the quantity of final nodes.
\medskip
{\bf Example.} For the above tree (4.1.2) we have $\g y_1=1$, $\phi_1=4$ (final ramification nodes having one or two final branches of length 1 are nodes No. 8, 12, 17, 19). The end of $\varsigma$ is $$23=...+1+1+1+1$$ and $T_{(2)}$ is the following:

$${\underset{23}\to{\bullet}}^{\nearrow^{{\underset{8}\to{\circ}}-{\underset{6}\to{\circ}-{\underset{5}\to{\diamond}}<^{\ast2}_{\ast4}}}}_ {\searrow_{\underset{22}\to{\circ}-{\underset{21}\to{\diamond}}^{\nearrow^{\circ13-\circ12-\ast10}}_{\searrow_{\underset{20}
\to{\circ}-{\underset{19}\to{\circ}-\underset{18}\to{\circ}-{\underset{17}\to{\ast}}}  }    }     }}\eqno{(5.14.4.2.1)}$$
We have $\g y_2=1$, $\phi_2=1$ (final ramification node having one or two final branches of length 1 is the node No. 5). The end of $\varsigma$ is $$23=...+2+1+1+1+1$$ and $T_{(3)}$ is the following:
$${\underset{23}\to{\bullet}}^{\nearrow^{{\underset{8}\to{\circ}}-{\underset{6}\to{\circ}-{\underset{5}\to{\ast}}}}}_ {\searrow_{\underset{22}\to{\circ}-{\underset{21}\to{\diamond}}^{\nearrow^{\circ13-\ast12}}_{\searrow_{\underset{20}
\to{\circ}-{{\underset{19}\to{\circ}}-\underset{18}\to{\ast}}}      }     }}\eqno{(5.14.4.2.2)}$$
We have $\g y_3=2$, $\phi_3=1$ (final ramification node having one or two final branches of length 2 is the node No. 21). The end of $\varsigma$ is $$23=...+4+2+1+1+1+1$$ and $T_{(4)}$ is the following:
$${\underset{23}\to{\bullet}}^{\nearrow^{\overset{8}\to{\ast}
}}_ {\searrow_{\underset{22}\to{\circ}-\underset{21}\to{\circ}-\underset{20}\to{\ast}     }}\eqno{(5.14.4.2.3)}$$
We have $\g y_4=1$, $\phi_4=1$, (final ramification node having one or two final branches of length 1 is the node No. 23). The end of $\varsigma$ is $$23=...+5+4+2+1+1+1+1$$ and $T_{(5)}$ is the following (a simple tree --- end of the process): $${\underset{23}\to{\bullet}}-\underset{22}\to{\circ}-\underset{21}\to{\ast}\eqno{(5.14.4.2.4)}$$ We have $\g y_5=3$, $\phi_5=1$, and the partition of the tree (4.1.2) is $$23=8+5+4+2+1+1+1+1$$

Finally, for a forest $F$ we have: $\g s(F)=\sqcup_\al \ \g s(T_\al)$, the union is in terms of (5.14.3.3).
\medskip
{\bf Proposition 5.14.4.3.} The first and the second methods are equivalent.
\medskip
{\bf Proof.} First, according the definition of the conjugate partition, we should show that $$\g h(0)=\g h(1)-\g h(0)=\g h(2)-\g h(1)=...=\g h(\g y_1-1)-\g h(\g y_1-2)\eqno{(5.14.4.3.1)}$$ (here $\g h(-1):=0$) and $$[\g h(\g y_1)-\g h(\g y_1-1)] -  [\g h(\g y_1-1)-\g h(\g y_1-2)]=\phi_1\eqno{(5.14.4.3.2)}$$
(5.14.4.3.1) holds because $\g y_1$ is the minimal length of final branches. (5.14.4.3.2) holds because any final ramification node having one or two final branches of length $\g y_1$ gives a contribution 1 to the difference $[\g h(\g y_1)-\g h(\g y_1-1)] -  [\g h(\g y_1-1)-\g h(\g y_1-2)]$, hence the whole such difference is $\phi_1$.
\medskip
Second, we use induction process. We use notation $F_{(1)}=F$, $F_{(2)}$ is the disjoint union of $T_{\al,(2)}$. We can assume that for the forest $F_{(2)}$ the equivalence is already proved. Let $u$ be a node of $F_{(2)}$. We denote its height as a node of $F_{(2)}$ (resp. as a node of $F_{(1)}$) by $h_{(2)}(u)$ (resp. $h_{(1)}(u)$). The operation of cutting of $\g y_1$ nodes of all final branches of $T$ gives us $h_{(2)}(u)=h_{(1)}(u)-\g y_1$. Hence, the partition $\tilde \pi(F_{(2)})$ (see (5.14.4.1.1)) is obtained from $\tilde \pi(F_{(1)})$ by elimination of the first $\g y_1$ summands (each of them is equal $\g h(0)$ ). If we pass to the conjugate partition, this corresponds to subtraction of $\g y_1$ from all summands of $\g s(T)$. Since by the induction assumption two methods coincide for $F_{(2)}$, they coincide for $F_{(1)}$ as well. $\square$
\medskip
Finally, the above considerations imply
\medskip
{\bf Proposition 5.14.4.4.} Let $F$ be a forest, $w$ its weight, $\im \vf(F,w)$ its irreducible component belonging to $Irr(i,j)$, and $\{i=d_1+...+d_l\}$ its partition. We have: the present $l$ is $l$ of (4.2.1), (4.2.3), and the present $d_1$ is $d+1$ where $d$ is from (4.2.3). $\square$
\medskip
{\bf 5.14.5. From a partition to forests.} Here we consider the process which is inverse to the above Second method. We have seen that the partition does not depend on weight, but only on the forest. Moreover, different forests can give the same partition: $\g s$ is not injective. So, for a given partition $\varsigma$ we describe here in combinatorial terms the whole set of forests $\g s^{-1}(\varsigma)$.

The numbering of $_*\phi$ and of $_{(*)}T$ is opposite with respect to the Section 5.14.4.2, because the process runs in the opposite direction. Hence, we shall use left subscripts. Let $$\varsigma: \{ i=\underset{_1\phi\hbox{ times }}\to {\underbrace{_1\varrho+...+{_1\varrho}}}+\underset{_2\phi\hbox{ times }}\to {\underbrace{_2\varrho+...+{_2\varrho}}}+...+ \underset{_\vk\phi\hbox{ times }}\to {\underbrace{_\vk\varrho+...+{_\vk\varrho}}}\}\eqno{(5.14.5.1)}$$ be a partition of $i$, where $_1\varrho> {_2\varrho}>...> {_\vk\varrho}$, and we let $_{\vk+1}\varrho=0$. First, we describe the subset of $\g s^{-1}(\varsigma)$ formed by trees; we denote it by $\g s^{-1}(\varsigma)_1$.
\medskip
{\bf 5.14.5.1a.} We must have $_1\phi=1$: if $_1\phi\ne1$ then $\g s^{-1}(\varsigma)_1=\emptyset$, see 5.14.4.2.0.
\medskip
The tree $_{(1)}T$ is a simple tree of length $_1\varrho-{_2\varrho}$. Now, we must choose $_2\phi$ nodes of $_{(1)}T$ (this can be made by $\binom{_1\varrho-{_2\varrho}}{_2\phi}$ ways; if $_1\varrho-{_2\varrho}<{_2\phi}$ then $\g s^{-1}(\varsigma)_1=\emptyset$). To get $_{(2)}T$, we join $_2\phi$ simple trees of length $_2\varrho-{_3\varrho}$ to chosen nodes of $_{(1)}T$, and one complementary simple tree of length $_2\varrho-{_3\varrho}$ to the final node of $_{(1)}T$. Particularly, the final node of $_{(1)}T$ becomes a ramification node of $_{(2)}T$ iff it belongs to the set of chosen nodes. The below picture gives examples for $_{(2)}T$ for the case $_1\varrho-{_2\varrho}=3, \ \ _2\phi=2, \ _2\varrho- {_3\varrho}=3$ (nodes of $_{(1)}T$ are denoted by numbers 1 (root node), 2, 3 (final node)).
\medskip
$^{{\underset{1}\to{}}^{\nearrow^{\circ-\circ-\circ}}_{-2_{\searrow_{\circ-\circ-\circ}}^{-3-\circ-\circ-\circ}}}$ (nodes 1, 2 of $_{(1)}T$ are chosen);
\medskip
$^{{\underset{1}\to{}}^{\nearrow^{\circ-\circ-\circ}}_{-2-3<^{\circ-\circ-\circ}_ {\circ-\circ-\circ}}}$ (nodes 1, 3 of $_{(1)}T$ are chosen);
\medskip
$^{{\underset{1-2}\to{}}^{\nearrow^{\circ-\circ-\circ}}_{-3<^{\circ-\circ-\circ}_ {\circ-\circ-\circ}}}$ (nodes 2, 3 of $_{(1)}T$ are chosen).
\medskip
To get $_{(3)}T$, we must choose $_3\phi$ nodes of $_{(2)}T$ which are not ramification nodes. $_{(2)}T$ has $_1\varrho- {_2\varrho}$ nodes of $_{(1)}T$ and $(_2\phi+1)(_2\varrho-{_3\varrho})$ nodes joined while we extended $_{(1)}T$ to $_{(2)}T$. Also, $_{(2)}T$ has $_2\phi$ ramification nodes. Therefore, a choice can be made by $$\binom{(_1\varrho-{_3\varrho})+(_2\varrho-{_3\varrho}-1)_2\phi}{_3\phi}$$ ways; if $(_1\varrho-{_3\varrho})+(_2\varrho-{_3\varrho}-1)_2\phi< {_3\phi}$ then $\g s^{-1}(\varsigma)_1=\emptyset$. To get $_{(3)}T$, we join $_3\phi$ simple trees of length $_3\varrho-{_4\varrho}$ to chosen nodes of $_{(2)}T$, and one complementary simple tree of length $_2\varrho-{_3\varrho}$ to each of the $_2\phi+1$ final nodes of $_{(2)}T$. Like above, a final node of $_{(2)}T$ becomes a ramification node of $_{(3)}T$ iff it belongs to the set of chosen nodes.
\medskip
Continuing the process, we get that $\forall \ \ga$ \  $_{(\ga)}T$ has
$$(_1\varrho-{_{\ga+1}\varrho})_1\phi+(_2\varrho-{_{\ga+1}\varrho})_2\phi+...+(_\ga\varrho-{_{\ga+1}\varrho})_\ga\phi\hbox{ nodes,}$$ $_1\phi+{_2\phi}+...+{_\ga\phi}$ final nodes and $_2\phi+{_3\phi}+...+{_\ga\phi}$ ramification nodes. We let $T={_{(\vk)}T}$. We see that $\#(\g s^{-1}(\varsigma)_1)$ is $\le$
$$\binom{_1\varrho-{_2\varrho}}{_2\phi}\cdot\binom{(_1\varrho-{_3\varrho})+(_2\varrho-{_3\varrho}-1)_2\phi}{_3\phi}\cdot$$ $$\cdot\binom{(_1\varrho-{_4\varrho})+(_2\varrho-{_4\varrho}-1)_2\phi+(_3\varrho-{_4\varrho}-1)_3\phi}{_4\phi}\cdot \ \dots $$ $$\cdot\binom{(_1\varrho-{_\vk\varrho})+(_2\varrho-{_\vk\varrho}-1)_2\phi+(_3\varrho-{_\vk\varrho}-1)_3\phi+... +(_{\vk-1}\varrho-{_\vk\varrho}-1)_{\vk-1}\phi}{_\vk\phi}$$
(it can be less, because there are symmetries --- see the below example 5.14.5.2b).
\medskip
{\bf Example 5.14.5.2a.} Let us consider the partition (5.14.4.1.3), and show how to get the tree (4.1.2) from it. We have $_1\rho=8$, $_2\rho=5$, hence $_{(1)}T$ is a simple tree of length 3 --- this is the tree (5.14.4.2.4) ($=T_5$). We have $_2\phi=1$, hence we choose 1 node of the tree (5.14.4.2.4). For our case, this is the node 23 (the root). We have $_3\rho=4$, hence we join one simple tree of length $_2\rho-{_3\rho}=1$ to the node 23, and one simple tree of the same length 1 to the final node 21. We get the tree (5.14.4.2.3), it is $_{(2)}T=T_4$.
\medskip
Now we have $_3\phi=1$, hence we choose 1 node of the tree (5.14.4.2.3). For our case, this is the node 21. We have $_4\rho=2$, hence we join one simple tree of length $_3\rho-{_4\rho}=2$ to the chosen node 21, and one simple tree of the same length 2 to the final nodes 20 and 8. We get the tree (5.14.4.2.2), it is $_{(3)}T=T_3$.
\medskip
Next step of the process: we have $_4\phi=1$, hence we choose 1 node of the tree (5.14.4.2.2). For our case, this is the node 5. We have $_5\rho=1$, hence we join two simple trees of length $_4\rho-{_5\rho}=1$ to the chosen node 5 (because it is a final node), and one simple tree of the same length 1 to other final nodes 12 and 18. We get the tree (5.14.4.2.1), it is $_{(4)}T=T_2$.
\medskip
The final step of the process: we have $_5\phi=4$, hence we choose 4 nodes of the tree (5.14.4.2.1). For our case, these are the nodes 8, 12, 19, 17, where node 17 is a final node. We have $_6\rho=0$ (by definition), hence we join two simple trees of length $_5\rho-{_6\rho}=1$ to the final node 17, one simple tree of the same length 1 to other chosen nodes 8, 12, 19, and one simple tree of the same length 1 to other final nodes 2, 4, 10. We get the tree (4.1.2), it is $_{(5)}T=T_1=T$.
\medskip
Clearly at any step of the present process \{the set of chosen nodes\} coincides with \{the set of final ramification nodes having one or two final branches of minimal length\} of the corresponding step of the second method 5.14.4.2. For example, these nodes of the first step for the tree (4.1.2) are the same nodes 8, 12, 19, 17 as the chosen nodes of the above last step.
\medskip
{\bf Example 5.14.5.2b.} Let us consider partitions $\varsigma_n: \{\binom{n+1}{2}=n+...+2+1\}$. The corresponding trees $T\in \g s^{-1}(\varsigma_n)_1$ for small $n$ are given in the below Table 5.14.5.3. For $n\le4$ all $T\in \g s^{-1}(\varsigma_n)_1$ are given, for $n=5$ only 3 of 15 elements of $\g s^{-1}(\varsigma_5)_1$ are given.
\medskip
\centerline{Table 5.14.5.3.}
\medskip
\settabs 5 \columns
\+$T\in \g s^{-1}(\varsigma_1)_1$&$T\in \g s^{-1}(\varsigma_2)_1$&$T\in \g s^{-1}(\varsigma_3)_1$&$T\in \g s^{-1}(\varsigma_4)_1$&$T\in \g s^{-1}(\varsigma_4)_1$\cr
\medskip
\+$^{\bullet}$&$^{\bullet<^\circ_\circ}$&$^{\bullet<^{\circ<^\circ_\circ}_{\circ-\circ}}$& $^{\bullet^{\nearrow^{\circ<^{\circ-\circ}_{\circ-\circ}}}_{\searrow_{\circ<^\circ_{\circ-\circ}}}}$&
$^{\bullet^{\nearrow^{\circ<^{\circ-\circ}_{\circ-\circ}}}_{\searrow_{\circ-\circ<^\circ_{\circ}}}}$& \cr
\medskip
\+1&1&3&5&4\cr
\medskip
\medskip
\medskip
\+$T\in \g s^{-1}(\varsigma_4)_1$&$T\in \g s^{-1}(\varsigma_5)_1$&$T\in \g s^{-1}(\varsigma_5)_1$&$T\in \g s^{-1}(\varsigma_5)_1$&\cr
\medskip
\+$^{\bullet^{\nearrow^{\circ<^{\circ<^\circ_{\circ}}_{\circ-\circ}}}_{\searrow_{\circ-\circ-\circ}}}$&   $^{\bullet^{\nearrow^{\circ^{\nearrow^{\circ<^\circ_{\circ-\circ}}}_{\searrow_{\circ-\circ-\circ}}}}_ {\searrow_{\circ<^{\circ-\circ}_{\circ-\circ-\circ}}}}$& $^{\bullet^{\nearrow^{\circ<^{\circ-\circ-\circ}_{\circ-\circ-\circ}}}_ {\searrow_{\circ^{\nearrow^\circ}_{\searrow_{\circ<^{\circ-\circ}_{\circ-\circ}}}}}}$& $^{\bullet^{\nearrow^{ \circ^{\nearrow^{\circ<^{\circ<^\circ_{\circ}}_{\circ-\circ}}}_{\searrow_{\circ-\circ-\circ}}}} _{\searrow_{\circ-\circ-\circ-\circ}}}$\cr
\medskip
\+6&11&6&10\cr
\medskip
The third line of the table indicates the quantity of the trees $_{(*+1)}T$ over a given $T$ considered as $_{(*)}T$. It is the quantity of non-ramification nodes of $T$ up to symmetries. Recall that the maximal value of this quantity over a $T\in \g s^{-1}(\varsigma_n)_1$ is $\binom{n}{2}+1$.
We see that for all trees of the table, except the first and the seventh trees (we count from the left to the right) having no symmetries, this quantity is less than the maximal value. We get that there are 15 (=5+4+6) trees $T\in \g s^{-1}(\varsigma_5)_1$. It is easy to calculate that there are 126 trees $T\in \g s^{-1}(\varsigma_6)_1$ and $\sim 2000$ trees $T\in \g s^{-1}(\varsigma_7)_1$. It is unlikely that there is a formula (even a recurrent formula like (4.1.4.1) ) for this quantity.
\medskip
{\bf 5.14.5.4.} Now we consider the general case. Let $\varsigma$ be a partition. We consider all partitions $\{\varsigma_\al\}$ of $\varsigma$ (see (5.14.3.2)), and we consider $F=\sqcup_\al \ T_\al$ (a disjoint union) where $T_\al\in \g s^{-1}(\varsigma_\al)_1$. This $F$ belongs to $\g s^{-1}(\varsigma)$, and all elements of $\g s^{-1}(\varsigma)$ are obtained by this manner.
\medskip
{\bf Remark.} The construction of the present Section 5.14.5 is inverse to the construction of Section 5.14.4.2. Really, the construction of 5.14.5 consists of consecutive gluing together of simple trees. The construction of 5.14.4.2 consists of consecutive cutting of the same simple trees.
\medskip
{\bf Example 5.14.5.5.} The only partition $\varsigma$ such that $\#(\g s^{-1}(\varsigma))=1$ are partitions $\varsigma_e(i,d):=\ \{i=d+d+...+d\}$ of equal parts (subscript $e$ means equal). Really, $\g s^{-1}(\varsigma)_1=\emptyset$ if $\phi_1$ of (5.14.5.1) is $\ne1$. This means that the only partition of $\varsigma_e(i,d)$ giving a forest is the lower trivial partition $\{d=d\}, \ \{d=d\}, \ \dots, \{d=d\}$. The above construction shows that $\g s^{-1}(\{d=d\})$ consists of one element --- the simple tree of length $d$. See Table A2.2, components $C_{i11}, \ C_{ii1}$ (any $i$), $C_{422}, \ C_{623}, \ C_{633}$.
\medskip
Obviously for other partitions $\varsigma$ we have $\#(\g s^{-1}(\varsigma))>1$. Really, if $\varsigma: \{i=d_1+...+d_l\}$ is a partition and $\exists \ \vk$ such that $d_\vk>d_{\vk+1}$ then $\varsigma$ has at least two partitions: the lower trivial and a partition such that $\Cal Z_1 = \{\vk,\vk+1\}$ and other $\Cal Z_\al$ consist of one element. For $\varsigma_1$ from (5.14.3.2) for this partition of $\varsigma$ we have $\g s^{-1}(\varsigma_1)\ne\emptyset$, hence $\g s^{-1}(\varsigma)$ contains at least two elements.

\medskip
{\bf 5.15. Irreducible components, field(s) of definition and Galois orbits.}
\medskip
Let $F$ be a fixed forest. We have $W(F) \subset \g F_{ij}$. We denote the set $\tilde\vf_{ij}(W(F))$ (it is a subset of $Irr(i,j)$; see (5.12) for a definition of $\tilde\vf_{ij}$) by $C(F)$. According (4.5), groups $\g G(F)\rtimes\Aut(F)$ and $\Ga_d$ act on $W(F)$. Conjecture 5.10 affirms that the set of orbits of the action of $\g G(F)\rtimes\Aut(F)$ on $W(F)$ is isomorphic to $C(F)$. We shall give here some propositions concerning the action of $\Ga_d$ on $C(F)$. Recall that an action of a group $G$ on a set $S$ is called a free action if $\forall \ g\in G, \ s\in S$ we have $\{g(s)=s\} \Rightarrow \{g=1\}$, i.e. every orbit has $\#G$ elements.
\medskip
Apparently, for a general case there is no simple description of Galois action on $C(F)$ (see Example 5.15.4), hence we give only some particular results.
\medskip
First, let us consider the case $F=T$ a tree such that $\Aut(T)=1$. We choose a longest simple subtree $B$ (a simple tree whose root coincide with the root of $T$, and whose final node is a final node of $T$) of $T$. We denote by $k$ the depth of the forest $T-B$.
\medskip
{\bf Proposition 5.15.1.} Conjecture 5.10 implies that for $T$ such that $\Aut(T)=1$, for all $w\in W(T)$ we have: $\im \vf(T,w)$ is defined over $\n Q(k+1)$ (see (4.5.4) for the notation). The action of $\Ga_{k+1}$ on $C(T)$ is free.
\medskip
{\bf Proof.} For any $w\in W(T)$ there exists the only element $\g g\in \g G(T)$ such that $\g g(w)_a$ is 0 on all nodes of $B$, i.e. the action of $\g G(T)$ on $W(T)$ is free. Hence, $\#C(T)=2^{i-d-l}$. These $w\in W(T)$ (such that $w_a(u)=0$ for all $u\in B$) form a set of representatives $R$ of orbits of the $\g G(T)$-action. For $w\in R$ we have:
\medskip
numbers $w_\mu(u)$, $u\in T$, belong to $\n Q(k+1)$;
\medskip
$\Ga_{k+1}$ acts on $R$, and this action is free --- this is obvious.
\medskip
Hence, $\Ga_{k+1}$ acts freely on $C(T)$. $\square$
\medskip
An example of such $T$ is given in (6.3), see the tree (6.3.1) and the Remark 6.3.1a.
\medskip
Let us consider now the case $F=T$, $\Aut(T)$ arbitrary. Let $T'$ be the minimal contraction of $T$ (see (5.11.5a)) such that the depth of $T'$ is $d':=d-\be$, and  $i'$, $l'$ the quantities of nodes, resp. final nodes of $T'$. Let $\vf: W(T)\to W(T')$ be a map of forgetting the values of $w$ on nodes of $T-T'$. We have also natural surjections $\g G(T) \to \g G(T')$, $\Aut(T) \to \Aut(T')$. They define a map of the set of orbits $\bar \vf: (\g G(T) \times\Aut(T))\backslash W(T) \to (\g G(T') \times\Aut(T'))\backslash W(T')$.
\medskip
{\bf Proposition 5.15.2.} $\bar \vf$ is an isomorphism; $\g G(T')\times \Aut(T')$ acts freely on $W(T')$. Hence, Conjecture 5.10 implies: the quantity of elements of the set $\tilde \vf(T,W(T))$ is $2^{i'-d'-l'}/\#\Aut(T')$.
\medskip
{\bf Proof.} First, we prove a lemma:
\medskip
{\bf Lemma 5.15.2.1.} Let $T=T'$ be a non-contractible tree, $g\in \Aut(T)$, $\g g\in \g G(T)$, $w\in W(T)$. If $g(w)=\g g(w)$ then $g=1$, $\g g=0$ (i.e. $\g G(T) \times\Aut(T)$ acts freely on $W(T)$).
\medskip
{\bf Proof.} Let us fix a final node $u$ of $T$ of the maximal depth $d$, and let $v=v(u,g)$ be the rightmost common predecessor of $u$, $g(u)$, i.e. a ramification node such that $u\in D_1(v)$, $g(u)\in D_2(v)$ (such $v$ exists and unique). Let $\ga=\ga(u,g)$ be the quantity of ramification nodes in the (only) way joining $v$ and $u$, including $v$ itself. We have $g^{2^\ga}(u)=u$. This is proved by induction by $\ga$. Really, $g$ interchanges $D_1(v)$ and $D_2(v)$, hence $g^2$ stabilizes both $D_1(v)$, $D_2(v)$. This means that $g^2(u)\in D_1(v)$, $\ga(u,g^2)<\ga(u,g)$, and the induction argument holds.
\medskip
Replacing $w$ by $\tilde \g g(w)$ for some $\tilde \g g\in \g G$ we can assume that $w_a(u)=0$, hence $w_a(v)=0$. Let the depth of $v$ is $d-k$, obviously $k\ge \ga$. We have $w_a(g(u))\equiv 2^{d-k}$ mod $2^{d-k+1}$ (because $w_a(v'_2)=2^{d-k}$ where $v'_1\in D_1(v)$, $v'_2\in D_2(v)$ are right neighbors of $v$).
\medskip
We consider $\g g$ as an element in $\n Z/2^d$. We have: $g(w)=\g g(w)$ implies $\g g=w_a(g(u))$. $g^{2^\ga}(u)=u$ implies $2^\ga\g g=0$. This equality and $k\ge \ga$, $\g g\equiv 2^{d-k}$ mod $2^{d-k+1}$ imply $k=\ga$, all nodes $g^\de(u)$ for $0\le\de<2^\ga$ are different and form a complete subtree of $T$ of node $v$ and depth $k$.
\medskip
Let $\bar u$ be another final node of $T$ of the maximal depth $d$. Numbers $\bar k=\bar\ga$ for $\bar u$ are the same, because $\g g$ is the same. If $k\ne 0$ we get a contradiction to the condition that $T$ is a non-contractible. $\square$
\medskip
So, if $T$ is non-contractible, then the group $\Aut(T)\times \g G(T)$ acts freely on $W(T)$, hence for this case (5.15.2) is proved. Let us consider the general case. We have a trivial
\medskip
{\bf Lemma 5.15.2.2.} Let $T$ be a complete tree. In this case, the group $\Aut(T)$ acts simply transitively on $W(T)$.
\medskip
{\bf Proof.} Induction by $k$ --- the depth of $T$. Let $r$ be the root of $T$ and $r_1$, $r_2$ its right neighbors. They are roots of $D_1(r)$, $D_2(r)$ (see (4.6)). Let $w_1$, $w_2$ be two weights of $T$. Let us consider the case $w_1(r_1)=w_2(r_1)$. We denote by $w_{\al\be}$ ($\al, \ \be=1, \ 2$) the restrictions of $w_\al$ to $D_\be(r)$. There exist (by induction supposition) automorphisms $\g g_\be$ of $D_\be(r)$ such that $\g g_\be(w_{1\be})=w_{2\be}$ for $\be=1, \ 2$. The pair $(\g g_1, \g g_2)\rtimes 1 \in (G_{k-1}\times G_{k-1}) \rtimes \n Z/2=G_k$ belongs to $\Aut(T)$, satisfies $(\g g_1, \g g_2)(w_1)=w_2$, and it is the only element of $\Aut(T)$ having this property.
\medskip
If $w_1(r_1)\ne w_2(r_1)$ we use the same arguments, but now we consider an element $(\g g_1, \g g_2) \rtimes \s\in (G_{k-1}\times G_{k-1}) \rtimes \n Z/2=G_k$ where $\s$ is the non-zero element of $\n Z/2$.  $\square$
\medskip
We denote by $\Aut(T/T')$ the kernel of the natural map $\Aut(T)\to\Aut(T')$. Lemma 5.15.2.2 implies that if $\vf(w_1)=\vf(w_2)$ then $\exists \ g\in \Aut (T/T')$ such that $w_2=g(w_1)$. Using Lemma 5.15.2.1, we get the proposition. $\square$
\medskip
To consider the case of a forest, we need one trivial lemma more. Namely, let $G$ be any group acting on any set $\g S$, and $\eta$ a number. The group $G^\eta\rtimes S_\eta$ acts naturally on $\g S^\eta$ ($G^\eta$ acts on $\g S^\eta$ coordinatewise, and $S_\eta$ interchange coordinates in $\g S^\eta$).
\medskip
{\bf Lemma 5.15.3.} In the above notations, let $\g O$ be the set of orbits of the action of $G$ on $\g S$. Then the set of orbits of the  action of $G^\eta\rtimes S_\eta$ on $\g S^\eta$ is $S^\eta(\g O)$ (the symmetric product).
\medskip
{\bf Proof.} Let $\tilde \g S\subset \g S$ be a set of representatives of orbits of the action of $G$ on $\g S$. We have: $\tilde \g S^\eta\subset \g S^\eta$ is a set of representatives of orbits of the action of $G^\eta$ on $\g S^\eta$. Let $\tilde \g S^\eta \to S^\eta(\g S)$ be the canonical projection (here $S^\eta$ means the $\eta$-symmetric product of a set). Let $\bar \g S^\eta\subset \tilde \g S^\eta$ be a set of representatives in fibers of this projection. It is a set of representatives of orbits of the action of $G^\eta\rtimes S_\eta$ on $\g S^\eta$, because any element of $\g S^\eta$ is $G^\eta\rtimes S_\eta$-equivalent to an element of $\bar \g S^\eta$, and no two elements of $\bar \g S^\eta$ are $G^\eta\rtimes S_\eta$-equivalent. Finally, $\bar \g S^\eta$ is isomorphic to $S^\eta(\g O)$ as a set. $\square$
\medskip
First, we apply this lemma to the case $F=T\sqcup T\sqcup...\sqcup T$ --- a disjoint union of $\eta$ copies of $T$. (5.15.2), (5.15.3) give us for this case a description of the set of orbits of the action of $\g G(F)\rtimes \Aut(F)$ on $W(F)$ and hence, by Conjecture 5.10, of the set $C(F)$.
\medskip
For a general case $F=\sqcup_{\be=1}^\de \eta_\be T_\be$ where $T_\be$ are different (see (4.2.2)), we have that the set of orbits of the action of $\g G(F)\rtimes \Aut(F)$ on $W(F)$ is the product of the corresponding sets for $F_\be:= \eta_\be T_\be$. This gives a complete description of $C(F)$.
\medskip
There is no simple description of the Galois orbits, i.e. the action of $\Ga_{k+1}$ on $C(F)$.
We give here an example showing that Galois orbits of elements of $C(T)$ can be different for different $w\in W(T)$.
\medskip
{\bf Example 5.15.4.} Let $T$ be the below tree (it is the tree $T_1$ of $F$ of Example 4.5.3). Its final nodes 1,...,4 are marked.
\settabs 20 \columns
\medskip
\+&&&$^{\bullet^{\nearrow^{\circ<^{\circ-\circ1}_{\circ-\circ-\circ2}}}_{\searrow_{\circ<^{\circ-\circ-\circ3}_{\circ-\circ4}  }}}$& &
\cr
\medskip
We denote an additive weight $w$ such that $w(i$-th final node)$=a_i$ (clearly a weight is defined uniquely by its values on final nodes) by $(a_1,...,a_4)$, where $a_1, a_4\in \n Z/8$, $a_2, a_3\in \n Z/16$. We choose the longest simple subtree $B$ from the root to the third node. There are 32 orbits of $\g G$ on $W(T)$,\footnotemark \footnotetext {Because for $T$ we have $i=13$, $l=4$, hence according (4.4.1) $\#(W(T))=512$. Since $d=4$, we have $\#\g G=16$. The action of $\g G$ on $W(T)$ is free, hence we have 32 orbits.} each of them has exactly one representative having $a_3=0$. We have $k=3$. The group $Aut(T)$ is $\n Z/2$ (we denote by $\bh$ its non-trivial element; it interchanges nodes 1 and 4, 2 and 3). We have: $\#C(T)=16$. Let $\s$ be the complex conjugation.
\medskip
Let us consider $w_1=(3,1,0,6)$. We have $\bh(w_1)=(6,0,1,3)$ (it interchanges values of $w_1$ on nodes 1 and 4, 2 and 3), $\s(3,1,0,6)=(5,15,0,2)$ (because $(\s(w_a))(x)=-w_a(x)$ ) and $(6,0,1,3) \sim (5,15,0,2)$ where $\sim$ means the same $\g G$-orbit (because $(6,0,1,3)=\g g(5,15,0,2)$ for $\g g=1 \ \in \n Z/16$ identified with $0,1,\dots,15$). This means that $\s$ stabilizes $\im \vf(T,w_1)$. Hence, it is defined over $\n Q(4)^+$, where $\n Q(k)^+:=\n Q(k)\cap \n R$, and the Galois orbit of $\im \vf(T,w_1)$ consists of 4 elements.
\medskip
Let us now consider $w_2=(3,1,0,2)$. We have $\bh(w_2)=(2,0,1,3)$, $\s(w_2)= (5,15,0,6)$ and $(2,0,1,3)\not\in \g G(5,15,0,6)$ (because differences $2-5, \ 0-15, \ 1-0, \ 3-6$ are not images of one element of $\n Z/16$). This means that the complex conjugation acts non-trivially of the irreducible component of corresponding to $w_2$. Considering the action of other elements of Gal$(\n Q(4)/\n Q)$ on this irreducible component we get that they all act non-trivially on it, i.e. it is defined over $\n Q(4)$, and its Galois orbit consists of 8 elements.
\medskip
{\bf 6. Examples.}
\medskip
{\bf 6.1. $T$ is the simple tree of length $i$}. We have $j=1$, $m=i+1$, $d=m-2$. There is one $\g G$-orbit on $W(T)$, hence we can assume that $w_a(u)=0$ for all $u\in T$. We choose the order $O$ of nodes from the right to the left. $A$ of (5.4.3) is the following $m\times (m+1)$-matrix (here $c=c_1$):
$$A=\left(\matrix 0&1&0&c&0&c^2&0&c^3&\dots\\1&0&c&0&c^2&0&c^3&0&\dots \\1&-c&c^2&-c^3&c^4&&\dots&&(-1)^mc^m \\1&-c^2&c^4&-c^6&c^8&&\dots&&(-1)^mc^{2m}\\1&-c^4&c^8&-c^{12}&c^{16}&&\dots&&(-1)^mc^{4m} \\1&-c^8&c^{16}&-c^{24}&c^{32}&&\dots&&(-1)^mc^{8m} \\ \dots & \dots & \dots & \dots & \dots & \dots & \dots & \dots & \dots &\\1&-c^{2^{m-3}}&c^{2\cdot2^{m-3}}&-c^{3\cdot2^{m-3}}&c^{4\cdot2^{m-3}}&&\dots&&(-1)^mc^{m\cdot2^{m-3}}\endmatrix \right)\eqno{(6.1.1)}$$
$\g c_*$ of (5.4.4) is $(\sqrt{c},-\sqrt{c},-c,-c^2,-c^4,-c^8,\dots, -c^{2^{m-3}})$
and $V_{m,m+1}(\g c_*)$ of (5.4.6) is $V_{m,m+1}(\sqrt{c},-\sqrt{c},-c,-c^2,-c^4,-c^8,\dots, -c^{2^{m-3}})$, hence

$$\la_\vk=(-1)^{m-\vk}\sigma_{m-\vk}(\sqrt{c},-\sqrt{c},-c,-c^2,-c^4,-c^8,\dots, -c^{2^{m-3}})$$ $$=(-1)^m[\ \sigma_{m-\vk}(c,c^2,c^4,c^8,\dots, c^{2^{m-3}})-c\ \sigma_{m-\vk-2}(c,c^2,c^4,c^8,\dots, c^{2^{m-3}})\ ]\eqno{(6.1.2)}$$
In homogeneous coordinates $(c:c')$ we have
\medskip
$$\la_\vk= (-1)^m[\ (c')^{2^{m-2}}(\sigma_{m-\vk}(c/c',(c/c')^2,(c/c')^4,(c/c')^8,\dots, (c/c')^{2^{m-3}})$$ $$-c/c'\ \sigma_{m-\vk-2}(c/c',(c/c')^2,(c/c')^4,(c/c')^8,\dots, (c/c')^{2^{m-3}}))\ ]\eqno{(6.1.3)}$$

{\bf 6.1.4.} Let us calculate its $\g d$, $\mu$ and $\g s(T)$ using (5.11.11), (5.13.2), (5.14.4.2). We have $\ga_1=0$ (notation of (5.11.6)), $d_1=d=i-1$, $f_*=1$, hence according (5.11.11) we have $\g d=2^{i-1}$. Since $T$ has no ramification points, the product (5.13.3) is empty, hence $\mu=1$.
\medskip
We have $\g y_1=i$, $\phi_1=1$ (notations of (5.14.4.2)), hence $\g s(T)$ is the trivial partition $i=i$. According 5.14.5.5, a simple tree $T$ is the only forest whose $\g s(T)$ is the trivial partition $i=i$.
\medskip
The only irreducible component in Table A2.2 having the above characteristics ($j$, $\g d$, $\mu$ and partition) is $\bar C_{i11}$. This confirms Conjecture 5.13.1.
\medskip
{\bf 6.2. Construction of components with maximal $j$.}
\medskip
Let $\varsigma=\{i=d_1+...+d_l\}$ be a fixed partition. Here we describe\footnotemark \footnotetext{We describe the subset of such $C\in Irr(i,j)$ in combinatorial terms, and we give parametric equations of the corresponding irreducible components.} $C\in Irr(i,j)$ such that $\pi(i,j)(C)=\varsigma$ and $j=j(C)$ is the maximal possible (for this $\varsigma$). This means that a partition of $\varsigma$ corresponding to $C$ is the lower trivial partition, $j=l$ and all trees $T_\al$ are simple trees of length $d_\al$. We have $F=\sqcup_{\al=1}^l T_\al$. As above there is one $\g G$-orbit on $W(F)$, hence we can assume that $w_a(u)=0$ for all $u\in F$. We get that this $F$ is the only forest consisting of $l$ trees such that $\g s(F)=\s$, and if another forest $\bar F$ satisfies $\g s(\bar F)=\s$ then $j(\bar F)<l$.

We choose the order $O$ of nodes: First, from the right to the left of the tree $T_1$, second --- of the tree $T_2$ etc. until the tree $T_l$. The analog of (6.1.2) for the present case is (here $m=i+l$):
$$\la_\vk=(-1)^{m-\vk}\sigma_{m-\vk}(\sqrt{c_1},-\sqrt{c_1},-c_1,-c_1^2,\dots, -{c_1}^{2^{d_1-2}}, $$
$$\sqrt{c_2},-\sqrt{c_2},-c_2,-c_2^2,\dots, -{c_2}^{2^{d_2-2}}, \dots,
\sqrt{c_l},-\sqrt{c_l},-c_l,-c_l^2,\dots, -{c_l}^{2^{d_l-2}}).\eqno{(6.2.1)}$$

{\bf 6.2.2.} Let $\{i=f_1+...+f_{d_1}\}$ be a partition conjugate to $\varsigma$. We have formulas: $$\g d(C)=\frac{2^{i-l}\cdot l!}{(f_1-f_2)!\cdot(f_2-f_3)!\cdot...\cdot(f_{d_1-1}-f_{d_1})!\cdot f_{d_1}!}; \ \ \ \mu(C)=\frac{i!}{d_1!\cdot d_2! \cdot ... \cdot d_l!}$$
This follows immediately from the above (5.11.11), (5.13.2). Table A2.2 confirms these results.
\medskip
{\bf 6.2.3. Example.} Let $\varsigma=\{4=2+2\}$, $m=6$ (this is $C_{422}$ of Appendix, A2.2). (6.2.1) becomes
$$\la_\vk=(-1)^\vk\sigma_{6-\vk}(\sqrt{c_1},-\sqrt{c_1},-c_1,\sqrt{c_2},-\sqrt{c_2},-c_2), \hbox{ i.e.}$$
$$\vf(F,w)((c_1:1);(c_2:1))=(\la_0:...:\la_6)$$ $$=(\tilde c_2^2:\tilde c_1\tilde c_2:\tilde c_2-\tilde c_1\tilde c_2:-\tilde c_1^2:-\tilde c_1+\tilde c_2:\tilde c_1:1)$$ where $\tilde c_1:=c_1+c_2$, $\tilde c_2:=c_1c_2$. It is obvious that $\dim \im \vf(F,w)=2$ (the projection to the plane $(\la_4:\la_5:\la_6)$ is surjective) but not 0 or 1, i.e. Conjecture 5.8 trivially holds for this case.
\medskip
{\bf 6.3. Construction of components $C$ having the partition $i=d_1+d_2$.} Number $j=j(C)$ can be 1 or 2. The case $j=2$ is treated in (6.2), so we consider the case $j=1$. We have $m=i+1$. If $d_1=d_2$ then $j$ cannot be 1 (see 5.14.4.2.0, 5.14.5.1a), so we let $d_1>d_2$. The set of trees $T$ is parameterized by a number $\g n$ satisfying $d_1-d_2-1\ge \g n \ge 0$. The corresponding tree $T(\g n)$ is the following (numbers at nodes are their depths):
\medskip
$$\circ\to\circ\to ... \to\circ\to\underset{\g n}\to{\circ}<^{\circ\to\circ\to ... \to\circ\to\overset{\g n+d_2}\to{\circ}}_{\circ\to\circ\to ...\to \circ\to\circ \to \circ\to\underset{d_1-1}\to{\circ}}\eqno{(6.3.1)}$$
\medskip
Any $\g G$-orbit of a weight contains exactly one weight $w$ taking additive values 0 at all nodes of the long branch. This $w$ is uniquely defined by the value of $w(u)$ where $u$ is the final node of the short branch (of depth $\g n+d_2$). We have $w_\mu(u)=\g z_{d_2}^\vk$ where $\vk\in \n Z/2^{d_2}$ is odd. This means that for a fixed $\g n$ the quotient set of $W(T(\g n))$ by the action of $\g G(T(\g n))$ contains $2^{{d_2}-1}$ elements. They form one Galois orbit. If $\g n<d_1-{d_2}-1$ then $\Aut(T(\g n))=1$ and we have $2^{{d_2}-1}$ components of $X(m,i)$ defined over $\n Q({d_2})$ forming a $\Ga_{d_2}$-orbit. If $\g n=d_1-{d_2}-1$ (two final branches are of equal length) then $\Aut(T(\g n))=\n Z/2$, we have $2^{{d_2}-2}$ components of $X(m,i)$ defined over $\n Q({d_2})^+$ forming a Gal$(\n Q({d_2})^+/\n Q)$-orbit.
\medskip
This supports Conjecture 2.2.
\medskip
{\bf Remark 6.3.1a.} The above tree (case $\g n<d_1-{d_2}-1$) is a particular case of a tree of Proposition 5.15.1. We have: $B$ of (5.15.1) is a subtree of (6.3.1) from the root to the rightmost low branch, $T-B$ of (5.15.1) is the upper branch (from the upper right neighbor of the ramification node, until the upper final node), $k$ of (5.15.1) is $d_2-1$.
\medskip
{\bf Example 6.3.2.} Let the partition be $6=4+2$. If $\g n=1$ then the corresponding irreducible component is $C_{616}$ of Table A2.2, it is defined over $\n Q$. If $\g n=0$ then the two corresponding irreducible components are $C_{618a}, \ C_{618b}$ of Table A2.2. They are defined over $\n Q(\g z_2)$ and they are Gal$(\n Q(\g z_2)/\n Q)$-conjugate. This is the simplest example of components defined over a field larger than $\n Q$ (recall that the notion of a field of definition is from (1.7). It can happen that $C_{618a}, \ C_{618b}$ have models over $\n Q$, compare with $\n P^1$ of Example 1.7.1. We do not consider this problem. Since the complex conjugation interchanges $C_{618a}$ and $C_{618b}$ we get that $\n Q(\g z_2)$ is the minimal field of definition of them).
\medskip
{\bf 6.3.3.}
Let us calculate $\g d$ and $\mu$ of the above components using formulas (5.11.11), (5.13.2). We have $\ga_1=1$ if $d_2=1$, $\g n=d_1-2$ and $\ga_1=0$ otherwise. We have $d=d_1-1$, hence $\g d=2^{d_1-2}=2^\g n$ if $d_2=1$, $\g n=d_1-2$ and $\g d=2^{d_1-1}$ otherwise. Since there is only one ramification point, the product of (5.13.3) gives us $\mu=\binom{d_2+d_1-1-\g n}{d_2}$. Table A2.2 confirms these results (see Table A2.3 for a list of these components and their $\g n$, $d_1$, $d_2$).
\medskip
{\bf 6.4. Case of a forest consisting of equal complete trees.} Let $F=\eta \ T$ where $T$ is the complete tree of depth $k$. $\g G(T)$ acts transitively on $W(T)$ (see Lemma 5.15.2.2), and $\g G(F)$ acts transitively on $W(F)$. We have $\ga=k$, where $\ga$ is from (5.11.6), hence $$\deg \im \vf(F,w)=1\eqno{(6.4.1)}$$ It is easy to see that these are the only $(F,w)$ such that $\deg \im \vf(F,w)=1$. For $\eta=1$ (5.13.3) gives
$$\mu(\im \vf(T,w))=\binom{2}{1}^{2^{k-1}}\binom{6}{3}^{2^{k-2}}\binom{14}{7}^{2^{k-3}}\dots \binom{2^{k+1}-2}{2^{k}-1}\eqno{(6.4.2)}$$
and for any $\eta$ (5.13.2) gives
$$\mu(\im \vf(F,w))=\frac{[(2^{k+1}-1)\eta]!}{[(2^{k+1}-1)!]^\eta}\ \ \mu(\im \vf(T,w))^\eta\eqno{(6.4.3)}$$
The partition $\g s(F)$ is the conjugate to $i=(2^{k+1}-1)\eta=\{2^k\eta+2^{k-1}\eta+...+2\eta+\eta\}$. See 11.3 and Table A7 for more details on these $\im \vf(F,w)$.
\medskip
{\bf 6.5.} Variety $X(5,4)$. It is a curve in $\n P^5$. All its components are the minimal ones (see Remark 2.5.1). According Table A2.2, is consists of 3 irreducible components. We give formulas for $\vf(c)$, $c\in \n A^1$, for all these components. The order $O$ is from the right to the left (the final results do not depend on the order, see (5.5.1.1b)).
\medskip
{\bf 6.5.1.} Component $\bar C_{411}$ (see (6.1)). We have $$\g c_*=(\sqrt{c},-\sqrt{c},-c,-c^2,-c^4)$$
$$\vf(c)=(-c^{8}:-c^7-c^6-c^4:c^7-c^5-c^3-c^2:c^6+c^5+c^3-c:c^4+c^2+c:1)$$

{\bf 6.5.2.} Component $\bar C_{412}$. We have $$\g c_*=(\sqrt{c},-\sqrt{c},\sqrt{-c},-\sqrt{-c},-c^2)$$
$$\vf(c)=(-c^{4}:-c^2:0:0:c^2:1)$$

{\bf 6.5.3.} Component $\bar C_{413}$. Let the weight of the nodes of the long branch be 1. We have $$\g c_*=(\sqrt{c},-\sqrt{c},-c,\sqrt{-c^2},-\sqrt{-c^2})$$
$$\vf(c)=(-c^{4}:-c^3:c^3-c^2:c^2-c:c:1)$$
It is possible to check by a direct calculation that the above equations really give us points of  $X(5,4)$. According the definitions 1.2, 1.3, $(a_0:...:a_5)\in X(5,4) \iff $ the characteristic polynomial of $\g M(5)(a_0,...,a_5)$ is equal to $U^4$. This is really so:
$$\left|\matrix -c^7-c^6-c^4-U& c^6+c^5+c^3-c&1&0\\-c^{8}&c^7-c^5-c^3-c^2-U&c^4+c^2+c&0\\0&-c^7-c^6-c^4& c^6+c^5+c^3-c-U&1 \\ 0&-c^{8}&c^7-c^5-c^3-c^2&c^4+c^2+c-U\endmatrix \right|=$$ $$=\left|\matrix -c^2-U& 0&1&0\\-c^{4}&-U&c^2&0\\0&-c^2&-U&1 \\ 0&-c^{4}&0&c^2-U\endmatrix \right|=\left|\matrix -c^3-U& c^2-c&1&0\\-c^{4}&c^3-c^2-U&c&0\\0&-c^3&c^2-c-U&1 \\ 0&-c^{4}&c^3-c^2&c-U\endmatrix \right|=U^4$$
(equalities in $\n C[U]$). See (7.5.2), (7.5.3) for other examples.
\medskip
{\bf 7. The second construction: from $\bar C_{ij*}$ to $C_{ij*}(m)$.}\footnotemark \footnotetext{The first construction is in Section 5.}
\medskip
We use notations of (2.5). The second construction shows how --- starting from the minimal irreducible component $\bar C_{ij*}$ --- we can get all components of its series.
\medskip
For given $m>0$, $\vartheta>0$ we define a map $\nu=\nu_{m,\vartheta}: \n P^{m}\times \n P^{\vartheta}\to \n P^{m+\vartheta}$ like a product of polynomials. Namely, let
\medskip
$((\la_0:\la_1:...:\la_{m}); (b_0:...:b_{\vartheta}))\in \n P^{m}\times \n P^{\vartheta}$.
\medskip
We associate $(\la_0:\la_1:...:\la_{m})$ with $P_1:=\sum_{i=0}^m \la_ix^i$ and $(b_0:...:b_{\vartheta})$ with $P_2:=\sum_i b_ix^i$, then $\nu(\la_*; b_*)$ is associated with $P_1P_2$. Explicitly, the coordinates $(a_0:...:a_{m+\vartheta})\in \n P^{m+\vartheta}$ of $\nu(\la_*; b_*)$ are defined as follows: for $s\in[0,...,m+\vartheta]$ we let
$$a_s=\sum_{\ga\in \n Z} \la_\ga b_{s-\ga}$$
where $\la_*=0$, resp. $b_*=0$, if $*\not\in [0,...,m]$, resp. $[0,...,\vartheta]$.
\medskip
{\bf Proposition 7.1.} If Conjecture 0.2.4 is true\footnotemark \footnotetext{The reader can think that Conjecture 0.2.4 is in characteristic 2. Really, it is over $\n C$, see Section 0.2.} for $n=1$ then the following holds: $\forall \ m,i,\vartheta$
$$\nu(X(m,i)\times \n P^{\vartheta})\subset X(m+\vartheta,i)\eqno{(7.1.1)}$$

{\bf Proof.} We can represent the above $P_2$ as a product of linear polynomials: $P_2=\prod_{\vk=1}^\vartheta P_{2\vk}$ where $P_{2\vk}=b_{0\vk}+b_{1\vk}x$. Now, we define inductively numbers $\la_{\vk \ga}$ as follows:
$$\la_{0\ga}=\la_\ga,\ \ \ga=0,\dots, m$$
$$\nu_{m+\vk,1}((\la_{\vk 0}:...:\la_{\vk ,m+\vk});(b_{0,\vk+1}:b_{1,\vk+1}))=(\la_{\vk+1, 0}:...:\la_{\vk+1 ,m+\vk+1})$$
We have $\nu_{m,\th}((\la_0:\la_1:...:\la_{m}); (b_0:...:b_{\vartheta}))=(\la_{\vt 0}:...:\la_{\vt ,m+\vt})$. This implies that if (7.1.1) holds for $\vartheta=1$ then it holds for any $\vartheta$. Hence, we shall prove the proposition for $\vt=1$.
\medskip
Let us consider a point $(\la_0:...:\la_m)\in X(m,i)$,
a point $(b_0:b_1)\in \n P^1$, and let $\nu((\la_0:...:\la_m),(b_0:b_1))=(a_0:...:a_{m+1})\in \n P^{m+1}$.
\medskip
We have: the matrix $\g M(m+1)(a_0,...,a_{m+1})$ becomes $\Cal M_{nt}(\la_0,...,\la_{m};\ 1,m)$ (see (0.1.8) for the notation) after the substitution $b_0 \to t$, $b_1 \to -1$, hence
$$D(m+1,i)(a_0:...:a_{m+1})=\sum_{\g j=0}^{m-i}(-1)^{m-i-\g j}H_{i\g j,21}(\la_0,...,\la_m)\ b_0^\g j b_1^{m-i-\g j}\eqno{(7.1.2)}$$
This implies the proposition. $\square$
\medskip
{\bf Remark 7.2.} The idea used in the above proof of (7.1) also can be used for a proof of the following (7.2.1):
\medskip
{\bf Proposition 7.2.1.} If Conjecture 0.2.4 holds for $n=1$ for all $m$, then it holds for all $n$.
\medskip
{\bf Proof.} The matrix $\g M(m+\vartheta)(a_0,...,a_{m+\vartheta})$ becomes $\Cal M_{nt}(\la_0,...,\la_{m};\vartheta;m+\vartheta-1)$ after the substitution
$b_\g j \to (-1)^\g j \binom{\vartheta}{\g j} t^{\vartheta-\g j}$, hence $$D((m+\vartheta),i)(a_0,...,a_{m+\vartheta})=\sum_{\g J}\tilde H_{\g J,\vartheta}(\la_0,...,\la_m)b^\g J$$ where
$\g J=(\g j_0,...,\g j_{\vartheta})$ is a multiindex satisfying $\sum_{\varkappa=0}^{\vartheta} \g j_\varkappa=m+\vartheta-1-i$, $b^\g J:=\prod_{\varkappa=0}^{\vartheta} b_\varkappa^{\g j_\varkappa}$ and $\tilde H_{\g J,\vartheta}$ are some polynomials. This implies that
$H_{i\g j,2\vartheta}$ are linear combinations of the corresponding $\tilde H_{\g J,\vartheta}$.
\medskip
If Conjecture 0.2.4 holds for $n=1$ for all $m$, then (7.1.1) holds, and hence some powers of $\tilde H_{\g J,\vartheta}$ belong to $<D(m,0),\dots,D(m,i-1)>$.
Therefore, the same holds for $H_{i\g j,2\vartheta}$.  $\square$
\medskip
Let us fix $i$, $j$ and an element $C_{ij\g k}\in Irr(i,j)$. We consider $\bar C_{ij\g k}\subset X(i+j,i)\subset \n P^{i+j}$ from (2.5), and let $m\ge i+j$ be arbitrary. We let $\vartheta:=m-(i+j)$.
\medskip
{\bf Conjecture 7.3.} $C_{ij\g k}(m)=\nu_{i+j,\vartheta}(\bar C_{ij\g k}\times \n P^{\vartheta})$.
\medskip
This conjecture should be understood as follows. According Proposition 7.1, we get that $\nu_{i+j,\vartheta}(\bar C_{ij\g k}\times \n P^{\vartheta})$ is (conjecturally) an irreducible component of $X(m,i)$. We conjecture that if we denote this irreducible component of $X(m,i)$ by $C_{ij\g k}(m)$ then all affirmations of (2.4) hold.
\medskip
For example, we have the following
\medskip
{\bf 7.4. Justification of (2.4.2).} We have dim $\bar C_{ij\g k}=j$, dim $\nu_{i+j,\vartheta}(\bar C_{ij\g k}\times \n P^{\vartheta})=j+\vartheta=m-i$, hence
$$\deg \nu_{i+j,\vartheta}(\bar C_{ij\g k}\times \n P^{\vartheta})=$$ $$=\deg \bar C_{ij\g k} \cdot \hbox{ (the degree of the Segre embedding }\n P^j \times \n P^{\vartheta} \to \n P^{j\cdot \vartheta+j+\vartheta}\ )$$
Since the degree of the Segre embedding is $\binom{m-i}{j}$ we get (2.4.2).
 \medskip
Proof (or, at least, a justification) of the fact that the multiplicity and
the Jordan form of $\nu(\bar C_{ij\g k}\times \n P^{\vartheta})$ coincide with the ones of $\bar C_{ij\g k}$, as well as other affirmations of Conjectures 2.4, is a subject of future research. A simple particular case is treated in (11.3.4).
\medskip
{\bf 7.5.} Example: Variety $X(6,4)$. It is a surface in $\n P^6$. According (2.4.1), $Irr(X(6,4))$ consists of 6 elements (lines $C_{411}-C_{423}$ of Table A2.2), components corresponding to the lines $C_{421}-C_{423}$ are the minimal ones. Let us give parametric equations of some of them. We use parameters $(c:c')\in \n P^1$, $(b_0:b_1)\in \n P^1$ for $C_{412}(6)$, $(c_1:1)\in \n P^1$, $(c_2:1)\in \n P^1$ for $\bar C_{421}$,  $\bar C_{423}$. According Proposition 7.1, parametric equations of $C_{412}(6)$ are
$$(-c^{4}b_0: -c^{4}b_1   -c^3c'b_0: -c^3c'b_1  + (c^3c'-c^2{c'}^2)b_0:  (c^3c'-c^2{c'}^2)b_1+   (c^2{c'}^2-c{c'}^3)b_0: $$ $$ :(c^2{c'}^2-c{c'}^3 )b_1+  c{c'}^3b_0:  c{c'}^3b_1+   {c'}^4b_0 : {c'}^4b_1)\eqno{(7.5.1)}$$
Parametric equations for $C_{411}(6)$, $C_{413}(6)$ can be obtained from (6.5.1), (6.5.3) using the same method.

{\bf 7.5.2.} Component $\bar C_{421}$ (see (6.2)). We have $$\g c=(\sqrt{c_1},-\sqrt{c_1},-c_1,c_1^2,\sqrt{c_2},-\sqrt{c_2})$$
$$\vf((c_1:1), (c_2:1))=(-c_1^4c_2 : c_1^3c_2+ c_1^2c_2 :-c_1^4 -c_1^3c_2 +c_1c_2 : $$
$$ :-  (c_1+c_2 )(  c_1^2  + c_1):   c_1^3-   c_1-c_2 :  c_1^2+c_1:1)$$

{\bf 7.5.3.} Component $\bar C_{422}$ is described in (6.2.3). For the component $\bar C_{423}$ we have
$$\g c=(\sqrt{c_1},-\sqrt{c_1},\sqrt{-c_1},-\sqrt{-c_1},\sqrt{c_2},-\sqrt{c_2})$$
$$\vf((c_1:1), (c_2:1))=(-c_1^4c_2:0:-c_1^4:0:-c_2:0:1)$$
\medskip
\medskip
{\bf Lifts.} There exist elementary constructions: irreducible components of $X(m,i)$ can give rise another irreducible components of $X(\bar m,\bar i)$ for some $\bar m, \bar i$. We consider two such constructions which are called the even lift and the odd lift (according the parity of $m$), respectively. They can be useful for generalizations of the present theory to other cases. Namely, $X(m,i)$ of the present paper are a special case of the varieties $X_\g c(q,n,m,i)$, see the Introduction (also the base field can be either of finite characteristic or of characteristic 0). Analogs of constructions of lifts for all these cases can help us to get an analog of the present theory for other $q$ and $n$. \medskip
Interpretation of the lifts in terms of weighted forests is an exercise for a student. Some constructions of the present part depend on Conjectures 2.1, 2.4.
\medskip
\newpage
{\bf 8. Even lift.}
\medskip
{\bf Proposition 8.1.} Let $m$ be even. The set of points $(a_0:...:a_m)\in \n P^m$ satisfying the condition:\footnotemark \footnotetext{Here and below we consider a representative $(a_0,...,a_m)\in \n A^{m+1}$ of $(a_0:...:a_m)\in \n P^m$. All statements do not depend on its choice.}
$$\g M(m)(a_0,...,a_m)\hbox{ has rank }\le m/2-1$$
is the linear subspace of $\n P^m$ defined by the equations $a_1=a_3=...=a_{m-1}=0$.
\medskip
{\bf Proof.}\footnotemark \footnotetext{The authors are grateful to an anonymous reviewer who indicated them the present simplification of the original version of the proof.} If $a_1=a_3=...=a_{m-1}=0$ then clearly the rank of $\g M(m)(a_0,...,a_m)$ is $\le m/2-1$. If the rank of $\g M(m)(a_0,...,a_m)$ is $\le m/2-1$ then all the odd rows of $\g M(m)(a_0,...,a_m)$ are linearly dependent. Let us assume that not all $a_1, \ a_3, \ ... \ a_{m-1}$ are 0, and let $\vk$ be the minimal odd index such that $a_\vk\ne0$. The $\frac m2 \times \frac m2$-submatrix of $\g M(m)(a_0,...,a_m)$ formed by all odd rows and $\frac{\vk+1}2$-th and subsequent columns is a triangular matrix with non-zeros on diagonal --- a contradiction. $\square$
\medskip
{\bf 8.2.} We denote the linear subspace $a_1=a_3=...=a_{m-1}=0$ by $\La_1(m)$. We have an obvious
\medskip
{\bf Lemma 8.3.} For an element $(a_0:0:a_2:0:...:0:a_m)\in \La_1(m)$ we have: $Ch(\g M(m)(a_0:0:a_2:0:...:0:a_m))=\pm U^{m/2}Ch(\g M(m/2)(a_0:a_2:...:a_m))$ (see (1.2) for the definition of $Ch$). $\square$
\medskip
Let us fix $i,j,\g k,m$ and the corresponding elements $C_{ij\g k}\in Irr(i,j)$, $C_{ij\g k}(m)\subset X(m,i)$.
\medskip
{\bf Definition 8.4.} Let us assume the truth of Conjecture 2.1. The even lift of $C_{ij\g k}(m)$ is an irreducible component $L_{e;m}(C_{ij\g k})$ (subscript $e$ means even) of $X(2m,m+i)$ defined as follows:
$(a_0:...:a_{2m})\in L_{e;m}(C_{ij\g k}) \iff a_1=a_3=...=a_{2m-1}=0, \ \  (a_0:a_2:...:a_{2m})\in C_{ij\g k}(m)$.
\medskip
The fact that it is really an irreducible component of $X(2m,m+i)$ is obvious, because $\dim C_{ij\g k}(m)=m-i=\dim X(2m,m+i)$.
\medskip
We need a version of definitions of Section 2. Namely, let $Irr(i):=\sqcup_{j=1}^i Irr(i,j)$ and $\be_{m,i}: Irr(i)\to Irr(X(m,i))$ the disjoint union of the maps $\be_{m,i,j}$.
\medskip
{\bf Definition 8.5.} Let us assume the truth of Conjectures 2.1, 2.4. The element $\g L_{e;m}(C_{ij\g k})\in Irr(m+i)$ is defined as $\be^{-1}_{2m,m+i}(L_{e;m}(C_{ij\g k}))$.
\medskip
Particularly, $\forall \ \bar m\ge2m$ there exist elements $\g L_{e;m}(C_{ij\g k})(\bar m)$ of $Irr(X(\bar m,m+i))$. Warning: for $m_1\ne m_2$ the elements $ L_{e;m_1}(C_{ij\g k})$, $L_{e;m_2}(C_{ij\g k})$ belong to different series. Each of them gives rise its own series.
\medskip
{\bf Proposition 8.6.} Let the partition of $C_{ij\g k}$ be $i=d_1+...+d_l$. Then the partition of $L_{e;m}(C_{ij\g k})$ is $m+i=(d_1+1)+(d_2+1)+...+(d_l+1)+1+1+...+1$ ($m-l$ one's at the end).
\medskip
{\bf Proof.} This is equivalent to the following formula. Let $i=d'_1+...+d'_{l'}$ be the conjugate partition. Then the conjugate partition of $L_{e;m}(C_{ij\g k})$ is $m+i=m+d'_1+...+d'_{l'}$.
\medskip
Taking into consideration (5.14.2) we get that to prove the proposition, it is sufficient to prove
\medskip
{\bf Lemma 8.6.1.} For an element $(a_0:0:a_1:0:...:0:a_m)\in \La_1(2m)$, for any $\varkappa\ge0$ we have:
$$\rank (\g M(m)(a_0,...,a_m))^\varkappa=\rank (\g M(2m)(a_0:0:a_1:0:...:0:a_m))^{\varkappa+1}$$
\medskip
{\bf Proof.} A submatrix formed by all even rows and even columns of a matrix $\g N$ will be called its even submatrix, it is denoted by $\g N_{even}$. First, we have
\medskip
{\bf 8.6.2.} $\forall \varkappa$ the even submatrix of $(\g M(2m)(a_0:0:a_1:0:...:0:a_m))^{\varkappa}$ is equal to $(\g M(m)(a_0,...,a_m))^\varkappa$.
\medskip
This is proved immediately by induction along $\vk$. Really, for $\vk=1$ this follows from the definition of $\g M$. The induction step from $\vk$ to $\vk+1$ follows from the following fact. Let $\g N_1$, $\g N_2$ be two square matrices of odd size such that all odd lines of both $\g N_1$, $\g N_2$ are zero lines. Then all odd lines of $\g N_1\g N_2$ are zero lines, and $(\g N_1\g N_2)_{even}=(\g N_1)_{even}(\g N_2)_{even}$.
\medskip
Further, let $\g r=(\g r_1,\g r_2,...,\g r_{m-1})$ be a row matrix such that $\g r\cdot (\g M(m)(a_0,...,a_m))^\vk=0$. Let us prove that $\bar \g r:=(0,\g r_1,0,\g r_2,0,...,0,\g r_{m-1},0)$ satisfies $$\bar \g r\cdot (\g M(2m)(a_0:0:a_1:0:...:0:a_m))^{\vk+1}=0$$ Really, (8.6.2) implies that $$\bar \g r\cdot (\g M(2m)(a_0:0:a_1:0:...:0:a_m))^{\vk}=(*,0,*,0,*,...,0,*)$$
(this is because $$[\bar \g r\cdot (\g M(2m)(a_0:0:a_1:0:...:0:a_m))^{\vk}]_{1\be}=$$ $$\sum_{\al=1}^{2m-1}[\bar \g r]_{1\al}[(\g M(2m)(a_0:0:a_1:0:...:0:a_m))^{\vk}]_{\al\be}$$ and if $\al$ is odd then $[\bar \g r]_{1\al}=0$, if $\al$ is even then for even $\be$ (8.6.2) implies  $$[(\g M(2m)(a_0:0:a_1:0:...:0:a_m))^{\vk}]_{\al\be}=[(\g M(m)(a_0,...,a_m))^\varkappa]_{\frac\al2\frac\be2}$$ and the result follows from the condition $\g r\cdot (\g M(m)(a_0,...,a_m))^\vk=0$).
\medskip
Further, $(*,0,*,0,*,...,0,*) \cdot \g M(2m)(*:0:*:0:*:...:0:*)=0$.

Since for a set of linearly independent $\g r$'s the corresponding $\bar \g r$'s are also linearly independent, we get that $$\rank (\g M(m)(a_0,...,a_m))^\vk\ge \rank (\g M(2m)(a_0:0:a_1:0:...:0:a_m))^{\vk+1}.$$

Let now $\g r=(\g r_1,\g r_2,...,\g r_{2m-1})$ be a row matrix such that $$\g r\cdot (\g M(2m)(a_0:0:a_1:0:...:0:a_m))^{\vk+1}=0$$
Since the odd rows of $\g M(2m)(a_0:0:a_1:0:...:0:a_m)$ are 0, we have that $\g r_{0,even}:=(0,\g r_2,0,\g r_4,0,...,\g r_{2m-2},0)$ satisfies
$$\g r_{0,even}\cdot (\g M(2m)(a_0:0:a_1:0:...:0:a_m))^{\vk+1}=0$$
Let us denote $\g r_{0,even}\cdot (\g M(2m)(a_0:0:a_1:0:...:0:a_m))^{\vk}$ by $\g r':=(\g r'_1,\g r'_2,...,\g r'_{2m-1})$. We have $\g r'\cdot \g M(2m)(a_0:0:a_1:0:...:0:a_m)=0$. Since the submatrix of $\g M(2m)(a_0:0:a_1:0:...:0:a_m)$ formed by the even rows and all columns is of the maximal rank, we get $\g r'_{\be}=0$ for all even $\be$. But (8.6.2) implies that $\g r_{even}:=(\g r_2,\g r_4,...,\g r_{2m-2})$ satisfies $\g r_{even}\cdot (\g M(m)(a_0,...,a_m))^\vk= (\g r'_2,\g r'_4,...,\g r'_{2m-2})$, i.e. $\g r_{even}\cdot (\g M(m)(a_0,...,a_m))^\vk=0$. This gives us the opposite inequality $$\rank (\g M(m)(a_0,...,a_m))^\vk\le \rank (\g M(2m)(a_0:0:a_1:0:...:0:a_m))^{\vk+1}\ \ \ \  \square\  \square$$
\medskip
Let us denote by $i_e=i_e(m)$, $j_e=j_e(m)$, $\g k_e=\g k_e(m)$ the $i,j,\g k$-parameters of $\g L_{e;m}(C_{ij\g k})$, and by $\g d_e=\g d_e(m)=\g d_{e,i_e,j_e,\g k_e}$ its $\g d$-coefficient of (2.4.2). We have $i_e=m+i$.
\medskip
{\bf Conjecture 8.7.} $j_e=m-i=i_e-2i$, $\g d_e=\g d(C_{ij\g k})\binom{m-i}{j}$.
\medskip
{\bf Justification.} We have $\deg L_{e;m}(C_{ij\g k})=\deg C_{ij\g k}(m)= \g d(C_{ij\g k})\binom{m-i}{j}$. From another side, $\deg \g L_{e;m}(C_{ij\g k})(\bar m) = \g d_e \binom{\bar m-i_e}{j_e}$. Substituting $\bar m=2m$ we get
$$\g d_e\binom{m-i}{j_e}=\g d(C_{ij\g k})\binom{m-i}{j}\eqno{(8.7.1)}$$
This is one equation with two unknowns $\g d_e$, $j_e$, hence it has many solutions. The most natural solution is $\g d_e=\g d(C_{ij\g k})$, $j_e=j$
or $m-i-j$, but this is wrong. Really, let us consider $\g L_{e,i-1}(C_{111})\in Irr(i)$. (8.7.1) becomes
(here $i\mapsto1$, $m\mapsto i-1$) $$\g d_e\binom{i-2}{j_e}=i-2$$ which has 3 solutions ($\g d_e$, $j_e$ are integer $\ge1$):
\medskip
(1) $\g d_e=i-2$, $j_e=i-2$;
\medskip
(2) $\g d_e=1$, $j_e=i-3$;
\medskip
(3) $\g d_e=1$, $j_e=1$.
\medskip
According 8.6, the partition of $\g L_{e,i-1}(C_{111})$ is $i=2+1+1+...+1$, hence (see Table A3) $j_e$ cannot be 1. We let $i=4$.  Table A2.2 has no element corresponding $\g d_e=1$, $j_e=i-3$, hence for $i=4$ we have
$\g L_{e,i-1}(C_{111})=C_{423}$ of Table A2.2, and $j_e=m-i=i_e-2i$, $\g d_e=\g d(C_{ij\g k})\binom{m-i}{j}$ holds for this case.
\medskip
The following lemma gives us the image of $\g L_e$.
\medskip
{\bf Lemma 8.8.} Let $C_{ij\g k}\in Irr(i,j)$, and let $\pi(i,j)(C_{ij\g k})$ be a partition $i=d_1+...+d_l$. There exists $m$ such that $C_{ij\g k}(m)$ is an even lift of an element of $Irr(X(m/2, i-m/2))$ iff $l\ge(i+j)/2$.
\medskip
{\bf Proof.} According Proposition 8.1, a matrix $\g M(m)(a_0,...,a_m)\in C_{ij\g k}(m)$ belongs to the image of the lift $\iff$ its rank $\le m/2-1$. From another side, the rank of a generic matrix $\g M(m)(a_0,...,a_m)\in C_{ij\g k}(m)$ is $m-1-l$ (because of its Jordan form). Hence, $l\ge \frac m2 \iff C_{ij\g k}(m)$ belongs to the image of the lift. Finally, the minimal possible value of $m$ such that $C_{ij\g k}(m)\ne\emptyset$ is $i+j$. Comparing these results we get the desired. $\square$
\medskip
{\bf Example 8.9.} Let us consider the case $i=j=1$ (and hence $\g k=1$ --- the only possibility), $m=\g q-1$. We have $L_{e;\g q-1}(C_{111})\in Irr(X(2\g q-2),\g q))$. The corresponding element $\g L_{e;\g q-1}(C_{111})\in Irr(\g q)$  is
$C_{\g q,\g q-2,3}$ of Table A3 ($i$ of the corresponding entry of Table A3 is $\g q$ of the present notation), and $$L_{e;\g q-1}(C_{111})=C_{\g q,\g q-2,3}(2\g q-2)\eqno{(8.9.1)}$$
because according 8.6, its partition is $\g q=2+1+1+...+1$, and the value $\g k=3$ is only value of $\g k$ such that the partition of $C_{\g q,\g q-2,\g k}$ is $\g q=2+1+1+...+1$ (see Table A3). Values of $\g d$ and $\mu$ for $C_{\g q,\g q-2,3}$ from Table A3 are in concordance with Conjecture 8.7.
\medskip
{\bf Example 8.10.} Let us consider the case $i=2$. Components $\g L_{e;\g q-2}(C_{211})$, resp. $\g L_{e;\g q-2}(C_{221})$ are denoted by $C_{\g q,\g q-4,*_1}$, resp.  $C_{\g q,\g q-4,*_2}$ in Table A6.
\medskip
{\bf 9. Odd lift.}  There is an odd version of varieties $X(m,i)$. To avoid confusion of notations, we use $\g m$ instead of $m$, and $\g a_0,...,\g a_\g m$ instead of $a_0,...,a_m$.
\medskip
The initial variety (analog of $\n P^m$) is $\n P^1\times \n P^\g m$. Let

$$\bar \g a=((\g l_0:\g l_1); (\g a_0:...:\g a_\g m))\in \n P^1\times \n P^\g m$$

An analog of the matrix $\g M$ is the following $(2\g m \times 2\g m)$-matrix $\g M'_{odd}=\g M'_{odd}(\bar \g a)$:
\medskip
$\g M'_{odd}=\left(\matrix \Cal A_1&\Cal A_0\\ \g l_1 I_{\g m}& -\g l_0 I_{\g m} \endmatrix \right)$ where $\Cal A_1$, $\Cal A_0$ are the following $(\g m\times \g m)$-blocks (here $\ve=0,1 $ is the residue of $ \g m \mod 2$):
\medskip
$$\Cal A_1=\left(\matrix \g a_0 & \vdots& {\matrix \g a_2&\g a_4&...&\g a_{\g m-\ve}&0&...&0\endmatrix} \\ \cdots \cdots \cdots &&\cdots \cdots \cdots \cdots \cdots \cdots \cdots \cdots \cdots  \\ 0_{(\g m-1)\times1}&\vdots&\g M(\g m)(\g a_0,...,\g a_{\g m}) \endmatrix \right)\eqno{(9.1)}$$ of diagonal block sizes $1, \g m-1$,
$$\Cal A_0=\left(\matrix \g M(\g m)(\g a_0,...,\g a_{\g m}) &  \vdots&0_{(\g m-1)\times1}\\ \cdots \cdots \cdots \cdots \cdots \cdots \cdots \cdots \cdots  \cdots && \cdots \cdots \cdots  \\ {\matrix0&...&0& \g a_\ve&\g a_{2+\ve}&...&\g a_{\g m-2}\endmatrix} & \vdots&\g a_ {\g m}\endmatrix \right)\eqno{(9.2)}$$ of diagonal block sizes $ \g m-1,1$. $\Cal A_0$ is the original Hurwitz matrix ([H], p. 274).
\medskip
For the reader's convenience, we give explicitly the matrix $\g M'_{odd}$ for $\g m=4$:

$$\g M'_{odd}=\left(\matrix \g a_0 &  \g a_2&\g a_4&0&\g a_1&\g a_3&0&0\\ 0&\g a_1 &\g a_3&0& \g a_0 &  \g a_2&\g a_4&0 \\ 0& \g a_0 &  \g a_2 &\g a_4& 0&\g a_1&\g a_3&0\\ 0&0&\g a_1&\g a_3 &0& \g a_0 &  \g a_2&\g a_4\\ \g l_1&0&0&0&-\g l_0 &0&0&0\\ 0&\g l_1&0&0&0&-\g l_0 &0&0\\ 0&0&\g l_1&0&0&0&-\g l_0 &0\\ 0&0&0&\g l_1&0&0&0&-\g l_0 \endmatrix \right)$$

\medskip
{\bf Definition 9.3.} $Ch(\g M'_{odd})(\bar \g a)$ --- the odd characteristic polynomial of $\g M'_{odd}(\bar \g a)$ --- is $|\g M'_{U;odd}(\bar \g a)|$ where
\medskip
$\g M'_{U;odd}(\bar \g a):=\left(\matrix \Cal A_1-\frac{U}{\g l_0} I_{\g m}&\Cal A_0\\ \g l_1 I_{\g m}& -\g l_0 I_{\g m} \endmatrix \right)$  \ \ ( $\Cal A_1$, $\Cal A_0$ as above).
\medskip
We define $D_{odd}(\g m,i)$ by the following formula where $Ch(\g M'_{odd})(\bar \g a)$ is considered as a polynomial in $U$:
$$Ch(\g M'_{odd})(\bar \g a)=D_{odd}(\g m,0)+D_{odd}(\g m,1)U+D_{odd}(\g m,2)U^2+...+U^\g m$$

$D_{odd}(\g m,i)$ are bihomogeneous polynomials in $((\g l_0, \g l_1)$; $(\g a_0, ... ,\g a_\g m))$ of bidegree $(\g m-i), (\g m-i)$.
\medskip
{\bf Definition 9.4.} $X_{odd,S}(\g m,i)\subset \n P^1(\n C)\times \n P^\g m(\n C)$ is a projective scheme corresponding to the ideal generated by the first $i$ polynomials $D_{odd}(\g m,i')$, $i'=0,1,...,i-1$ (here the subscript $S$ means scheme). Its support (i.e. the set of closed points = the set of zeroes of $D_{odd}(\g m,i')$) is denoted by $X_{odd}(\g m,i)$.
\medskip
{\bf Conjecture 9.5.} $X_{odd}(\g m,i)$ is a complete intersection and hence has codimension $i$ in $\n P^1\times \n P^\g m$.
All its irreducible components have the same codimension.
\medskip
Let us show that there exists a relation between $X_{odd}(\g m,i)$ and $X(m,i)$. Let $m=2\g m+1$ be odd. We denote by $\psi$ the Segre embedding $\n P^1\times \n P^{\g m}\hookrightarrow \n P^m$.
\medskip
{\bf Lemma 9.6.} Let $\g M(m)(a_0,...,a_m)$ have the rank $\le \g m$. Then

$\exists$ $\bar \g a=((\g l_0:\g l_1);  (\g a_0:...:\g a_{\g m}))\in \n P^1\times \n P^\g m$ such that $(a_0,...,a_m)=\psi(\bar \g a)$, i.e.
$a_{2j}=\g l_1 \g a_j$, $a_{2j+1}=\g l_0 \g a_j$, $j=0,...,\g m$.
\medskip
{\bf Proof.} Similar to the proof of Proposition 8.1. If $a_0\ne0$ then we consider the $2,4,6,...,2\g m-2, 2\g m-1,2\g m$-th rows of $\g M(m)$. They are linearly dependent. The coefficients of the linear dependence for the $2,4,6,...,2\g m-2$-th rows are 0 (because the leftmost non-zero element of the $2\vk$-th row of $\g M(m)$ is $a_0$, and it is on the $\vk$-th position). Hence, $2\g m-1$-th and $2\g m$-th rows are linearly dependent, hence the lemma. If $a_0=0$ and $a_1\ne0$ then we consider the $1,3,5,...,2\g m-1,2\g m$-th rows of $\g M(m)$. They are linearly dependent. The coefficients of the linear dependence for the $1,3,5,...,2\g m-1$-th rows are 0 (because the leftmost non-zero element of the $2\vk-1$-th row of $\g M(m)$ is $a_1$, and it is on the $\vk$-th position). This implies that the $2\g m$-th row is the zero row, hence the lemma. Finally, if $a_0=a_1=0$ then we repeat the process considering the case $a_2\ne0$ etc. $\square$
\medskip
{\bf Lemma 9.7.} Let $\bar \g a\in \n P^1\times \n P^{\g m}$. Then the characteristic polynomial of $\g M(m)(\psi(\bar \g a))$ is given by the formula $$Ch(\g M(m)(\psi(\bar \g a)))=\pm Ch(\g M'_{odd}(\bar \g a))\cdot U^\g m\eqno{(9.7.1)}$$

{\bf Proof.} Elementary matrix transformations. We denote the $\vk$-th row of $\g M(m)(\psi(\bar \g a))-U\cdot I_{m-1}$ by $r_\vk$.
We apply the following elementary
transformations to rows of $\g M(m)(\psi(\bar \g a))-U\cdot I_{m-1}$ (denoted by the symbol $\mapsto$):
$$\hbox{First --- for even $\vk$: $r_\vk \mapsto (\g l_0 r_\vk-\g l_1 r_{\vk-1})/U$; second --- for odd $\vk$: $r_\vk \mapsto r_\vk/\g l_0$.}\eqno{(9.7.2)}$$
We denote the obtained matrix by $\g M_{U;odd}=\g M_{U;odd}(\bar \g a)$. It is the following:
\medskip
$$\g M_{U;odd}=\left(\matrix \g a_0-U/\g l_0& \g a_1&\g a_2& ... & \g a_{\g m}&0& 0&  ... &0\\ \g l_1&-\g l_0&0& 0& 0& 0&0&  ... &0 \\  0&  \g a_0& \g a_1-U/\g l_0&\g a_2& ... & \g a_{\g m}& 0&  ... &0 \\ 0& 0&\g l_1&-\g l_0& 0&0& 0&  ... &0 \\ ... & ... &... &... &... &... &... &... &... \\  0& 0&... & 0& \g a_0& \g a_1& ... & \g a_{\g m-1}-U/\g l_0& \g a_{\g m} \\ 0& 0&0&0& 0&  ... &0& \g l_1&-\g l_0 \endmatrix \right)$$
The matrix $\g M'_{U;odd}$ is a permutation $\left(\matrix 1&2&...&\g m&\g m+1&\g m+2&...&2\g m\\ 1&3&...&2\g m-1& 2&4&...&2\g m\endmatrix \right)$ of both rows and columns of $\g M_{U;odd}$, hence the lemma (the factor $U^\g m$ in (9.7.1) appears because of division by $U$ in (9.7.2)). \ \ $\square$
\medskip
{\bf Corollary 9.8.} $X_{odd}(\g m,i)=\psi(\n P^1\times \n P^{\g m})\cap X(2\g m+1, \g m+i)$ (intersection in $\n P^m$; up to $\psi$).
\medskip
{\bf Remark 9.9.} It is more convenient to consider $X_{odd}(\g m,i)$ as subvarieties of $\n P^1\times \n P^{\g m}$ and not of $\n P^{2\g m+1}$. For example, its Chow class in $\n P^1\times \n P^{\g m}$ is a more thin invariant than the degree in $\n P^{2\g m+1}$.
\medskip
Since (conjecturally) $\dim X_{odd}(\g m,i)=\dim X(2\g m+1, \g m+i)$ we have that any irreducible component of $X_{odd}(\g m,i)$ is an irreducible component of $X(2\g m+1,\g m+i)$. This means that we have an inclusion of the set of irreducible components: $Irr(X_{odd}(\g m,i))\to Irr(X(2\g m+1,\g m+i)) $.
\medskip
{\bf Definition 9.10.} This map is called the odd lift, and it is denoted by $L_o$.
\medskip
{\bf Remark 9.10a.} A reader can think that always $\psi(X_{odd}(\g m,i))=X(2\g m+1, \g m+i)$, but this is not true. Really, let $(b_0:...:b_m) \in X(2\g m+1, \g m+i)$. This means that the first $\g m+i$ coefficients of the characteristic polynomial of $\g M(m)(b_0:...:b_m)$ are zeros, but this does not mean that the rank of $\g M(m)(b_0:...:b_m)$ is $\le \g m$. Hence, it can happen that $(b_0:...:b_m) \not\in \psi(\n P^1\times \n P^{\g m})$.
\medskip
The image of $L_o$ is described by the following
\medskip
{\bf Lemma 9.11.}  Let $C_{ij\g k}\in Irr(i)$ be a type of irreducible components, and $\{i=d_1+...+d_l\}$ its partition $\pi_i(C_{ij\g k})$. There exists $m$ such that $C_{ij\g k}(m)$ is an odd lift of an element of $Irr(X_{odd}(\frac{m-1}2, i-\frac{m-1}2))$ iff $l\ge(i+j-1)/2$.
\medskip
{\bf Proof.} It repeats the proof of Lemma 8.8.
According Lemma 9.6, a matrix $\g M(m)(a_0,...,a_m)\in C_{ij\g k}(m)$ belongs to the image of the odd lift $\iff$ its rank $\le \frac{m-1}2$. From another side, the rank of a generic matrix $\g M(m)(a_0,...,a_m)\in C_{ij\g k}(m)$ is $m-1-l$ (because of its Jordan form). Hence, $l\ge \frac {m-1}2 \iff C_{ij\g k}(m)$ belongs to the image of the odd lift. Finally, the minimal possible value of $m$ such that $C_{ij\g k}(m)\ne\emptyset$ is $i+j$. Comparing these results we get the desired. $\square$
\medskip
{\bf 10. Odd lift and first power of Carlitz modules.}
\medskip
Recall that we consider only the case $q=2$. We have the following relation between $\Cal M(\g a_*,1,\g m)$ from (0.1.4a) and $\g M'_{odd}$ of Section 9:
\medskip
{\bf 10.1.} $\Cal M(\g a_*,1,\g m)=\Cal A_0t-\Cal A_1$ where $\Cal A_1$, $\Cal A_0$ are from (9.1), (9.2).
\medskip
Let $H_{i\g j}:=H_{i\g j,21}(\g m)(\g a_0,...,\g a_\g m)$ be from (0.2.2). We have
\medskip
{\bf Proposition 10.2.} $$D_{odd}(\g m,i)=\sum_{\g j=0}^{\g m-i} \pm H_{i\g j} \g l_0^{\g m-i-\g j}\g l_1^\g j\eqno{(10.2.1)}$$
\medskip
{\bf Proof.} Immediately follows from a well-known equality
$\left|\matrix \Cal A_1&\Cal A_0\\ \g l_0\cdot I_n & \g l_1\cdot I_n  \endmatrix \right|\sim |\g l_1 \Cal A_1+\g l_0 \Cal A_0|$ where $\sim$ means the equality up to signs of coefficients of monomials $\g l_0^i\g l_1^{n-i}$ for any $n\times n$-matrices $\Cal A_1$, $\Cal A_0$. $\square$
\medskip
{\bf 10.3.} Let $pr: \n P^1\times \n P^{\g m}\to \n P^{\g m}$ be the projection. Let us fix a number $i_0$. Proposition 10.2 implies
\medskip
{\bf Corollary 10.4.} If Conjecture 0.2.4 holds for $n=1$, $m=\g m$, any $\be$ and any $i\le i_0-1$ then
\medskip
$pr^{-1}(X(\g m,i_0))\subset X_{odd}(\g m,i_0)$.
\medskip
Particularly, in this case for $Z$ an irreducible component of $X(\g m,i_0)$ there exists

$L_o(\n P^1\times Z)\in Irr(2\g m+1,\g m+i_0)$.
\medskip
Conjecture 0.2.4 holds for $i=0$, $n=1$. This is a particular case of [GL16], Theorem III. Since its proof is long and complicated, we give here a short proof of this particular case.
\medskip
{\bf Proposition 10.5.} $H_{0\g j,21}(\g m)=\pm \g a_\g j D(\g m,0)$.
\medskip
{\bf Proof.} We denote the result of substitution $U=0$ in $\g M_{U;odd}$ by $\g M_{odd}(0)$. According (10.2.1), the formula (10.5) is equivalent to $$|\g M_{odd}(0)|=\pm |\g M(\g m)(\g a_0,...,\g a_{\g m})|\cdot(\sum_{\g j=0}^{\g m} \g a_\g j\g l_0^{\g m-\g j}\g l_1^\g j)$$ (recall that $D(\g m,0)= |\g M(\g m)(\g a_0,...,\g a_{\g m})|$\ ).
\medskip
Let $V$ be a column matrix $\left(\matrix \g l_0^{2\g m-1}& \g l_0^{2\g m-2}\g l_1& \g l_0^{2\g m-3}\g l_1^2& ... & \g l_1^{2\g m-1}\endmatrix \right)^{tr}$.
\medskip
If $\sum_{\g j=0}^{\g m} \g a_\g j\g l_0^{\g m-\g j}\g l_1^\g j=0$ then $\g M_{odd}(0) \cdot V=0$, hence in this case $|\g M_{odd}(0)|=0$. This means that
$\sum_{\g j=0}^{\g m} \g a_\g j\g l_0^{\g m-\g j}\g l_1^\g j$ is a factor of $|\g M_{odd}(0)|$. Comparing bidegrees we see that the second factor is a homogeneous polynomial in $\g a_*$ of degree $\g m-1$. To find it,
we substitute $\g l_1=0$. After elimination of rows and columns containing $\g l_0$ in $\g M_{odd}(0)|_{\g l_1=0}$ we get the
matrix $\Cal A_1$ from 9.1.
This gives us the desired. $\square$
\medskip
 Slightly modifying the notations of (5.11.8), we denote by $L_\g m$ the class of a hyperplane in $\n P^\g m$ in its Chow ring $CH(\n P^\g m)$. We have $CH(\n P^1\times \n P^\g m)$ is generated by $L_1$, $L_\g m$.
\medskip
{\bf Corollary 10.6.} $X_{odd}(\g m,1)$ is a union of 2 irreducible components, one of them (denoted by $X_{odd}(\g m,1)_1$) is
$pr^{-1}(X(\g m,1))$. Its Chow class is $(\g m-1)L_\g m$. Another component (denoted by $X_{odd}(\g m,1)_2$) has the equation  $\sum_{\g j=0}^{\g m} \g a_\g j\g l_0^{\g m-\g j}\g l_1^\g j=0$. Its Chow class is $\g m L_1+L_\g m$.
\medskip
Study of $X_{odd}(\g m,i)$, i.e. finding of their irreducible components, their Chow classes etc. is a subject of further research.
\medskip
{\bf Remark 10.7.} Conjecture 0.2.4 states that the set of zeroes of $H_{i'\g j}$, where $i'=0,...,i-1$, is $X(\g m,i)$ as a set of points, but not
as a scheme. This means that irreducible components of $X(\g m,i)$ can have different multiplicities being considered as
\medskip
(a) Intersection of hypersurfaces $\{D(\g m,i')=0\}$, $i'=0,...,i-1$;
\medskip
(b) Intersection of hypersurfaces $\{H_{i'\g j}=0\}$, $j=0,...,i-1$.
\medskip
See, for example, [GL16], Tables 9.7.7, 9.7.8, lines: Multiplicities in $X(m,i)$; Multiplicities in $X(2,1,m,i)$.
We see that these multiplicities are different.
\medskip
This phenomenon must be taken into consideration when we shall count the Chow class of $pr^{-1}(X(\g m,i))$ for $i>1$: it is not equal to
$(\g m-1)(\g m-2)\cdot...\cdot(\g m-i)L_\g m$.
\medskip
{\bf 11. More examples.}
\medskip
{\bf 11.1. Components having Jordan form $2+1+1+...+1$.}
\medskip
Here we describe these components using results of Sections 8 --- 10. Let us fix $C_{ij\g k}\in Irr(i)$ such that $\pi(i,j)(C_{ij\g k})=\{i=2+1+1+...+1\}$. We have $j\le i-1$.  Since the partition $i=2+1+1+...+1$ has $l=i-1$ we get that the condition $l\ge(i+j-1)/2$
of Lemma 9.11 is satisfied, hence for all such $C_{ij\g k}$ there exists $m$ --- namely, $m=2i-1$ --- such that $C_{ij\g k}(2i-1)$
is the odd lift of an element of $Irr(X_{odd}(i-1,1))$.

Since $Irr(X_{odd}(i-1,1))$ consists of two elements, we get that the above conjectures imply that $X(2i-1,i)$ has two irreducible components of Jordan type $2+1+1+...+1$. Hence, the same is true for $Irr(i)$, $i\ge 3$.
\medskip
Let us find degrees of these components, and identify them with entries of Tables A2.2, A3. The direct image map $CH(\psi)_*: CH(\n P^1\times \n P^\g n) \to CH(\n P^{2\g n+1})$ is the following:
$$CH(\psi)_*(L_1)=L_{2\g n+1}^{\g n+1}, \ \ \ \  CH(\psi)_*(L_\g n)=\g nL_{2\g n+1}^{\g n+1}$$
Applying these formulas for $\g n=i-1$ to $X_{odd}(i-1,1)_1$, $X_{odd}(i-1,1)_2$, we get that $$\deg L_o(X_{odd}(i-1,1)_1)=(i-1)(i-2), \ \ \ \ \  \deg L_o(X_{odd}(i-1,1)_2)=2(i-1)$$
Comparing these values with results of computations (see Table A2.2), we can strongly conjecture that $$L_o(X_{odd}(i-1,1)_1)=C_{i,i-2,3}(2i-1), \ \ \ \ L_o(X_{odd}(i-1,1)_2)=C_{i,i-1,1}(2i-1)$$

\newpage
{\bf 11.2. Remarks.}
\medskip
{\bf Remark 11.2.1.} We have: $L_{e;i-1}(C_{111})=C_{i,i-2,3}(2i-2)$ (see 8.9.1) having degree $i-2$, and
\medskip
$L_o(X_{odd}(i-1,1)_1)=C_{i,i-2,3}(2i-1)$ having degree $(i-1)(i-2)$ are two different elements of the same series $C_{i,i-2,3}(m)$ obtained by two different lifts.
\medskip
{\bf 11.2.2.} Computations (see Table A2.2) show that (conjecturally) $\g d_{i,i-1,1}=2i-2$, $\mu_{i,i-1,1}=i!/2$. Formula (2.8) for $j=i-1$ implies that $C_{i,i-1,1}$  is the only component having $j=i-1$. So, we get evidence that $\#Irr(i,i-1)=1$.
\medskip
{\bf 11.2.3.} We see that although Proposition 10.2 and Conjecture 0.2.4 show that $X_{odd}(\g m,i)$ are tightly related with $X(\g m,i)$,
complementary factors (like $\sum_{j=0}^{\g m} \g a_j\g l_0^{\g m-j}\mu^j$ for $i=1$) can appear, and they give new irreducible components in
$X_{odd}(\g m,i)$ which do not exist in $X(\g m,i)$. Therefore, study of $X_{odd}(\g m,i)$ is a non-trivial subject.
\medskip
{\bf 11.2.4.} Why there is a relation between $\n P^1\times \n P^\g m \hookrightarrow \n P^m$, and the first tensor power of Carlitz modules? Do exist objects related to the $n$-th tensor power of Carlitz modules?
\medskip
{\bf Example 11.3. Components having $\g d_{ij*}=1$. }
\medskip
They are obtained by consecutive application of the operator $L_e$ to $C_{ii1}(2i)$. Really, if $\g d_{ij*}=1$ and $m=i+j$, then (see 8.7) $\g d_e=1$ as well, and we can repeat the construction. Let us give its independent description and a proof of 8.7 for this case.
\medskip
Let us fix $k, \eta\ge1$ as in (6.4). It is convenient to let $\bar k:=k+1$. We choose $m=2^{\bar k}\eta$. We consider a linear subspace $\La_{\bar k}(m)\subset \n P^{m}$ defined by
the equations $\{a_\vk=0$ if $\vk$ is not a multiple of $2^{\bar k}\}$ (generalization of $\La_1(m)$ of 8.2). We have $\La_{\bar k}(m)\subset X(m,m-\eta)$
(consecutive application of Lemma 8.3), and because $\dim \La_{\bar k}(m)=\dim X(m,m-\eta)=\eta$,
it is an irreducible component of $X(m, m-\eta)$. We denote it by $C_1(\bar k,\eta)$.
\medskip
For a generic point $(a_0:...:a_{m})\in \La_{\bar k}(m)$ we have
\medskip
{\bf Lemma 11.3.1.} For $\vk \le\bar k$ the rank of $[\g M(m)(a_0:...:a_m)]^\vk$ is $2^{\bar k - \vk}\eta-1$, and for $\vk \ge\bar k$ the rank of $[\g M(m)(a_0:...:a_m)]^\vk$ is $\eta-1$.
\medskip
{\bf Proof.} This lemma can be proved by a consecutive application of Proposition 8.6; let us give a short proof. For $\vk \le\bar k$
the only non-zero rows of $[\g M(m)(a_0:...:a_m)]^\vk$ are rows $\be\cdot 2^\vk$, $\be=1,\dots, 2^{\bar k - \vk}\eta-1$ --- this is
easily proved by induction by $\vk$. For $\vk \ge\bar k$ the only non-zero rows
of $[\g M(m)(a_0:...:a_m)]^\vk$ are rows $\be\cdot 2^{\bar k}$, $\be=1,\dots, \eta-1$.  $\square$
\medskip
{\bf Corollary 11.3.2.} The conjugate of the partition $\pi_{m,m-\eta}(C_1(\bar k,\eta))$ is $$i=m-\eta=2^{k}\eta+2^{k-1}\eta+...+2\eta+\eta$$
\medskip
{\bf Proof.} The rank of the $\vk$-th power of the 0-Jordan part of $\g M(m)(a_0:...:a_m)$ is $(2^{\bar k - \vk}\eta-1)-(\eta-1)=(2^{\bar k - \vk}-1)\eta$, hence (8.6.1) implies the proposition. $\square$
\medskip
Let us change notations: we denote $2^{\bar k}\eta$ by $m_0$.
Conjectures 2.4 imply that for any fixed $\bar k\ge 1$, $\eta$, for any $m\ge m_0$ there exists an irreducible component
$C_{(2^{\bar k}-1)\eta,\eta,*}$ of $X(m,(2^{\bar k}-1)\eta)$ of degree $\binom{m-(2^{\bar k}-1)\eta}{\eta}$ and of partition conjugate to
$i=(2^{\bar k}-1)\eta=2^{k}\eta+2^{k-1}\eta+...+2\eta+\eta$. For $k=0$ this is the component $C_{ii1}$.
\medskip
The opposite particular case is $\eta=1$. This component is $\im \vf(T,w)$ where $T$ is the complete tree of depth $k$ and $w$ is any its weight. We denote it by $C_{2^{\bar k}-1,1,*}$, its degree is $m-(2^{\bar k}-1)$ and the conjugate of its Jordan form is $2^{\bar k}-1=2^{k}+2^{k-1}+...+2+1$.
\medskip
{\bf Conjecture 11.3.3.} These are the only irreducible components having $\g d(C_{ij\g k})=1$.
\medskip
{\bf Justification.} Let us consider the case $\eta=1$ and a component $C_{2^{\bar k}-1,1,*}$. We have $m=2^{\bar k}$. It has degree 1, hence it is a line $\n P^1\subset \n P^{2^{\bar k}}$. Its parametric equation is $(a_0:0:...:0:a_{2^{\bar k}})$. Let us show that for other $m$ there is no components of degree 1.
\medskip
{\bf Lemma 11.3.3.1.} Let us consider a point $(a_0:0:...:0:a_{m})\in \n P^m$ and its matrix $\g M(m)$. Its characteristic polynomial
$Ch(\g M(m))$ is equal to $\pm U^{m-1}$ iff $m=2^{\bar k}$.
\medskip
{\bf Proof} is straightforward. Let $m=2^\vk m_1$, where $m_1$ is odd. If $m=m_1$ then $|\g M(m)|=\pm (a_0a_m)^{(m-1)/2}$. If $m$ is even then $Ch(\g M(m))=\pm U^{m/2}Ch(\g M_{even}(m))$. Now we use induction by $\vk$. $\square$
\medskip
This means that for $m=2^{\bar k}$ the above $\n P^1$ is contained in $X(m,m-1)$. Since dim $X(m,m-1)=1$, it is its irreducible component.
\medskip
Let us prove a version of (7.1) for a case of $C_{2^{\bar k}-1,1,*}(2^{\bar k})$ and $dd=1$:
\medskip
{\bf Proposition 11.3.4.} $\nu(C_{2^{\bar k}-1,1,*}(2^{\bar k})\times \n P^1)\subset X(2^{\bar k}+1,2^{\bar k}-1)$.
\medskip
{\bf Proof.} $\nu(C_{2^{\bar k}-1,1,*}(2^{\bar k})\times \n P^1)\subset \n P^{2^{\bar k}+1}$ is the set of elements $a_0: ... : a_{2^{\bar k}+1}$ such that $a_i=0$ for
$i\ne 0,1, 2^{\bar k},2^{\bar k}+1$ and
\medskip
{\bf 11.3.4.1.} The rank of $\left(\matrix a_1 & a_{2^{\bar k}+1} \\ a_0 & a_{2^{\bar k}}\endmatrix \right)$ is 1.
\medskip
Particularly, $\nu(C_{2^{\bar k}-1,1,*}(2^{\bar k})\times \n P^1)$ is an odd lift of $\n P^1 \times C_{2^{k}-1,1,*}(2^{k})$. The characteristic polynomial is given by the formula:
$$Ch(\g M(2^{\bar k}+1)(a_0:a_1:0:0:...:0:a_{2^{\bar k}}: a_{2^{\bar k}+1})=\pm (a_1+a_{2^{\bar k}})U^{2^{\bar k}-1}\pm U^{2^{\bar k}}$$
It follows immediately by induction by $\bar k$. The induction step is the following:
\medskip
$$Ch(\g M(2^{\bar k}+1)(a_0:a_1:0:0:...:0:a_{2^{\bar k}}: a_{2^{\bar k}+1})=$$
$$=\pm U^{2^{k}} Ch(\g M(2^{k}+1)(b_0:b_1:0:0:...:0:b_{2^{k}}: b_{2^{k}+1})$$ where $b_0=a_0$, $b_1=a_1$,
$b_{2^{k}}=a_{2^{\bar k}}$, $b_{2^{k}+1}= a_{2^{\bar k}+1}$.
\medskip
Proof of this formula is similar to the proof of Lemma 9.7. We apply the following elementary transformations
to the rows $\g r_i$ of $\g M(2^{\bar k}+1)(a_0:a_1:0:0:...:0:a_{2^{\bar k}}: a_{2^{\bar k}+1})-U\cdot I_{2^{\bar k}}$ (denoted by the symbol $\mapsto$):
\medskip
First --- for even $j$: $\g r_j \mapsto (\g r_j-(a_0/a_1) \g r_{j-1})/U$; second --- for $j$ of the form $j=4k+3$:
$\g r_j \mapsto \g r_j+a_1 \g r_{(j+1)/2}+a_{2^{\bar k}+1} \g r_{(j+1)/2+2^{k}}$.
\medskip
We get a matrix whose even columns contain only one non-zero element at the diagonal which is equal to $-1$. Eliminating the even lines
and columns we get the matrix $\g M(2^{k}+1)(b_0:b_1:0:0:...:0:b_{2^{k}}: b_{2^{k}+1})-U\cdot I_{2^{k}}$. $\square$
\medskip
{\bf Proposition 11.4.} $X(m,1)$ is an irreducible hypersurface.
\medskip
{\bf Proof.}\footnotemark \footnotetext{The authors are grateful to A. Esterov who indicated them this proof.} More generally, any resultantal variety is irreducible. Let $P_1=\sum_{i=0}^m a_ix^i$, $P_2=\sum_{i=0}^n b_ix^i$ be two polynomials, $R\in \n Z[a_0,...,a_m,b_0,...,b_n]$ their resultant and $\g R\subset \n P^m\times \n P^n$ the set of zeroes of $R$. Let us prove that $\g R$ is irreducible. We consider a flag variety $F\subset  \n P^m\times  \n P^n\times  \n P^1 $ defined as follows:
$$(a_0,...,a_m; \ b_0,...,b_n; \ \la, \mu)\in F \iff \sum_{i=0}^m a_i\la^i\mu^{m-i}= \sum_{i=0}^n b_i\la^i\mu^{n-i}=0$$
Let $\pi_1,\pi_2$ be projections of $ \n P^m\times  \n P^n\times  \n P^1 $ to $ \n P^m\times  \n P^n$, $  \n P^1 $ respectively. The fiber of $\pi_2$ is a product $H_1\times H_2\subset  \n P^m\times  \n P^n$ where $H_1\subset  \n P^m$, $H_2\subset  \n P^n$ are hyperplanes, hence $F$ is irreducible. Since $\g R=\pi_1(F)$, we get that $\g R$ is also irreducible. $\square$
\medskip
\medskip
{\bf Appendix. Tables.}
\medskip
{\bf A1.} The below table illustrates Conjecture 0.2.4. Let us recall its notations. All polynomials are over $\n C$. We have: $a_0,\dots, a_m$ are abstract variables, $q=2$, $n\ge0$ a number such that $\kappa:=\frac{m+n}{q-1}-1$ is integer (i.e. for our case $q=2$ we have $\kappa$ is always integer, $\kappa=m+n-1$), and a $\kappa\times\kappa$-matrix $\Cal M_{nt}(a_*,n,\kappa)$ is given by the same formula (0.1.4a) as the matrix $\Cal M(a_*,n,\bar \kappa)$, with the difference that for the $\Cal M_{nt}$-case we have $1\le \g i, \g j \le \kappa$. Further, $Ch(\Cal M_{nt}(a_*,n,\kappa))$ - a version of its characteristic polynomial - is given by (0.2.1), and $H_{i\g j,qn}(m)$ are defined by (0.2.2). Finally, $D(m,i)=H_{i0,20}(m)$, see also (1.2.1).
\medskip
We consider only the case $n=1$, i.e. $m=\kappa$. For some values of $m,i,\g j$ the table gives the minimal value of $\vk$ such that $H_{i\g j,21}(m)^\vk$ belongs to the ideal generated by $D(m,0),\dots, D(m,i)$.
\medskip
Recall that there is a symmetry between $H_{i\g j,21}(m)$ and $H_{i,m-i-\g j,21}(m)$. Further, Conjecture 0.2.4 is proved for the simple cases (a) $i=0$, see (10.5); (b) $\g j=0$ ([GL16], 9.12); (c) $i=m-2$, $\g j=1$ ([ELS]); in all these cases $\vk=1$, and most likely for $n=1$ these are (up to the symmetry) the only cases having $\vk=1$. These cases are omitted. Hence, the table covers (for $n=1)$ all cases for $m\le6$, and all cases for $m=7$, $i\le2$. Existence of $\vk$ shows that Conjecture 0.2.4 is true for these cases.
\medskip
We got the below results using computer-assistant calculations. We used Magma computer system in order to find values of $\vk$. The Magma system has an option permitting to determine whether a polynomial $\g H\in \n C[x_1,\dots,x_n]$ belongs to an ideal $I\subset \n C[x_1,\dots,x_n]$, or not. To define an ideal $I$, Magma system requieres to input its set of generators as a $\n C[x_1,\dots,x_n]$-module.
\medskip
For simplicity, we did not use cycles (one program for all entries of Table A1), but we wrote Magma operators separately for consecutive entries (there are only 15 cases, it is more rapidly to input data separately for all these cases than to write a quadruple cycle program). Namely, for all combinations of $m, \ i, \ \g j, \ \vk$ of lines of A1 (15 combinations of $m, \ i, \ \g j$) we inputted the corresponding polynomials $D(m,i')$ (where $0\le i'\le i$) which are generators of an ideal, and the polynomial $H_{i\g j,21}^\vk$. The program replied whether $H_{i\g j,21}^\vk$ belongs to $<D(m,0),\dots,D(m,i)>$, or not. For the last case $m=7$, $i=3$, $\g j=1$, $\vk=6$ the calculation took too much time, so we do not know the answer.
\medskip
We see that within the limits of the table the value of $\vk$ does not depend on $\g j$, except the above simple cases. Particularly, we can conjecture that for $n=1$, $i=1$, $\forall \ m$, $\forall \ \g j\ne 0, \ m-1$ we have $\vk=2$.
\medskip
We did not consider the case $n>1$, because of Proposition 7.2.1.
\medskip
\settabs 10 \columns
\+&&&{\bf Table A1}\cr
\+ $m$ &$i$&$\g j$ &&& $\vk$\cr
\medskip
\+4&1&1    &&&2\cr
\medskip
\+5&1&1 - 2&&&2\cr
\medskip
\+5&2&1    &&&3\cr
\medskip
\+6&1&1 - 2&&&2\cr
\medskip
\+6&2&1 - 2&&&4\cr
\medskip
\+6&3&1    &&&6\cr
\medskip
\+7&1&1 - 3&&&2\cr
\medskip
\+7&2&1 - 2&&&4\cr
\medskip
\+7&3&1    &&&? \ ($\ge6$)\cr
\medskip
\medskip
{\bf A2. Table of irreducible components.} The below table A2.2 gives us conjectural properties of irreducible components of $X(m,i)$ for $i\le6$. It generalizes the corresponding tables of [GL16], 9.7.7, 9.7.8, 9.7.10. The first column of the table contains the designation of a component $C_{ij\g k}(m)$ (notation $m$ is omitted). The second column contains the value of $j$, the third column contains the degree of this component as a function of $m$. The 4-th -- 6-th columns contain, respectively, the multiplicity, partition of the Jordan form, and the forest of the component (they do not depend on $m$). The weight of the forest is not given, except the components $C_{618a}(m)$, $C_{618b}(m)$, because for all other entries of the table this weight is unique up to the action of $\g G(F)\rtimes\Aut(F)$. Within the limit of the table, all irreducible components are defined over $\n Q$, except the same components $C_{618a}(m)$, $C_{618b}(m)$. They are defined over $\n Q(\sqrt{-1})$ and they are conjugate over $\n Q$.
\medskip
Hence, in order to get information on $X(m,i)$, $i\le6$, we should consider lines corresponding to this value of $i$ (see the first column). If $m\ge2i$ we consider all these lines, if $m<2i$ we consider lines corresponding to $j=1, \dots, m-i$ (see the second column). Any line of this set is in 1 -- 1 correspondence with irreducible components of $X(m,i)$, and the 3-rd -- 6-th columns give us invariants of this component. See also Example 2.5.2.
\medskip
{\bf Comparison with results of [GL16].} Columns of the tables of [GL16], Section 9.7 are lines of the tables of the present paper (sorry). $C_{21}$, $C_{31}$, $C_{32}$, $C_{33}$, $C_{i1}$, $X(m,i)_{pr}$ of [GL16] are respectively $C_{211}(m)$, $C_{311}(m)$, $C_{312}(m)$, $C_{321}(m)$, $C_{i11}(m)$, $C_{ii1}(m)$ of the present paper. $c$ of Conjecture 9.7.10 of [GL16] and $c_{ijk}$ below --- are $\g d(C_{ij\g k})$ of the present paper. Supposition 9.7.12 of [GL16] is wrong: really, for $m\ge2i$ we have $\# Irr(X(m,i))=2^{i-1}$ for $i\le5$, but this is a chance coincidence, and there is nothing common between the ordered partitions of $i$ and elements of $Irr(X(m,i))$ as it was conjectured in [GL16].
\medskip
The simplest irreducible components are the following. First, $C_{111}(m)$ is $X(m,1)$ itself, it is an irreducible hypersurface in $\n P^m$ of degree $m-1$. Further, $C_{ii1}(m)$ is the principal irreducible component of [GL16], 9.7.3.
\medskip
{\bf A2.1. Method of computation.} Entries of Table A2.2 were obtained by computer-assistant computations. We used the Wolfram Mathematica computer system for the cases $i\le4$, and the Sage computer algebra system for $i\ge5$. We consider a random affine space $\Cal Y$ defined over $\n Q$ of dimension $i$ in $\n P^m$,
hence $\Cal Y\cap X(m,i)$ is 0-dimensional. We denote by $D_\Cal Y(m,0), \dots, D_\Cal Y(m,i-1)$ the restrictions of polynomials $D(m,0), \dots, D(m,i-1)$ to $\Cal Y$. Finding of a 0-dimensional variety $\Cal Y\cap X(m,i)$ is the same as solving a system $D_\Cal Y(m,0)=0, \dots, D_\Cal Y(m,i-1)=0$ of $i$ polynomials in $i$ unknowns (coordinates of $\Cal Y$). This is made by elimination of all coordinates of $\Cal Y$ except one of them, using the theory of a multivariate resultant (see [CLO], Chapter 3, Section 4).
\medskip
Let us introduce notations. We denote by $R_x(P,Q)$ the ordinary resultant of polynomials $P(x)$, $Q(x)$ in $x$; it is a number. Now, let $P_1,\dots, P_\vk$ be polynomials in $\vk$ variables $x_1,\dots, x_\vk$ of degrees $a_1,\dots, a_\vk$ respectively. Their multivariate resultant $MR(P_1,\dots, P_\vk)=MR_{x_2,\dots, x_\vk}(P_1,\dots, P_\vk)$ is a polynomial in $x_1$ of degree $a_1\cdot...\cdot a_\vk$ (for generic $P_1,\dots, P_\vk$). Its roots are $x_1$-coordinates of the $a_1\cdot...\cdot a_\vk$ solutions to the system $\{P_1=0, \dots, P_\vk=0\}$ counting with multiplicities (the degree of $MR(P_1,\dots, P_\vk)$ can be less than $a_1\cdot...\cdot a_\vk$ if some solutions to $\{P_1=0, \dots, P_\vk=0\}$ are at infinity).
\medskip
Now we consider a projection $\g x: \Cal Y\cap X(m,i)\to \n A^1$ defined over $\n Q$, here $\n A^1$ is an affine line corresponding to the only coordinate which is not eliminated. Let $\g x(\Cal Y\cap X(m,i))=\sum_\vk \g e_\vk t_\vk$ (equality of divisors), where $t_\vk$ are affine coordinates of points of the support of $\g x(\Cal Y\cap X(m,i))$, and $\g e_\vk$
are their multiplicities. We have: the polynomial $\Cal P=\prod_\vk(x-t_\vk)^{\g e_\vk}\in \n Q[x]$ whose roots are $t_\vk$ with multiplicities $\g e_\vk$ is $MR(D_\Cal Y(m,0), \dots, D_\Cal Y(m,i-1))$.
\medskip
The relation between $\Cal P$ and the degrees and multiplicities of the irreducible components of $X(m,i)$ is the following. Let $Z\subset X(m,i)$ be an irreducible component of degree $d$ and multiplicity $\mu$. Hence, the contribution of $\Cal Y\cap Z$ to $\sum_\vk \g e_\vk t_\vk$ is $\mu (t_{\al_1}+...+t_{\al_d})$ where $t_{\al_1},\dots,t_{\al_d}$ is a subset of $t_1,\dots,t_\vk$. Results of Section 5.15 show that for $i\le6$ all irreducible components of $X(m,i)$, except $C_{618a}$, $C_{618b}$, are defined over $\n Q$. Since $\Cal Y$ is a random $i$-dimensional linear space, we get (for a generic $\Cal Y$) that $t_{\al_1},\dots,t_{\al_d}$ form a set of $\n Q$-conjugate numbers. Hence, after finding $\Cal P$, we calculate its factorization in $\n Q[x]$: $$\Cal P=\prod_{\vk=1}^\g n {\Cal P_\vk}^{\mu_\vk}, \ \ \Cal P_\vk \in \n Q[x] \hbox { are irreducible} \eqno{(A2.1.1)}$$

According the above, we get (for almost all $\Cal Y$; except the case of $C_{618a}$, $C_{618b}$): the quantity of irreducible components of $X(m,i)$ is $\g n$, and the degree and multiplicity of the $\vk$-th irreducible component are $\deg(\Cal P_\vk)$ and $\mu_\vk$ respectively.
\medskip
For the exceptional case of $C_{618a}$, $C_{618b}$ we have that there exists one polynomial $\Cal P_*$ from (A2.1.1) for both these components; its degree is twice the degree of $C_{618*}$ ($*=a$ or $b$), and its $\mu_*$ is the multiplicity of $C_{618*}$.
\medskip
Let us indicate the beginning of the theory of the multivariate resultant. We consider a simple case of 3 polynomials $P_1(x,y,z)$, $P_2(x,y,z)$, $P_3(x,y,z)$ in 3 variables, of degrees $a_1, \ a_2, \ a_3$ respectively. Hence, $R_y(P_1,P_2)$ is a polynomial in $x, \ z$ of degree $a_1a_2$, and $R_y(P_1,P_3)$ is a polynomial in $x, \ z$ of degree $a_1a_3$. Now we consider $R_z(R_y(P_1,P_2), R_y(P_1,P_3))$; it is a polynomial in $x$ of degree $a_1^2a_2a_3$, i.e. it is the $MR_{y,z}(P_1,P_2,P_3)$ times a "parasite factor" of degree $a_1$. The simplest way to eliminate this "parasite factor" is to use a formula
$$MR_{y,z}(P_1,P_2,P_3)=GCD (\ R_z[R_y(P_1,P_2), R_y(P_1,P_3)]; \ R_z[R_y(P_1,P_2), R_y(P_2,P_3)]\ )\eqno{(A2.1.2)}$$
For $i=3$ this formula can be used directly for computing, but for $i\ge4$ it is necessary to eliminate the "parasite factors" manually at some intermediate steps of computation (if not the computation runs out of time and memory). This explains why we cannot represent the computations in a form of one program, but we have to mix computer calculations and hand operations.
\medskip
Finally, the Jordan form corresponding to any irreducible component can be easily got, if we know coordinates $(a_0:...:a_m)$ of points of $\Cal Y\cap X(m,i)$ corresponding to $\vk$-th factor of (A2.1.1). Really, for any $\g f$ we calculate the rank of $\g M(m)(a_0,...,a_m)^\g f$ and we use (5.14.2). Entries of Table A2.2, cases $i\le 4$, and $i=5$, $j=1$ were obtained by this method.
\medskip
For $i=5$, $j\ge2$ computation over $\n Q$ takes too much time and memory, hence we choose a large prime $\g q$ (we use values $\g q=31013$ and 16661) and we make all computations of $MR(D(m,0), \dots, D(m,i-1))$ modulo $\g q$.
Polynomials that are irreducible over $\n Q$ become, most likely, reducible over $\n F_{\g q}$, but it is very few likely that an irreducible over $\n Q$ polynomial becomes non-square-free after reduction modulo ${\g q}$ (the probability of this event is $\sim \frac{1}{\g q}$, i.e. very low). Results of Sections 5.13, 5.14 show that for $i=5$ there are no elements of $Irr(5)$ having the same multiplicity and Jordan form. So, we calculate $\bar\Cal P:=\overline{MR(D(m,0), \dots, D(m,i-1))}$ (bar means reduction modulo $\g q$) and we calculate its factorization $\bar\Cal P=\prod_{\vk=1}^{\tilde\g n} {\tilde\Cal P_\vk}^{\tilde\mu_\vk}$ in $\n F_{\g q}[x]$. If two factors ${\tilde \Cal P_\vk}^{\tilde\mu_\vk}$, ${\tilde \Cal P_{\vk'}}^{\tilde\mu_{\vk'}}$ have the same multiplicity (i.e. $\tilde\mu_\vk=\tilde\mu_{\vk'}$) and the same Jordan form, we think that they come from one polynomial over $\n Q$ in (A2.1.1). Hence, we get (with high probability) the degrees and multiplicities of irreducible components of $X(m,i)$ for this case.
\medskip
This method was used for $i=5$, $j\ge2$ and $i=6$, $j=1,2$. We get that in all cases the irreducible varieties and their numerical characteristics are exactly the same as they must be according conjectures of the present paper. The same result is obtained by two independent conjectural methods, hence we conclude that both methods are correct.
\medskip
Results for $i=6$, $j\ge3$ were not obtained by computer calculations, they were obtained using results of Section 5, hence they cannot be used for justification of Conjectures 5.11.10, 5.13.1. But these results are in concordance with the formula (2.8) for $i=6$, $j\ge3$, so they justify Conjecture 5.10 for these $i$, $j$.
\medskip
{\bf Table A2.2.} Irreducible components of $X(m,i)$, case $i\le 6$.
\medskip
\settabs 10 \columns
\+ Compo- &&Degree &&Multi- &Jordan form&&Forest&\cr
\+nent&$j$&&&plicity&(partition)\cr
\medskip
\+$i=1$&\cr
\medskip
\+$C_{111}$&1&$m-1$&&1&1=1&&$\circ$&\cr
\medskip
\+$i=2$\cr
\medskip
\+$C_{211}$&1&$2(m-2)$&&1&2=2&&$\circ-\circ$&\cr
\medskip
\+$C_{221}$&2&$\binom{m-2}{2}$&&2&2=1+1&&$\circ$ $\circ$&\cr
\medskip
\+$i=3$\cr
\medskip
\+$C_{311}$&1&$4(m-3)$&&1&3=3&&$\circ-\circ-\circ$&\cr
\medskip
\+$C_{312}$&1&$m-3$&&2&3=2+1&&$^{\circ<^\circ_\circ}$\cr
\medskip
\+$C_{321}$&2&$4\binom{m-3}{2}$&&3&3=2+1&&$\circ-\circ$ \ \ $\circ$&\cr
\medskip
\+$C_{331}$&3&$\binom{m-3}{3}$&&6&3=1+1+1&&$\circ$ $\circ$ $\circ$&\cr
\medskip
\+$i=4$\cr
\medskip
\+$C_{411}$&1&$8(m-4)$&&1&4=4&&$\circ-\circ-\circ-\circ$&\cr
\medskip
\+$C_{412}$&1&$2(m-4)$&&2&4=3+1&&$^{\circ-\circ<^\circ_\circ}$\cr
\medskip
\+$C_{413}$&1&$4(m-4)$&&3&4=3+1&&$^{\circ<^{\circ-\circ}_\circ}$\cr
\medskip
\+$C_{421}$&2&$8\binom{m-4}{2}$&&4&4=3+1&&$\circ-\circ-\circ$ \ \ $\circ$&\cr
\medskip
\+$C_{422}$&2&$4\binom{m-4}{2}$&&6&4=2+2&&$\circ-\circ$ \ \ $\circ-\circ$ \cr
\medskip
\+$C_{423}$&2&$2\binom{m-4}{2}$&&8&4=2+1+1&&$^{\circ<^\circ_\circ\ \ \circ}$\cr
\medskip
\+$C_{431}$&3&$6\binom{m-4}{3}$&&12&4=2+1+1&&$\circ-\circ$ \ \ $\circ$ \ \ $\circ$ &\cr
\medskip
\+$C_{441}$&4&$\binom{m-4}{4}$&&24&4=1+1+1+1&&$\circ$ $\circ$ $\circ$ $\circ$ &\cr
\medskip
\+$i=5$&\cr
\medskip
\+$C_{511}$&1&$16(m-5)$&&1&5=5&&$\circ-\circ-\circ-\circ-\circ$&&\cr
\medskip
\+$C_{512}$&1&$4(m-5)$&&2&5=4+1 &&$^{\circ-\circ-\circ<^\circ_\circ}$&&\cr
\medskip
\+$C_{513}$&1&$8(m-5)$&&3&5=4+1 &&$^{\circ-\circ<^{\circ-\circ}_\circ}$&&\cr
\medskip
\+$C_{514}$&1&$8(m-5)$&&4&5=4+1 &&$^{\circ<^{\circ-\circ-\circ}_\circ}$&&\cr
\medskip
\+ Compo- &&Degree &&Multi- &Jordan form&&Forest&\cr
\+nent&$j$&&&plicity&(partition)\cr
\medskip
\+$C_{515}$&1&$4(m-5)$&&6&5=3+2 &&$^{\circ<^{\circ-\circ}_{\circ-\circ}}$&&\cr
\medskip
\+$C_{516}$&1&$2(m-5)$&&8&5=3+1+1 &&$^{\circ<^{\circ<^\circ_\circ}_{\circ}}$&&\cr
\medskip
\+$C_{521}$&2&$16\binom{m-5}{2}$&&5&5=4+1&&$\circ-\circ-\circ-\circ$ \ \ $\circ$&\cr
\medskip
\+$C_{522}$&2&$16\binom{m-5}{2}$&&10&5=3+2&&$\circ-\circ-\circ$ \ \ $\circ-\circ$\cr
\medskip
\+$C_{523}$&2&$4\binom{m-5}{2}$&&10&5=3+1+1&&$^{\circ-\circ<^\circ_\circ\ \ \circ}$\cr
\medskip
\+$C_{524}$&2&$8\binom{m-5}{2}$&&15&5=3+1+1&&$^{\circ<^{\circ-\circ}_\circ\ \ \circ}$\cr
\medskip
\+$C_{525}$&2&$4\binom{m-5}{2}$&&20&5=2+2+1&&$^{\circ<^{\circ}_\circ\ \ \circ-\circ}$\cr
\medskip
\+$C_{531}$  &3&$12\binom{m-5}{3}$&&20&5=3+1+1&&$\circ-\circ-\circ$ \ \ $\circ$ \ \ $\circ$&\cr
\medskip
\+$C_{532}$  &3&$12\binom{m-5}{3}$&&30&5=2+2+1&&$\circ-\circ$ \ \ $\circ-\circ$ \ \ $\circ$&\cr
\medskip
\+$C_{533}$  &3&$3\binom{m-5}{3}$&&40&5=2+1+1+1&&$^{\circ<^{\circ}_\circ\ \ \circ \ \ \circ}$&\cr
\medskip
\+$C_{541}$&4&$8\binom{m-5}{4}$&&60&5=2+1+1+1&&$\circ-\circ$ \ \ $\circ$ \ \ $\circ$ \ \ $\circ$&\cr
\medskip
\+$C_{551}$&5&$\binom{m-5}{5}$&&120&5=1+...+1&&$\circ$ \ \ $\circ$ \ \ $\circ$ \ \ $\circ$ \ \ $\circ$&\cr
\medskip
$i=6$
\medskip
\+$C_{611}$&1&$32(m-6)$&&1&6=6&&$\circ-\circ-\circ-\circ-\circ-\circ$&\cr
\medskip
\+$C_{612}$&1&$8(m-6)$&&2&6=5+1 &&$^{\circ-\circ-\circ-\circ<^\circ_\circ}$&&\cr
\medskip
\+$C_{613}$&1&$16(m-6)$&&3&6=5+1 &&$^{\circ-\circ-\circ<^{\circ-\circ}_\circ}$&&\cr
\medskip
\+$C_{614}$&1&$16(m-6)$&&4&6=5+1 &&$^{\circ-\circ<^{\circ-\circ-\circ}_\circ}$&&\cr
\medskip
\+$C_{615}$&1&$16(m-6)$&&5&6=5+1 &&$^{\circ<^{\circ-\circ-\circ-\circ}_\circ}$&\cr
\medskip
\+$C_{616}$&1&$8(m-6)$&&6&6=4+2 &&$^{\circ-\circ<^{\circ-\circ}_{\circ-\circ}}$&&\cr
\medskip
\+$C_{617}$&1&$4(m-6)$&&8&6=4+1+1 &&$^{\circ-\circ<^{\circ<^\circ_\circ}_{\circ}}$&&\cr
\medskip
\+$C_{618a}$&1&$8(m-6)$&&10&6=4+2 &&${\overset{0}\to{\circ}<^{\overset{0}\to{\circ}-\overset{0}\to{\circ}- \overset{0}\to{\circ}}_{\underset{1}\to{\circ}-\underset{1}\to{\circ}}}$\cr
\medskip
\+$C_{618b}$&1&$8(m-6)$&&10&6=4+2 &&${\overset{0}\to{\circ}<^{\overset{0}\to{\circ}-\overset{0}\to{\circ}- \overset{0}\to{\circ}}_{\underset{1}\to{\circ}-\underset{3}\to{\circ}}}$\cr
\medskip
Remark: $C_{618a}$, $C_{618b}$ have the same tree but different weights. These weights are indicated (in additive form), they belong to different $\g G(T)\rtimes\Aut(T)$-orbits. $C_{618a}(m)$, $C_{618b}(m)$ are defined over $\n Q[\sqrt{-1}]$, they are $\n Q$-conjugate.
\medskip
\+$C_{619}$&1&$4(m-6)$&&10&6=4+1+1 &&$^{\circ<^{\circ-\circ<^\circ_\circ}_{\circ}}$\cr
\medskip
\+$C_{6,1,10}$&1&$8(m-6)$&&15&6=4+1+1 &&$^{\circ<^{\circ<^{\circ-\circ}_\circ}_{\circ}}$\cr
\medskip
\+$C_{6,1,11}$&1&$4(m-6)$&&20&6=3+2+1 &&$^{\circ<^{\circ<^\circ_\circ}_{\circ-\circ}}$\cr
\medskip
\+ Compo- &&Degree &&Multi- &Jordan form&&Forest&\cr
\+nent&$j$&&&plicity&(partition)\cr
\medskip
\+$C_{621}$&2&$32\binom{m-6}{2}$&&6&6=5+1&&$\circ-\circ-\circ-\circ-\circ$ \ \ $\circ$&\cr
\medskip
\+$C_{622}$&2&$8\binom{m-6}{2}$&&12 &6=4+1+1 &&$^{\circ-\circ-\circ<^\circ_\circ\ \ \circ}$&&\cr
\medskip
\+$C_{623}$&2&$32\binom{m-6}{2}$&&15 &6=4+2&&$\circ-\circ-\circ-\circ$ \ \ $\circ-\circ$&\cr
\medskip
\+$C_{624}$&2&$16\binom{m-6}{2}$&&18 &6=4+1+1 &&$^{\circ-\circ<^{\circ-\circ}_\circ\ \ \circ}$&&\cr
\medskip
\+$C_{625}$&2&$16\binom{m-6}{2}$&&20 &6=3+3&&$\circ-\circ-\circ$ \ \ $\circ-\circ-\circ$&\cr
\medskip
\+$C_{626}$&2&$16\binom{m-6}{2}$&&24 &6=4+1+1 &&$^{\circ<^{\circ-\circ-\circ}_\circ\ \ \circ}$&&\cr
\medskip
\+$C_{627}$&2&$8\binom{m-6}{2}$&&30 &6=3+2+1 &&$^{\circ-\circ<^\circ_\circ\ \ \circ-\circ}$&\cr
\medskip
\+$C_{628}$&2&$8\binom{m-6}{2}$&&36 &6=3+2+1 &&$^{\circ<^{\circ-\circ}_{\circ-\circ}\ \ \circ}$&&\cr
\medskip
\+$C_{629}$&2&$8\binom{m-6}{2}$&&40 &6=3+2+1 &&$^{\circ<^\circ_\circ\ \ \circ-\circ-\circ}$&\cr
\medskip
\+$C_{6,2,10}$&2&$16\binom{m-6}{2}$&&45 &6=3+2+1 &&$^{\circ<^{\circ-\circ}_\circ\ \ \circ-\circ}$&&\cr
\medskip
\+$C_{6,2,11}$&2&$4\binom{m-6}{2}$&&48 &6=3+1+1+1&&$^{\circ<^{\circ<^\circ_\circ}_{\circ}\ \ \circ}$&\cr
\medskip
\+$C_{6,2,12}$&2&$\binom{m-6}{2}$&&80 &6=2+2+1+1&&$^{\circ<^{\circ}_{\circ}\ \ \circ<^{\circ}_{\circ}}$&\cr
\medskip
The below data are not a result of computer calculation. They are obtained by application of the results of Section 5.
\medskip
\+$C_{631}$&3&$24\binom{m-6}{3}$&&30 &6=4+1+1 &&$\circ-\circ-\circ-\circ$ \ \ $\circ$ \ \ $\circ$&\cr
\medskip
\+$C_{632}$&3&$48\binom{m-6}{3}$&&60 &6=3+2+1 &&$\circ-\circ-\circ$ \ \ $\circ-\circ$ \ \ $\circ$&\cr
\medskip
\+$C_{633}$&3&$8\binom{m-6}{3}$&&90 &6=2+2+2 &&$\circ-\circ$ \ \ $\circ-\circ$ \ \ $\circ-\circ$&\cr
\medskip
\+$C_{634}$&3&$6\binom{m-6}{3}$&&60 &6=3+1+1+1 &&$^{\circ-\circ<^\circ_\circ\ \ \circ\ \ \circ}$&\cr
\medskip
\+$C_{635}$&3&$12\binom{m-6}{3}$&&90 &6=3+1+1+1 &&$^{\circ<^{\circ-\circ}_\circ\ \ \circ\ \ \circ}$&\cr
\medskip
\+$C_{636}$&3&$12\binom{m-6}{3}$&&120 &6=2+2+1+1 &&$^{\circ<^\circ_\circ\ \ \circ-\circ\ \ \circ}$&\cr
\medskip
\+$C_{641}$&4&$16\binom{m-6}{4}$&&120 &6=3+1+1+1 &&$\circ-\circ-\circ$ \ \ $\circ$ \ \ $\circ$ \ \ $\circ$&\cr
\medskip
\+$C_{642}$&4&$24\binom{m-6}{4}$&&180 &6=2+2+1+1 &&$\circ-\circ$ \ \ $\circ-\circ$ \ \ $\circ$ \ \ $\circ$ &\cr
\medskip
\+$C_{643}$&4&$4\binom{m-6}{4}$&&240 &6=2+1+...+1 &&$^{\circ<^{\circ}_{\circ}\ \ \circ\ \ \circ\ \ \circ}$&\cr
\medskip
\+$C_{651}$&5&$10\binom{m-6}{5}$&&360 &6=2+1+...+1 &&$\circ-\circ$ \ \ $\circ$ \ \ $\circ$ \ \ $\circ$ \ \ $\circ$&\cr
\medskip
\+$C_{661}$&6&$\binom{m-6}{6}$&&720 &6=1+...+1 &&$\circ$ \ \ $\circ$ \ \ $\circ$ \ \ $\circ$ \ \ $\circ$ \ \ $\circ$&\cr
\medskip
\newpage
{\bf A2.3. Components described in (6.3).} These are components having $j=1$ and the Jordan form $i=d_1+d_2$. We give here their list from the above table A2.2, and we indicate their $\g n$, $d_1$, $d_2$. Notation of the component is from A2.2.
\medskip
\+ Component &&$\g n$ &$d_1$&$d_2$&&Component &&$\g n$ &$d_1$&$d_2$\cr
\medskip
\+$C_{312}$&&0&2&1&&$C_{612}$&&3&5&1\cr
\medskip
\+$C_{412}$&&1&3&1&&$C_{613}$&&2&5&1\cr
\medskip
\+$C_{413}$&&0&3&1&&$C_{614}$&&1&5&1\cr
\medskip
\+$C_{512}$&&2&4&1&&$C_{615}$&&0&5&1\cr
\medskip
\+$C_{513}$&&1&4&1&&$C_{616}$&&1&4&2\cr
\medskip
\+$C_{514}$&&0&4&1&&$C_{618a}$&&0&4&2\cr
\medskip
\+$C_{515}$&&0&3&2&&$C_{618b}$&&0&4&2\cr
\medskip
\medskip

{\bf A3. Components having $j$ near $i$.} Let us fix $\vk$ and consider sets $Irr(i,j)$ for $j=i-\vk$, as $i\to\infty$. These sets are stable for $i\ge 2\vk$. For $\vk\le3$ they are described below. We see that (2.8) holds for these $i, \ j$, hence we get evidence of truth of Conjecture 5.10 for these cases.
\medskip

\+ Compo- &$\ \ j$&Degree &&Multi- &Jordan form&&&Forest&\cr
\+nent&&&&plicity&(partition)\cr

\medskip
\+$C_{i,i-3,1}$&$i-3$&$8(i-3)\binom{m-i}{i-3}$&&$(i!)/24$ &$i=4+1+...+1$&&&$^{\circ-\circ -\circ -\circ \ \ \circ \ ... \ \circ}$&&$i\ge4$\cr\medskip
\+$C_{i,i-3,2}$&$i-3$&$8(i-3)(i-4)\cdot$&&$(i!)/12$ &$i=3+2+1+...+1$&&&$^{\circ-\circ -\circ \ \ \circ-\circ  \ \circ \  ... \ \circ}$&&\cr
\+&&$\cdot\binom{m-i}{i-3}$&&&&&&&&$i\ge5$\cr
\medskip
\+$C_{i,i-3,3}$&$i-3$&$8\binom{i-3}{3}\binom{m-i}{i-3}$&&$(i!)/8$ &$i=2+2+2+1+...$&&&$^{\circ-\circ \ \ \circ-\circ \ \ \circ-\circ  \ \circ \ ... \ \circ}$&&\cr
\+&&&&&&$+...+1$&&&&$i\ge6$\cr
\medskip
\+$C_{i,i-3,4}$&$i-3$&$2(i-3)\binom{m-i}{i-3}$&&$(i!)/12$ &$i=3+1+...+1$&&&$^{\circ-\circ<^\circ_\circ\ \ \circ \ ... \ \circ}$&&$i\ge4$\cr
\medskip
\+$C_{i,i-3,5}$&$i-3$&$4(i-3)\binom{m-i}{i-3}$&&$(i!)/8$ &$i=3+1+...+1$&&&$^{\circ<^{\circ-\circ}_\circ \ \ \circ \ ... \ \circ}$&&$i\ge4$\cr\medskip
\+$C_{i,i-3,6}$&$i-3$&$2(i-3)(i-4)\cdot$&&$(i!)/6$ &$i=2+2+1+...+1$&&&$^{\circ<^\circ_\circ \ \ \circ-\circ  \ \ \circ \ ... \ \circ}$&&$i\ge5$\cr
\+&&$\cdot\binom{m-i}{i-3}$\cr
\medskip
\+$C_{i,i-2,1}$&$i-2$&$4(i-2)\binom{m-i}{i-2}$&&$(i!)/6$ &$i=3+1+...+1$&&&$^{\circ-\circ -\circ \ \ \circ \ ... \ \circ}$&&$i\ge3$\cr
\medskip
\+$C_{i,i-2,2}$&$i-2$&$2(i-2)(i-3)\cdot$&&$(i!)/4$ &$i=2+2+1+...+1$&&&$^{\circ-\circ \ \ \circ-\circ \ \ \circ \ ... \ \circ}$&&$i\ge4$\cr
\+&&$\cdot\binom{m-i}{i-2}$\cr
\medskip
\+$C_{i,i-2,3}$&$i-2$&$(i-2)\binom{m-i}{i-2}$&&$(i!)/3$ &$i=2+1+...+1$&&&$^{\circ<^\circ_\circ \ \ \circ \ ... \ \circ}$&&$i\ge3$\cr
\medskip
\medskip
\+$C_{i,i-1,1}$&$i-1$&$2(i-1)\binom{m-i}{i-1}$&&$(i!)/2$ &$i=2+1+...+1$&&&$\circ-\circ \ \ \circ \ ... \ \circ$&&&$i\ge2$\cr
\medskip
\+$C_{ii1}$&$i$&$\binom{m-i}{i}$&&$i!$ &$i=1+...+1$&&&$\circ \ ... \ \circ$&&&$i\ge1$\cr
\medskip
Component $C_{ii1}$ in [GL16] is called the principal component.
\medskip
\newpage
{\bf Table A4.} Quantities of irreducible components of $X(m,6)$ having given $j$ and given Jordan form.
\medskip
\+&&$j=$&1&&&2&3&4&5&6\cr
\medskip
\noindent
Jordan form
\medskip
\+6&&&1\cr
\medskip
\+5+1&&&4&&&1\cr
\medskip
\+4+2&&&3 (over $\n Q[\sqrt{-1}]$)&&&1\cr
\medskip
\+&&&2 (over $\n Q$)\cr
\medskip
\+3+3&&&&&&1\cr
\medskip
\+4+1+1&&&3&&&3&1\cr
\medskip
\+3+2+1&&&1&&&4&1\cr
\medskip
\+3+1+1+1&&&&&&1&2&1\cr
\medskip
\+2+2+2&&&&&&&1\cr
\medskip
\+2+2+1+1&&&&&&1&1&1\cr
\medskip
\+2+1+...+1&&&&&&&&1&1\cr
\medskip
\+1+...+1&&&&&&&&&&1\cr
\medskip
\+{\bf Table A5.} Known (conjectural) results for any $i$, and some $j=1,2$.&&&&&&&&& \cr
\medskip
\noindent
(Here numeration of $C_{i1*}$ differs from numeration in Table A2.2).
\medskip
\+ Component &&Degree &&Multi- &Jordan form&&&&Valid for\cr
\+&$j$&&&plicity&(partition)\cr
\medskip
\+$C_{i11}$&1&$2^{i-1}(m-i)$&&1&$i=i$&&&& any $i$\cr
\medskip
\+$C_{i12}$&1&$2^{i-3}(m-i)$&&2&$i=(i-1)+1$ &&&&$i\ge3$\cr
\medskip
\+$C_{i13}$&1&$2^{i-2}(m-i)$&&3&$i=(i-1)+1$ &&&&$i\ge4$\cr
\medskip
\+$C_{i14}$&1&$2^{i-2}(m-i)$&&4&$i=(i-1)+1$ &&&&$i\ge5$\cr
\medskip
\+$C_{i15}$&1&$2^{i-2}(m-i)$&&5&$i=(i-1)+1$ &&&&$i\ge6$\cr
\medskip
\+$\vdots$ etc., until&&$\vdots$&&$\vdots$&$\vdots$&&\cr
\medskip
\+$C_{i,1,i-1}$&1&$2^{i-2}(m-i)$&&$i-1$&$i=(i-1)+1$ &&&&$i\ge4$\cr
\medskip
$\dots \dots \dots \dots \dots \dots \dots \dots \dots $
\medskip
\+$C_{i1*_1}$&1&$2^{i-3}(m-i)$&&6&$i=(i-2)+2$  &&&&$i\ge5$\cr
\medskip
\+$C_{i1*_2}$&1&$2^{i-4}(m-i)$&&8&$i=(i-2)+1+1$ &&&&$i\ge5$\cr
\medskip
\+$C_{i1*_3}$&1&$2^{i-2}(m-i)$&&10&$i=(i-2)+2$ &&&&$i\ge6$\cr
\medskip
\newpage
\+ Component &&Degree &&Multi- &Jordan form&&&&Valid for\cr
\+&$j$&&&plicity&(partition)\cr
\medskip
\+$C_{i1*_4}$&1&$2^{i-4}(m-i)$&&10&$i=(i-2)+1+1$ &&&&$i\ge6$\cr
\medskip
\+$C_{i1*_5}$&1&$2^{i-3}(m-i)$&&15&$i=(i-2)+1+1$ &&&&$i\ge6$\cr
\medskip
\+$C_{i1*_6}$&1&$2^{i-4}(m-i)$&&20&$i=(i-3)+2+1$ &&&&$i\ge6$\cr
\medskip
\+$C_{i1*_7}$&1&$2^{i-4}(m-i)$&&12&&&&&$i\ge7$\cr
\medskip
\+$C_{i1*_8}$&1&$2^{i-2}(m-i)$&&15&&&&&$i\ge7$\cr
\medskip
\+$C_{i1*_9}$&1&$2^{i-3}(m-i)$&&18&&&&&$i\ge7$\cr
\medskip
\+$C_{i1*_{10}}$&1&$2^{i-3}(m-i)$&&20 &&&&&$i\ge7$\cr
\medskip
\+$C_{i1*_{11}}$&1&$2^{i-3}(m-i)$&&24 &&&&&$i\ge7$\cr
\medskip
\+$C_{i1*_{12}}$&1&$2^{i-4}(m-i)$&&30 &&&&&$i\ge7$\cr
\medskip
\+$C_{i1*_{13}}$&1&$2^{i-4}(m-i)$&&36 &&&&&$i\ge7$\cr
\medskip
\+$C_{i1*_{14}}$&1&$2^{i-4}(m-i)$&&40 &&&&&$i\ge7$\cr
\medskip
\+$C_{i1*_{15}}$&1&$2^{i-3}(m-i)$&&45 &&&&&$i\ge7$\cr
\medskip
\+$C_{i1*_{16}}$&1&$2^{i-5}(m-i)$&&48 &$i=(i-3)+1+1+1$ &&&&$i\ge7$\cr
\medskip
\+$C_{i1*_{17}}$&1&$2^{i-7}(m-i)$&&80 &$i=(i-4)+2+1+1$ &&&&$i\ge7$\cr
\medskip
$\dots \dots \dots \dots \dots \dots \dots \dots \dots $
\medskip
\+$C_{i21}$&2&$2^{i-1}\binom{m-i}{2}$&&$i$&$i=(i-1)+1$ &&&&$i\ge3$\cr
\medskip
{\bf Table A6.} Two irreducible components for any $i$ and $j=i-4$.
\medskip
\settabs 8 \columns
\+ Compo- &&Degree &&Multi- &Jordan form&&Valid for\cr
\+nent&$j$&&&plicity&(partition)\cr
\medskip
\+$C_{i,i-4,*_1}$&$i-4$&$2(i-4)\binom{m-i}{i-4}$&&$(i!)/15$  &$i=3+1+...+1$&&$i\ge5$\cr
\medskip
\+$C_{i,i-4,*_2}$&$i-4$&$\frac{(i-4)(i-5)}{2}\binom{m-i}{i-4}$&&$(i!)/9$  &$i=2+2+1+...+1$&&\cr
\+&&&&&&&$i\ge6$\cr
\medskip
They are lifts of some components with smaller $i$, see (8.10). $C_{i,i-4,*_1}$ is $C_{516}$ for $i=5$ and $C_{6,2,11}$ for $i=6$. $C_{i,i-4,*_2}$ is $C_{6,2,12}$ for $i=6$.
\medskip
{\bf Table A7. Components having $\g d(C_{ij\g k})=1$. }
\medskip
Let $\vk, \de\ge1$. For $i=(2^\vk-1)\de$, $j=\de$ there exists a component of $X(m,i)$ having $\g d(C_{ij\g k})=1$. They are described in the following table:
\settabs 12 \columns
\medskip
\+ $i$&&Component &&Degree &&Conjugate of the Jordan form&&&\cr
\medskip
\+ $(2^\vk-1)\de$&&$C_{(2^\vk-1)\de,\de,*}$&&$\binom{m-(2^\vk-1)\de}{\de}$&& $(2^\vk-1)\de=2^{\vk-1}\de+2^{\vk-2}\de+...+2\de+\de$\cr
\medskip
Particular case $\de=1$:
\settabs 10 \columns
\medskip
\+ $i$&&Component &&Degree &&Conjugate of the Jordan form&&&\cr
\medskip
\+ $2^\vk-1$&&$C_{2^\vk-1,1,*}$&&$m-(2^\vk-1)$&& $2^\vk-1=2^{\vk-1}+2^{\vk-2}+...+2+1$\cr
\medskip
For $\vk=1,2,3$ this is respectively $C_{111}$, $C_{312}$ of Table A2.2 and $C_{i1*_{17}}$ of Table A5 for $i=7$. See also Example 6.4.
\medskip
\medskip
{\bf References}
\medskip

[A86] Anderson, Greg W. $t$-motives.  Duke Math. J.  53  (1986),  no. 2, 457 -- 502
\medskip
[AT90] Anderson, Greg W.; Thakur, Dinesh S. Tensor powers of the Carlitz module and zeta values. Ann. of Math. (2) 132 (1990), no. 1, 159–191.
\medskip
[A00] Anderson, Greg W. An elementary approach to $L$-functions mod $p$. J. Number theory 80 (2000), no. 2, 291 -- 303.
\medskip
[B02]  B\"ockle, Gebhard. Global $L$-functions over function fields. Math. Ann. 323 (2002), no. 4, 737 -- 795.
\medskip
[B05] B\"ockle, Gebhard. Arithmetic over function fields: a cohomological approach. Number fields and function fields -- two parallel worlds, 1 -- 38, Progr. Math., 239, Birkhäuser Boston, Boston, MA, 2005.
\medskip
[BP] B\"ockle, Gebhard; Pink, Richard. Cohomological theory of crystals over function fields. EMS tracts in Mathematics, 9. European Mathematical Society (EMS), Z\"urich, 2009. viii+187 pp.
\medskip
[B12] B\"ockle, Gebhard. Cohomological theory of crystals over function fields and applications.
In "Arithmetic Geometry in Positive Characteristic". Advanced Courses in Mathematics. CRM Barcelona. Birkh\"auser Verlag, Basel (2012)
\medskip
[CLO] Cox, David A, Little, John, O'Shea, Donal. Using Algebraic Geometry, volume 185 of Graduate Texts in Mathematics 1998 - Springer-Verlag, New York.
\medskip
[ELS] Ehbauer, S., Logachev, D., Sarraff de Nascimento, M\'arcia. On some determinantal equalities related to L-functions of Carlitz modules. Communications in Algebra, 2018, vol. 46, No. 5, p. 2130 -- 2145. https://arxiv.org/pdf/1707.04339.pdf
\medskip
[FP] Fulton, William; Pragacz, Piotr. Schubert varieties and degeneracy loci. LNM 1689.
\medskip
[Ga] Gardeyn, F. t-motives and Galois representations. Thesis (Dr.) - Universiteit Gent (Belgium), ETH Z\"urch (Switzerland). 2001.
\medskip
[Ge] Gekeler, Ernst-Ulrich. The Galois image of twisted Carlitz modules. J. of Number Theory, 163 (2016), 316 -- 330.
\medskip
[G79] Goss, David. $v$-adic zeta functions, $L$-series and measures for function fields. Invent. Math. 55 (1979), 107 -- 116.
\medskip
[G92] Goss, David. $L$-series of $t$-motives and Drinfeld modules. The arithmetic of function fields (Columbus, OH, 1991), 313 -- 402, Ohio State Univ. Math. Res. Inst. Publ., 2, de Gruyter, Berlin, 1992.
\medskip
[G96] Goss, D. Basic structures of function field arithmetic. Springer-Verlag, Berlin, 1996. xiv+422 pp.
\medskip
[GL16] Grishkov, A.; Logachev, D. Resultantal varieties related to zeroes of $L$-functions of Carlitz modules. Finite Fields Appl. 38 (2016), 116-176.

arxiv.org/pdf/1205.2900.pdf
\medskip
[GL20] Grishkov, A.; Logachev, D. Introduction to Anderson t-motives: a survey. https://arxiv.org/pdf/2008.10657.pdf
\medskip
[GL21] Grishkov, A.; Logachev, D.; Zobnin, A. L-functions of Carlitz modules, resultantal varieties and rooted binary trees - I. 38 pages. Accepted for publication by J. of Number Theory. https://doi.org/10.1016/j.jnt.2021.08.013 A short version of the present paper. 
\medskip
[HJ] U. Hartl and A.-K. Juschka, Pink's theory of Hodge structures and the Hodge conjecture
over function fields, in: "t-motives: Hodge structures, transcendence and other motivic aspects", Editors G. B\"ockle, D. Goss, U. Hartl, M. Papanikolas, European Mathematical Society Congress Reports 2020.

https://arxiv.org/pdf/1607.01412.
\medskip
[H] Hurwitz, A. \"Uber die Bedingungen, unter welchen eine Gleichung nur Wurzeln mit negativen reellen Teilen besitzt. Mathematische Annalen, Leipzig, 46 (1895), 273 - 284. 
\medskip
[L] Lafforgue, Vincent. Valeurs sp\'eciales des fonctions $L$ en caract\'eristique $p$.  J. Number Theory  129  (2009),  no. 10, 2600 -- 2634
\medskip
[PPVW] Park, Jennifer; Poonen, Bjorn; Voight, John; Wood, Melanie Matchett. A heuristic for boundedness of ranks of elliptic curves. J. Eur. Math. Soc. (JEMS) 21 (2019), no. 9, 2859 -- 2903
\medskip
[Ta] Taelman, L. Special $L$-values of Drinfeld modules. Ann. of Math. (2) 175 (2012), no. 1, 369 -- 391.
\medskip
[TW] Taguchi, Y.; Wan, D. $L$-functions of $\vf$-sheaves and Drinfeld modules. J. Amer. Math. Soc. 9 (1996), no. 3, 755 -- 781
\medskip
[Th] Thakur, Dinesh S. On characteristic $p$ zeta functions. Compositio Math. 99 (1995), no. 3, 231 -- 247.
\medskip
[W1] Wikipedia, https://en.wikipedia.org/wiki/Partition\_(number\_theory)

\#Conjugate\_and\_self-conjugate\_partitions
\medskip
[W2] Wikipedia, https://en.wikipedia.org/wiki/Wedderburn-Etherington\_number

\enddocument
\medskip
\medskip
\medskip
\medskip
\medskip
\medskip
\medskip
\medskip
\medskip
\medskip

Let $m$, $i$ be two integer parameters, $m>i\ge 1$. Paper [GL16], Section 9 contains a definition of certain projective
varieties $X(m,i)\subset \n P^m(\n C)$ and some conjectures on their irreducible components, their degrees, multiplicities etc. The present paper continues this research.
\medskip
The origin of varieties $X(m,i)$ is the theory of $L$-functions of the twisted $n$-th tensor powers of Carlitz modules over a finite field $\n F_q$ where $q$ is a power of a prime $p$. We can expect that the study of $X(m,i)$ is the first step to solve the below Problem of boundedness of ranks of twists of Carlitz modules (Open question 0.0.3). This problem is a functional field analog of a famous problem of boundedness of ranks of twists of a fixed elliptic curve. See 0.4 below for further steps to solve (0.0.3).
\medskip

\medskip
Therefore, let us repeat some definitions of [GL16] (see (0.1) below for the details). Let $\g C^n$ be the $n$-th tensor power of the Carlitz module $\g C=\g C_q$ over $\n F_q$. Its twists are parameterized by polynomials $P\in \n F_q[\theta]$ (here $\th$ is an independent variable), see (0.1.2) for the definition of a twist. A twist corresponding to such $P$ is denoted by $\g C^n_P$. We fix $m$ --- the degree of $P$, and let
$$P=\sum_{\iota=0}^m a_\iota\theta^\iota \hbox{ where } a_\iota\in \n F_q\eqno{(0.0.1)}$$
The $L$-function $L(\g C^n_P,T)=L(P,T)=L((a_0,\dots, a_m),T)$ is an element of $(\n F_q[t])[T]$ where $t$ is an independent variable\footnotemark \footnotetext{For a general t-motive $M$ we have $L(M,T)\in (\n F_q[t])[[T]]$, but for these particular $M=\g C^n_P$ we have $L(\g C^n_P,T)\in(\n F_q[t])[T]$, see Remark 0.1.7, (1).} (we consider $n$ as a fixed number).
\medskip
{\bf Remark.} There are various versions of $L$-functions of Anderson t-motives. See (0.1.3A), (0.1.4), (0.1.5) for a definition of a version of the $L$-function used in the present paper, and (0.1.7A) for its relations with other types of $L$-functions.
\medskip
{\bf Definition 0.0.2.} The analytic rank of $\g C^n_P$ (or, by abuse of language, of $P$, or of $\{a_*\}=(a_0,\dots, a_m)$) is the multiplicity of the root $T=1$ to the equation $L(\g C^n_P,T)=0$; it is denoted by $r_1(P)=r_1(n,P)$ (here and below the subscript 1 means that we consider the behavior of $L(P,T)$ at $T=1$). We have $r_1(P)\ge0$.
\medskip
{\bf Remark.} We do not know whether this definition is a good analog of the notion of the analytic rank of elliptic curves, or not. Moreover, an analog of the algebraic rank of Anderson t-motives is not known yet, although [L], Proposition 2.1, p. 2604 can be considered as an analog of the strong form of the Birch and Swinnerton-Dyer conjecture for $L(M,T)$ at $T=1$ (here $M$ is an Anderson t-motive or even some its generalization like in [L]).
\medskip
But we have no better definition than the above one. Results of [GL16] show that the behavior of $r_1(P)$, while $P$ varies, resembles the behavior of ranks of twists of a fixed elliptic curve. Namely, for most $P$ we have $r_1(P)=0$; there exists a coset in a group of twists (see (0.1.1c)) such that for $P$ belonging to this coset we have $r_1(P)\ge1$ ([GL16], (4.4), (4.5)); rare jumps of $r_1(P)$ occur. From another side, there is no information on the parity of $r_1(P)$, because for Anderson t-motives there is no functional equation for their $L$-series.
\medskip
There is a general
\medskip
{\bf Open question 0.0.3.} --- Problem of boundedness of ranks of twists of Carlitz modules: Are $r_1(P)$ bounded? ($q$, $n$ are fixed, $P$ varies).
\medskip
{\bf 0.0.4.} There are arguments in favor of both possibilities. Namely, there exists an affine algebraic variety (denoted by $X_1(q,n,m,i)_{fc}$; the subscript $fc$ always means "finite characteristic") which is the set of $\bar \n F_q$-zeroes of some (non-homogeneous) polynomials $D_{**fc}\in\n F_p[a_0,\dots, a_m]$. [GL16], Theorem 3.3 shows that the set of $(a_0,\dots,a_m)\in A^{m+1}(\n F_q)$ such that $r_1(P)\ge i$, is $X_1(q,n,m,i)_{fc}(\n F_q)$ --- the set of $\n F_q$-points of $X_1(q,n,m,i)_{fc}$, i.e. the set of $\n F_q$-zeroes of $D_{**fc}$. See [GL16], between (6.2) and (6.3) for the definition of $D_{**fc}$ (in [GL16] polynomials $D_{**fc}$ are denoted by $D_{**}$); two subscripts ** mean that $D_{**fc}$ depend on 2 integer parameters (see (0.1.13) for details).
\medskip
From one side, there are too many $D_{**fc}$'s. If all polynomials $D_{**fc}$ were independent then the answer to  0.0.3 would be yes (see [GL16], Proposition 6.6).
\medskip
From another side, it seems that $D_{**fc}$ are dependent. For example, if they were independent then for $q=3$, $n=1$ we would have $r_1(P)\le3$ ([GL16], Proposition 6.6). Really, there are examples of $P$ such that $r_1(P)=6$ ([GL16], Table 6.11). Further, direct count of the quantity of points of $X_1(3,1,m,i)_{fc}(\n F_3)$ for small $m,i$ ([GL16], (6.8); Table 6.12 and below remarks) shows that most likely codim $X_1(3,1,m,i)_{fc}$ is less than the quantity of $D_{**fc}$, i.e. $D_{**fc}$ are dependent. Finally, $D_{**fc}$ are linear combinations of $H_{**fc}$ (see below) which are highly dependent.
\medskip
{\bf Remark.} It is not known what is even a conjectural answer to the analog of (0.0.3) in characteristic 0: are ranks of all twists of a fixed elliptic curve over $\n Q$, or, more widely, of all elliptic curves over $\n Q$, bounded or not? Earlier there was a conjecture that the ranks are not bounded. A recent paper [P] gives evidence that the ranks of all elliptic curves over $\n Q$ are bounded.
\medskip
{\bf 0.0.5.} It turns out that it is much easier to consider the behavior (as $P$ varies) of $L(P, T)$ not at $T=1$ but at $T=\infty$. The order of pole of $L(P, T)$ at $T=\infty$ (roughly speaking; see (0.1.11) for the exact definition) is called the analytic rank of $P$ (or of $a_*$) at infinity, it is denoted by $r_\infty(P)$ (or $r_\infty(a_0,\dots,a_m)$ ). Like for the case of rank at 1, the set of $(a_0,\dots,a_m)\in A^{m+1}(\n F_q)$ such that $r_\infty(P)\ge i$ is the set of $\n F_q$-points of an algebraic variety (its projectivization is denoted by $X_\infty(q,n,m,i)_{fc}$) over $\bar \n F_q$, i.e. the set of zeroes of some homogeneous polynomials $H_{**fc}\in\n F_p[a_0,\dots, a_m]$.
\medskip
It is easy to show that $D_{**fc}$ are linear combinations of $H_{**fc}$ with integer coefficients, so we can expect that study of $X_\infty(q,n,m,i)_{fc}$ can shed light to study of $X_1(q,n,m,i)_{fc}$ and hence to a solution of the Question 0.0.3.
\medskip
It is easy to prove ([GL16], Theorem 8.6) that $H_{**fc}$ are highly dependent, i.e. $X_\infty(q,n,m,i)_{fc}$ are of low codimension.
\medskip
{\bf 0.0.6.} There exists a natural lift of $H_{**fc}$ from $\n F_p\subset \n F_q$ to $\n Z$, see (0.2.0) --- (0.2.2). We denote these lifts by $H_{**}$. Hence we can consider the varieties of zeroes of $H_{**}$ in $\n P^m(\n C)$. They  are denoted by $X_\infty(q,n,m,i)$. These varieties are defined over $\n Q$, because coefficients of $H_{**}$ are in $\n Z$. It turns out that for this case of characteristic 0 the simplest possible value of $q$ is 2: although for $q=2$ there are no twists of $\g C^n$ (see (0.1.3)), the theory of $X_\infty(2,n,m,i)$ --- unlike the theory of $X_\infty(2,n,m,i)_{fc}$ over $\n F_2$ --- is non-trivial. Moreover, it turns out that for $q=2$ the varieties $X_\infty(2,n,m,i)$ (conjecturally) do not depend on $n$, see Conjecture 0.2.4, Remark 0.2.5 (1).
\medskip
{\bf 0.0.7.} We can consider formally polynomials $H_{**}$ for $n=0$: although $\g C^0$ does not exist, the polynomials $H_{**}$, and hence varieties $X_\infty(2,0,m,i)$ do exist, and (again conjecturally) $X_\infty(2,0,m,i)$ coincides with $X_\infty(2,n,m,i)$ for all $n>0$. Finally, the case $n=0$ is especially simple, because for this case the set of $H_{**}$ depends on 1 parameter and not on 2 parameters that simplifies their study. See (0.2) for a more detailed discussion.
\medskip
By definition, $X(m,i)$ are exactly $X_\infty(2,0,m,i)$. The above arguments show the importance of these objects.
\medskip
The subsequent subsections of the Introduction contain a more detailed exposition of the theory of $L$-functions of twisted Carlitz modules, the contents of the present paper and the possibilities of further research.
\enddocument